\documentclass[a4paper,10pt]{article}
\usepackage{amssymb}
\newtheorem{Th}{Theorem}
\newtheorem{Prop}{Proposition}

\newtheorem{Lm}{Lemma}
\newtheorem{Lma}{Lemma}[section]
\newtheorem{Dfi}{Definition}
\newtheorem{Rm}{Remark}

\newcommand{\be}{\begin{equation}}
\newcommand{\ee}{\end{equation}}
\newcommand{\R}{\mathbb{R}}
\newcommand{\N}{\mathbb{N}}
\newcommand{\C}{\mathbb{C}}
\newcommand{\Z}{\mathbb{Z}}

\newcommand\res{\mathop{\hbox{\vrule height 7pt width .5pt depth 0pt
\vrule height .5pt width 6pt depth 0pt}}\nolimits}

\newcommand{\reset}{\setcounter{equation}{0}\setcounter{Th}{0}\setcounter{Prop}{0}\setcounter{Co}{0}
\setcounter{Lm}{0}\setcounter{Rm}{0}}
\def\La{\Lambda}
\def\La{\Lambda}
\def\ti{\tilde}
\def\lf{\left}
\def\rg{\right}

\def\al{\alpha}
\def\la{\lambda}

\def\ep{\varepsilon}
\def\ds{\displaystyle}
\def\ov{\overline}

\def\om{\omega}
\def\p{\partial}

\def\si{\sigma}
\def\res{\mathop{\hbox{\vrule height 7pt width .5pt 
depth 0pt\vrule height .5pt width 6pt depth 0pt}}\nolimits}
\begin{document}
\title{Variational Principles for immersed Surfaces with $L^2$-bounded Second Fundamental Form.}
\author{Tristan Rivi\`ere\footnote{Department of Mathematics, ETH Zentrum,
CH-8093 Z\"urich, Switzerland.}}
%\date{ }
\maketitle

{\bf Abstract :} {\it In this work we present new fundamental tools for studying the variations of the Willmore functional
of immersed surfaces into ${\R}^m$.
This approach gives for instance a new proof of the existence of a Willmore minimizing embedding of an arbitrary closed surface  in
arbitrary codimension. We explain how the same approach  can solve constraint minimization problems for the Willmore functional.
We show in particular that, for a given closed surface and a given conformal class for this surface,  there is an immersion in ${\R}^m$, away possibly from
isolated branched points, which minimizes the Willmore energy among all possible Lipschitz  immersions in ${\R}^m$ having an $L^2-$bounded second
fundamental form and realizing this conformal class. This branched immersion is either a smooth Conformal Willmore branched immersion or an isothermic branched immersion.
We show that branched points do not exist whenever the minimal  energy in the conformal class is less than $8\pi$ and that these immersions extend to smooth conformal Willmore embeddings or global isothermic embeddings of the surface in that case.
Finally, as a by-product of our analysis, we establish that inside a compact subspace of the moduli space the following holds :  weak limit of Palais Smale { Willmore} sequences are { Conformal Willmore}, that weak limits of Palais Smale sequences
of { Conformal Willmore} are either { Conformal Willmore} or { Global Isothermic} and finally
we observe also that weakly converging Palais Smale sequences of { Global Isothermic Immersions} are { Global Isothermic}.
The analysis developped along the paper - in particular these last results - opens the door to the possibility of constructing new critical saddle points of the Willmore functional without or with constraints using min max methods.}

\medskip

\noindent{\bf Math. Class.} 30C70, 58E15, 58E30, 49Q10, 53A30, 35R01, 35J35, 35J48, 35J50.

\section{Introduction}

The goal of the present paper is to present a suitable framework to proceed to the calculus 
of variation of the Willmore functional for immersions.

\medskip

Let $\vec{\Phi}$ be a smooth immersion (rank$\,d\vec{\Phi}$ is equal to two at every point) from a closed oriented smooth 2-manifold $\Sigma$ into an euclidian space
${\R}^m$. The {\it first fundamental form} $g_{\vec{\Phi}}$ defined by this immersion on $\Sigma$ is the pull-back by ${\vec{\Phi}}$ of the metric induced by the restriction of the canonical metric of ${\R}^m$, $g_{{\R}^m}$, to the tangent planes $\vec{\Phi}_\ast T\Sigma$ of the immersed surface :
\[
g_{\vec{\Phi}}:=\vec{\Phi}^\ast g_{{\R}^m}\quad\quad.
\]
If there is no ambiguity we simply write $g$ instead of $g_{\vec{\Phi}}$. We shall denote by $\vec{n}_{\vec{\Phi}}$ the {\it Gauss map} of the immersion $\vec{\Phi}$ which to a point $p\in\Sigma$ assigns the oriented orthonormal $(m-2)-$plane to the tangent plane $\vec{\Phi}_\ast T_p\Sigma$ of the immersed surface at  $\vec{\Phi}(p)$. $\vec{n}_{\vec{\Phi}}$
will be seen as  a map into the Grassmanian $\ti{G}_{m-2}({\R}^m)$ (or equivalently $\ti{G}_2({\R}^m)$) of oriented $m-2-$planes (resp. $2-$planes) of ${\R}^m$. $\vec{n}_{\vec{\Phi}}$ is also a map into the unit simple \footnote{$\vec{n}_{\vec{\Phi}}\wedge\vec{n}_{\vec{\Phi}}=0$.} 2-vectors 
in ${\R}^m$ : $\vec{n}_{\vec{\Phi}}\in\wedge^2{\R}^m$. We also denote by $\pi_{\vec{n}_{\vec{\Phi}}}$ the orthonormal projections of vectors in ${\R}^m$ onto
the $m-2$-plane given by $\vec{n}_{\vec{\Phi}}$. With these notations the {\it second fundamental form}
\[
\forall X,Y\in T_p\Sigma\quad\quad\quad \vec{\mathbb I}_p(X,Y):=\pi_{\vec{n}_{\vec{\Phi}}}d^2\vec{\Phi}(X,Y)
\]\footnote{In order to define $d^2\vec{\Phi}(X,Y)$ one has to extend locally the vector $X$ or $Y$ by a vectorfield but it is not difficult to
check that $\pi_{\vec{n}_{\vec{\Phi}}}d^2\vec{\Phi}(X,Y)$ is independent of this extension.}
The {\it mean curvature vector} of the immersion at $p$ is given by
\[
\vec{H}:=\frac{1}{2}\,tr_g(\vec{\mathbb I})=\frac{1}{2}\,\lf[\vec{\mathbb I}(\ep_1,\ep_1)+\vec{\mathbb I}(\ep_2,\ep_2)\rg]\quad,
\]
where $(\ep_1,\ep_2)$ is an orthonormal basis of $T_p\Sigma$ for the metric $g_{\vec{\Phi}}$. 

\medskip
In the present paper we are mainly interested with the Lagrangian given by 
the $L^2$ norm of the second fundamental form :
\[
E(\vec{\Phi}):=\int_{\Sigma} |\vec{\mathbb I}|^2_g\ dvol_g\quad,
\]
where $dvol_g$ is the volume form induced by the metric $g_{\vec{\Phi}}$. An elementary computation gives 
\[
E(\vec{\Phi}):=\int_{\Sigma} |\vec{\mathbb I}|^2_g\ dvol_g=\int_{\Sigma} |d\vec{n}_{\vec{\Phi}}|^2_g\ dvol_g\quad.
\]
This energy $E$ can be hence seen as being the {\it Dirichlet Energy} of the Gauss map $\vec{n}_{\vec{\Phi}}$ with respect to the induced
metric $g_{\vec{\Phi}}$.
The Gauss Bonnet theorem implies that 
\be
\label{0.1}
E(\vec{\Phi}):=\int_{\Sigma} |\vec{\mathbb I}|^2_g\ dvol_g=4\ \int_{\Sigma}|\vec{H}|^2\ dvol_g-4\pi\ \chi(\Sigma)\quad,
\ee
where $\chi(\Sigma)$ is the {\it Euler characteristic} of the surface $\Sigma$.
The energy
\[
W(\vec{\Phi}):=\int_\Sigma|\vec{H}|^2\ dvol_g\quad,
\]
is the so called {\it Willmore energy} and has been extensively studied since the early 20th century due in one hand to it's  rich mathematical
signification but also to it's importance in other area of science (in general relativity, mechanics, biology, optics...etc). 
Probably the main property it satisfies which makes this lagrangian so universal is the conformal invariance : For any conformal diffeomorphism $\vec{U}$
of ${\R}^m$ one has (see \cite{Bla})
\be
\label{0.2}
W(\vec{\Phi})=W(\vec{U}\circ\vec{\Phi})\quad.
\ee
For a fixed surface $\Sigma$ ,
because of (\ref{0.1}) studying the variations of the $L^2-$norm of the second fundamental form or the variations of Willmore energy is identical.

\medskip

Since the lower bound to $W(\vec{\Phi})$ among all possible immersions of closed surfaces is non zero and equal to $4\pi$ (see for instance \cite{Wi})
it is natural to look at the existence of optimal immersions which minimize $W$ for a given surface $\Sigma$. When $\Sigma$ is a sphere it is well known that $W(\vec{\Phi})$ achieves it's minimal value $4\pi$
for the standard unit $S^2$ in ${\R}^3\subset{\R}^m$ and only for this submanifold. When $\Sigma$ is a genus 1 surfaces the existence of a smooth immersion into ${\R}^m$
minimizing $W$ was established by L.Simon in \cite{Si}. It has been conjectured by T.J.Willmore that the minimizing configuration should be achieved 
by the torus of revolution in ${\R}^3$ obtained by rotating around the $z-$axis the vertical circle included in the $Oxz$ plane, of center $(\sqrt{2},0,0)$ and
radius $1$ and the minimal energy would then be $2\pi^2$. This conjecture is still open at this stage. The existence result of L.Simon has been extended
to surfaces of arbitrary genus by M.Bauer and E.Kuwert in \cite{BK} : they proved that for an arbitrary given closed oriented surface $\Sigma$ there is an immersion into ${\R}^m$ that minimizes the Willmore energy among all smooth immersions of that surface.
The result of Bauer and Kuwert was using the result of L.Simon whose proof is quite involved. One of the characteristic of this proof is not to work directly with the immersions
$\vec{\Phi}$ but mostly instead with its image $\vec{\Phi}(\Sigma)$ \footnote{ Or more precisely with the rectifiable current $\vec{\Phi}_\ast[\Sigma]$ : the push forward by $\vec{\Phi}$ of the integration current over $\Sigma$} ..

\medskip

In the present paper we first present a new proof of Simon-Bauer-Kuwert's result that will based on the analysis of the immersions themselves . This new proof 
 will be ''transposable'' to the minimization of the Willmore functional under various constraints, as it arises in several applications,
( prescribed effective volume, prescribed conformal class...etc ) without having to change the main lines of the proof. Moreover the arguments and tools that we will develop
in this work should be sufficiently generic in order to generate new critical points of the Willmore functional under various constraints by applying fundamental
principle of the calculus of variation such as the mountain pass lemma...etc as we shall present it in a forthcoming work \cite{Ri3}.

\medskip 

The first difficulty encountered while working with immersions $\vec{\Phi}$ instead of working with their image $\vec{\Phi}(\Sigma)$ is the huge invariance group
of the functional : the space of diffeomorphisms of $\Sigma$, Diff$(\Sigma)$. Taking for instance a minimizing sequence $\vec{\Phi}_k$ 
of the Willmore functional (without or with constraints) one can always compose $\vec{\Phi}_k$ with diffeomorphisms that makes the sequence degenerate completely and not reaching 
an immersion at all ! There is then a ''choice of gauge'' to be made. By pulling back the standard
metric $g_{{\R}^m}$ of ${\R}^m$ onto $\Sigma$, $\vec{\Phi}$ defines then a metric $g_{\vec{\Phi}}$ on $\Sigma$ and hence a conformal structure on $\Sigma$.
There exists then a constant scalar curvature metric $h$ on $\Sigma$ and a conformal diffeomorphism $\Psi$ from $(\Sigma,h)$ into $(\Sigma,g_{\vec{\Phi}})$
such that the immersion $\vec{\Phi}\circ\Psi$ is conformal. The space of constant scalar curvature metrics on $\Sigma$ identifies (modulo dilations) to the space
of conformal structures on $\Sigma$ and hence is finite dimensional see \cite{Jo}. We have then broken the ''gauge degeneracy'' by replacing $\vec{\Phi}$ by $\vec{\Phi}\circ\Psi$ which satisfies the
{\it Coulomb gauge condition} :
\[
\lf\{
\begin{array}{l}
\ds div(\vec{e}_1,\nabla\vec{e}_2)=0\quad\quad\mbox{ where }\quad\quad\vec{e}_j:=e^{-\la}\ \p_{x_j}(\vec{\Phi}\circ\Psi)\quad,\\[5mm]
\ds\mbox{and }\quad\quad e^{\la}:=|\p_{x_1}(\vec{\Phi}\circ\Psi)|=|\p_{x_2}(\vec{\Phi}\circ\Psi)|
\end{array}
\rg.
\]
and the operator $div$ and $\nabla$ are the standard operators : $divX=\p_{x_1}X_1+\p_{x_2}X_2$ and $\nabla\cdot:=(\p_{x_1}\cdot,\p_{x_2}\cdot)$
taken in arbitrary complex coordinates $z=x_1+ix_2$ with respect to the conformal structure given by $(\Sigma,h)$. 

At this stage however the possible perplexity of the reader regarding this choice of gauge is totally justified because what we have gained by fixing the gauge that way
is not clear at all at this stage.
Indeed, looking again at a minimizing sequence $\vec{\Phi}_k$ of the Willmore functional (without or with constraints) and composing by $\Psi_k$ in order
to have a conformal immersion, first one does not have a-priori a control of the conformal class defined by $\Phi_k$ : $(\Sigma,h_k)$ may degenerate to the boundary
of the moduli space. A first task in our proof is to exclude this eventuality. More seriously, as $k$ goes to infinity we have a-priori no control at all of the conformal
factor $e^{\la_k}$ that could either go to $+\infty$ or $0$ at some points and then we would be out of the class of immersions at the limit.
This cannot be excluded easily. The problem is that the control of the $L^2$ norm of the second fundamental form does not provide a global pointwise control of the conformal factor $e^\la$ - counter-examples are easy to manufacture. This is just critical : an $L^{2+\ep}-$control of the second fundamental form would have done it. However, below a certain threshold this control exists.
This phenomenon has been discovered in a series of works by T.Toro \cite{To1} \cite{To2}, S.M\"uller-V.Sverak \cite{MS} and F.H\'elein \cite{Hel}.
Precisely one has
\begin{Th}
\label{th-0.1} {\bf [Control of local isothermal coordinates]}
Let $\vec{\Phi}$ be a conformal immersion of the disc $D^2$ such that
\be
\label{0.3}
\int_{D^2}|\nabla\vec{n}_{\vec{\Phi}}|^2<8\pi/3\quad\quad\mbox{ and }\quad\quad M(\vec{\Phi}_\ast[D^2])=\int_{D_2}e^{2\la}\ dx_1\,dx_2<+\infty
\ee
where $M(\vec{\Phi}_\ast[D^2])$ is the mass \footnote{$M(\vec{\Phi}_\ast[D^2]):=\sup\{\int_{D^2}\vec{\Phi}^\ast\om\ ; \ \|\om\|_\infty\le1\}$} of the current $\vec{\Phi}_\ast[D^2]$. Then for any $0<\rho<1$ there exists a constant $C_\rho$ independent of $\vec{\Phi}$ such that
\be
\label{0.3a}
\sup_{p\in D^2_\rho}e^\la(p)\le C_\rho\ \lf[ M(\vec{\Phi}_\ast[D^2])\rg]^{1/2}\ \exp\lf(\int_{D^2}|\nabla\vec{n}_{\vec{\Phi}}|^2\rg)\quad.
\ee
Moreover, for two given distinct points $p_1$ and $p_2$ in the interior of $D^2$ and again for $0<\rho<1$ there exists a constant $C>0$
independent of $\vec{\Phi}$ such that
\be
\label{0.3b}
\begin{array}{rl}
\ds\|\la\|_{L^\infty(D^2_\rho)}\le & \ds C_\rho\ \int_{D^2}|\nabla\vec{n}_{\vec{\Phi}}|^2+C_\rho\ \log|\vec{\Phi}(p_1)-\vec{\Phi}(p_2)|^{-1}\\[5mm]
 &\ds+C_\rho\ \log M(\vec{\Phi}_\ast[D^2])\quad.
 \end{array}
\ee
\hfill$\Box$
\end{Th}
\begin{Rm}
\label{rm-0.1}
The existence of two distinct points $p_1$ and $p_2$ such that, in the minimization process, $|\vec{\Phi}_k(p_1)-\vec{\Phi}_k(p_2)|$ is not 
converging to zero - ie the maintenance of the {\bf non-collapsing condition} - will be obtained -see the 3-point normalization lemma~\ref{lm-III.1} - by the composition with an ad-hoc M\"obius transformation of ${\R}^m$ which does not
affect the Willmore energy - see (\ref{0.2}) - and hence the minimizing nature of the sequence.\hfill$\Box$
\end{Rm}

This Theorem is only implicit in the above mentioned works and therefore we give a proof of it in section III. The main ingredients for proving theorem~\ref{th-0.1} are the following.
First, under the energy assumption (\ref{0.3}), one constructs a  controlled energy orthonormal frame ''lifting'' the Gauss map $\vec{n}_{\vec{\Phi}}$. Precisely one has
\be
\label{0.4}
\begin{array}{l}
\ds\forall \vec{n}\in W^{1,2}(D^2,\ti{G}_2({\R}^m))\quad\mbox{ s.t. }\quad\int_{D^2}|\nabla\vec{n}|^2<8\pi/3\quad,\\[5mm]
\exists\, (\vec{f}_1,\vec{f}_2)\in (S^{m-1})^2\quad\mbox{s.t. }\quad \vec{f}_1\cdot\vec{f}_2=0\quad\quad,\quad\quad\vec{n}_{\vec{\Phi}}=\vec{f}_1\wedge\vec{f}_2\quad,\\[5mm]
\ds\quad{and }\quad\quad\int_{D^2}|\nabla\vec{e}_1|^2+|\nabla\vec{e}_2|^2\le C\ \int_{D^2}|\nabla\vec{n}|^2\quad.
\end{array}
\ee
where $C>0$ is some universall constant. Next one observes that the logarithm of the conformal factor ,$\la$, satisfies the following equation\footnote{In fact equation (\ref{0.5}) is satisfied by any such a lifting and in particular by $(\vec{e}_1,\vec{e}_2)=e^{-\la}(\p_{x_1}\vec{\Phi},\p_{x_2}\vec{\Phi})$. We have
indeed  $\p_{x_1}\vec{e}_1\cdot\p_{x_2}\vec{e}_2-\p_{x_2}\vec{e}_1\cdot\p_{x_1}\vec{e}_2=\p_{x_1}\vec{f}_1\cdot\p_{x_2}\vec{f}_2-\p_{x_2}\vec{f}_1\cdot\p_{x_1}\vec{f}_2$ however
in the present equation the advantage of $(\vec{f}_1,\vec{f}_2)$ over $(\vec{e}_1,\vec{e}_2)$ comes from the fact that we control it's $W^{1,2}-$energy by
the $L^2$ norm of the second fundamental form which is not the case a-priori for $(\vec{e}_1,\vec{e}_2)$.}
:
\be
\label{0.5}
-\Delta\la=\p_{x_1}\vec{f}_1\cdot\p_{x_2}\vec{f}_2-\p_{x_2}\vec{f}_1\cdot\p_{x_1}\vec{f}_2\quad.
\ee

The main estimate to exploit this equation and its Jacobian structure in the r.hs. is given by the following theorem of H.Wente which has shown to play a central role in 
the analysis of 2-dimensional conformally invariant problems (see \cite{Ri1}).
\begin{Th}
\label{th-0.2}\cite{We} {\bf [Regularity by compensation]} Let $a$ and $b$ be two functions in $W^{1,2}(D^2)$ and $\varphi$ be the solution to the following
equation
\be
\label{0.6}
\lf\{
\begin{array}{l}
\ds-\Delta \varphi=\p_{x_1}a\,\p_{x_2}b-\p_{x_2}a\,\p_{x_1}b\quad\quad\mbox{ in }D^2\\[5mm]
\ds\varphi=0\quad\quad\quad\mbox{ on }\p D^2
\end{array}
\rg.
\ee
then the following estimates holds
\be
\label{0.7}
\|\varphi\|_{L^\infty(D^2)}+\|\nabla\varphi\|_{L^2(D^2)}\le C\ \|\nabla a\|_{L^2(D^2)}\ \|\nabla b\|_{L^2(D^2)}\quad.
\ee
where $C$ is some universal constant.\hfill$\Box$
\end{Th}
$\la$ is then the sum of solutions to equations of the form (\ref{0.6}) and some harmonic rest. Combining this decomposition, Wente theorem,
and Harnack inequalities for the harmonic rest are the main arguments in the proof of theorem~\ref{th-0.1}. The operation of finding a lifting
of the Gauss map $\vec{n}_{\vec{\Phi}}$, $(\vec{f}_1,\vec{f}_2)$, whose energy is controlled by the $L^2$ norm of the second fundamental form (as in (\ref{0.4})
is the main limitation for having to restrict to energy below $8\pi/3$. This construction was proved in \cite{Hel} lemma 5.1.4. It is not difficult to construct a counter-example
to the statement (\ref{0.4}) when $8\pi/3$ is replaced by any number strictly larger than $8\pi$. F.H\'elein conjectured however that $8\pi/3$ should be replaced by $8\pi$
and this would make the statement (\ref{0.4}) necessarily optimal. 

\medskip

The previous discussion explains how, while minimizing the Willmore functional (without or with constraints), the problem of the indeterminacy due to the huge invariance
group Diff$(D^2)$ is locally solved  and, as a consequence of theorem~\ref{th-0.1}, beside possibly at most isolated points where the second fundamental form is 
concentrating at least $8\pi/3$ energy, the conformal factor cannot degenerate in the minimization process. However the assumption of having a smooth immersion
(beside these isolated points) at the limit could be lost a-priori since locally the $L^2-$norm of the second fundamental form cannot control more than
the $L^\infty$ norm of the conformal factor\footnote{unless our limit is known to satisfy some special equation of course but we will come to that later.}. 
It is then necessary, following a classical approach in calculus of variations, to ''embed'' the problem in a weak class of immersions. 

\medskip

Let $g_0$ be a reference smooth metric on $\Sigma$. One defines the Sobolev spaces $W^{k,p}(\Sigma,{\R}^m)$ of measurable maps from $\Sigma$ into 
${\R}^m$ in the following way
\[
W^{k,p}(\Sigma,{\R}^m)=\lf\{f\ \mbox{ meas. } {\Sigma}\rightarrow {\R}^m\mbox{ s.t. }\sum_{l=0}^k\int_{\Sigma}|\nabla^l f|_{g_0}^p\ dvol_{g_0}<+\infty\rg\}
\]
Since $\Sigma$ is assumed to be compact it is not difficult to see that this space is independent of the choice we have made of $g_0$.

\medskip

First we need to have a weak first fundamental form that is we need  $\vec{\Phi}^\ast g_{{\R}^m}$ to define an $L^\infty$ metric with a bounded inverse.
The last requirement is satisfied if we assume that $\vec{\Phi}$ is in $W^{1,\infty}(\Sigma)$ and if $d\vec{\Phi}$ has maximal rank 2 at every point with some uniform quantitative control
of ''how far'' $d\vec{\Phi}$ is from being degenerate : there exists $c_0>0$ s.t.
\be
\label{I.1}
|d\vec{\Phi}\wedge d\vec{\Phi}|_{g_0}\ge c_0>0 \quad.
\ee
where $d\vec{\Phi}\wedge d\vec{\Phi}$ is a 2-form on $\Sigma$ taking values into 2-vectors from ${\R}^m$ and given in local coordinates
by $2\,\p_x\vec{\Phi}\wedge\p_y\vec{\Phi}\ dx\wedge dy$.         
The condition (\ref{I.1}) is again independent of the choice of the metric $g_0$ .
For a Lipschitz immersion satisfying (\ref{I.1}) we can define the Gauss map as being the following measurable map in $L^\infty(\Sigma)$
\[
\vec{n}_{\vec{\Phi}}:=\star\frac{\p_x\vec{\Phi}\wedge\p_y\vec{\Phi}}{|\p_x\vec{\Phi}\wedge\p_y\vec{\Phi}|}\quad.
\] 
We then introduce the space ${\mathcal E}_{\Sigma}$ of Lipschitz immersions of $\Sigma$ 
with bounded second fundamental form as follows :
\[
\mathcal{E}_\Sigma:=\lf\{
\begin{array}{c}
\ds\vec{\Phi}\in W^{1,\infty}(\Sigma,{\R}^m)\quad\mbox{ s.t. } \vec{\Phi} \mbox{ satisfies }(\ref{I.1})\mbox{ for some }c_0\\[5mm]
\ds\mbox{ and }\quad\int_{\Sigma}|d\vec{n}|_g^2\ dvol_g<+\infty
\end{array}
\rg\}\quad .
\]
Any Lipschitz immersion $\vec{\Phi}$  in ${\mathcal E}_\Sigma$ defines a smooth conformal structure on $\Sigma$. This comes again from the works of T.Toro \cite{To1} \cite{To2}, S.M\"uller-V.Sverak \cite{MS} and F.H\'elein \cite{Hel} :

\begin{Th}
\label{th-I.1} (\cite{To1},\cite{To2},\cite{MS},\cite{Hel} theorem 5.1.1) {\bf [Existence of local isothermal coordinates]}
Let $\vec{\Phi}\in {\mathcal E}_{D^2}$ satisfying
\be
\label{I.6}
\int_{D^2}|d\vec{n}_{\vec{\Phi}}|^2_g\ dvol_g<\frac{8\pi}{3}\quad,
\ee
then there exists a bilipschitz homeomorphism of the disk $\zeta\in W^{1,\infty}(D^2,D^2)$ such that 
\be
\label{I.7}
\lf\{
\begin{array}{l}
\ds |\p_x(\vec{\Phi}\circ\zeta)|^2-|\p_y(\vec{\Phi}\circ\zeta)|^2=0\quad\quad\mbox{in }D^2\\[5mm]
\ds \p_x(\vec{\Phi}\circ\zeta)\cdot\p_y(\vec{\Phi}\circ\zeta)=0\quad.
\end{array}
\rg.
\ee
\hfill$\Box$
\end{Th}
Hence for any lipschitz immersion $\vec{\Phi}$ in ${\mathcal E}_\Sigma$ one takes a finite covering of $\Sigma$ by disks $(U_j)$ such that
$\int_{U_j}|d\vec{n}_{\vec{\Phi}}|^2_{g}\ dvol_g<8\pi/3$, one gets bilipschitz homeomorphisms $\zeta_j$ for which $\vec{\Phi}\circ\zeta_j$ satisfies
(\ref{I.7}) and hence the transition functions $\zeta^{-1}_k\circ\zeta_j$ are holomorphic. $(U_j,\zeta_j)$ defines then a smooth conformal structure on $\Sigma$.
Let $h$ be a constant scalar curvature associated to this conformal structure and the smooth diffeomorphism $\Psi$ of $\Sigma$ such that the maps $\zeta_j^{-1}\circ\Psi$
are conformal from $(\Sigma,h)$ into $D^2$, then we have that $\vec{\Phi}\circ\Psi$ is a conformal $W^{1,\infty}\cap W^{2,2}$ immersion of $(\Sigma,h)$.
Using theorem~\ref{th-0.1} we can construct local isothermal coordinates for $(\Sigma,g_{\vec{\Phi}\circ\Psi})$ \underbar{with estimates} -i.e. satisfying  (\ref{0.3a})
and (\ref{0.3b}) - and work
with maps in ${\mathcal E}_\Sigma$ like with smooth embeddings.

\medskip

The next main difficulty encountered while working with the immersion $\vec{\Phi}$  instead of it's image in the minimization process of Willmore equation
comes from the Euler Lagrange equation as it has been written originally in the early 20th century in the works of W.Blaschke\cite{Bla}, G.Thomsen \cite{Tho} (in codim 1, i.e. $m=3$) and J. Weiner \cite{Wei} (arbitrary $m$).
A smooth immersion $\vec{\Phi}$ is a critical point to Willmore functional :
\[
\forall \vec{\xi}\in C^\infty(\Sigma,{\R}^m)\quad\quad\frac{d}{dt}W(\vec{\Phi}+t\vec{\xi})_{t=0}=0
\]
if and only if $\vec{\Phi}$ satisfy the {\it Willmore equation}
\be\label{0.8} 
\Delta_\perp\vec{H}\,-\,2\,|\vec{H}|^2\vec{H}\,+\,\ti{A}(\vec{H})\:=\:0\:,
\ee
where $\Delta_\perp$ is the negative covariant Laplacian for the connection\footnote{Namely, for every section $\sigma$ of $N_{\vec{\Phi}}(\Sigma)$\,, one has $D_{X}\sigma:=\pi_{\vec{n}_{\vec{\Phi}}}(\sigma_\ast X)$} $D$ in the normal bundle $N_{\vec{\Phi}}(\Sigma)$\, derived from the ambient scalar product in ${\R}^m$  and where $\ti{A}_p(\vec{L})=\sum_{i,j} \vec{\mathbb I}_p(\vec{e}_i,\vec{e}_j)\ \vec{\mathbb I}_p(e_i,e_j)\cdot\vec{L}$ for $\vec{L}\in{\R}^m$.

In \cite{Ri2} we explained why the Euler Lagrange written in the form (\ref{0.8}) seems not compatible from the lagrangian it is coming from in the sense that $W(\vec{\Phi})$ only
controls the $L^2-$norm of the mean curvature whereas in order to give a distributional meaning to the non linearities in the equation like  $|\vec{H}|^2\vec{H}$
one needs more information on the regularity of $\vec{\Phi}$ ($\vec{H}\in L^3$ for instance for this term). \footnote{This is a bit like writting the Euler Lagrange of the Dirichlet energy $\int_{D^2}|\nabla u|^2$ - that is the Laplace equation $\Delta u=0$ -
in the form $\Delta u^2-|\nabla^2u|^2=0$ which requires $u$ to have at least two derivatives in $L^2$ though the lagrangian gives only a-priori a $W^{1,2}-$control !}
One of the main achievement in \cite{Ri2} was to find a new formulation of the Willmore equation as a conservation law which makes sense for immersions in ${\mathcal E}_\Sigma$.
\begin{Th}
\label{th-0.2b} \cite{Ri2} {\bf[The conservative Willmore Equation]}
The Willmore equation (\ref{0.8}) is equivalent to 
\be
\label{0.9}
d\lf(\ast_g\, d\vec{H}-3\ast_g\pi_{\vec{n}_{\vec{\Phi}}}\big(d\vec{H}\big)\rg)\,-\,d\star\lf( d\vec{n}_{\vec{\Phi}}\wedge \vec{H}\rg)\:=\:0\:,
\ee
 where $\ast_g$ is the Hodge operator on $\Sigma$ associated with the induced metric $g_{\vec{\Phi}}$, and $\,\star$ is the usual Hodge operator on forms. \\[1ex]
\noindent
In particular, a conformal immersion $\vec{\Phi}$ from the flat unit-disc $D^2$ into ${\R}^m$ is Willmore if and only if
\be
\label{0.10}
\Delta\vec{H}\,-\,3\ div\big(\pi_{\vec{n}_{\vec{\Phi}}}(\nabla \vec{H})\big)\,+\,div\star\lf(\nabla^\perp\vec{n}_{\vec{\Phi}}\wedge\vec{H}\rg)\:=\:0\:,
\ee
where the operators $\nabla$, $\nabla^\perp$, $\Delta$, and $div$ are understood with respect to the flat metric on $D^2$. Namely, $\nabla=(\p_{x_1},\p_{x_2})$, $\nabla^\perp=(-\p_{x_2},\p_{x_1})$, $\Delta\cdot=\p^2_{x_1^2}\cdot+\p^2_{x_2^2}$, and \,$div\, X=\p_{x_1}X_1+\p_{x_2}X_2$.
\end{Th}
This conservative form of the Willmore equation  and more conservation laws attached to it permits to pass to the limit in local Palais Smale sequences to
the Willmore Lagrangian. The following result is established in \cite{BR}
\begin{Th}
\label{th-0.3} \cite{BR}
{\bf [Convergence of Willmore Palais Smale sequences.]}
Let $\vec{\Phi}_k$ be a sequence of conformal immersions in ${\mathcal E}_{D^2}$. Assume 
\begin{itemize}
\item[i)]
$$\int_{D^2}|\nabla\vec{n}_{\vec{\Phi}_k}|^2<8\pi/3\quad,$$
\item[ii)]
$$\vec{\Phi}_k(D^2)\subset B_1^m(0)$$ 
\item[iii)]
$$
\exists\,p_1,\, p_2\in D^2\quad\mbox{ s.t. }\quad \liminf |\vec{\Phi}_k(p_1)-\vec{\Phi}_k(p_2)|>0\quad.
$$
\end{itemize}
Under these assumptions if  
\be
\label{0.11}
\begin{array}{l}
\Delta\vec{H}_k\,-\,3\ div\big(\pi_{\vec{n}_{\vec{\Phi}_k}}(\nabla \vec{H}_k)\big)\,+\,div\star\lf(\nabla^\perp\vec{n}_{\vec{\Phi}_k}\wedge\vec{H}_k\rg)\\[5mm]
\quad\quad\quad\longrightarrow 0\quad\quad\mbox{ in }(W^{2,2}\cap W^{1,\infty}(D^2))^\ast
\end{array}
\ee
then\footnote{$\vec{H}_k$ denotes the mean-curvature vector of the immersion $\vec{\Phi}_k$.} there exists a subsequence $\vec{\Phi}_{k'}$ converges weakly in $W^{2,2}_{loc}$ to an analytic immersion $\vec{\Phi}_\infty$ of the disc satisfying
the {\it Conformal Willmore equation}
\be
\label{0.12}
\Delta\vec{H}_\infty-3\ div\big(\pi_{\vec{n}_{\vec{\Phi}_\infty}}(\nabla \vec{H}_\infty)\big)+div\star\lf(\nabla^\perp\vec{n}_{\vec{\Phi}_\infty}\wedge\vec{H}_\infty\rg)=\Im\lf[f(z)\, \vec{H}_{0,\infty}\rg]
\ee
where  $f(z)$ is an holomorphic function of the disc $D^2$ and $\vec{H}_{0,\infty}$ is the {\it Weingarten map} of the immersion $\vec{\Phi}_\infty$ given by
\[
\vec{H}_{0,\infty}:= -2^{-1}\ e^{-2\la_\infty}\ \lf[\p_x\vec{n}_{\vec{\Phi}_\infty}\res\p_x\vec{\Phi}_\infty-\p_y\vec{n}_{\vec{\Phi}_\infty}\res\p_y\vec{\Phi}_\infty+2i\ \p_x\vec{n}_{\vec{\Phi}_\infty}\res\p_y\vec{\Phi}_\infty\rg]
\]
where $e^{\la_\infty}=|\p_x\vec{\Phi}_\infty|=|\p_y\vec{\Phi}_\infty|$ and $\res$ is the standard contraction operator in ${\R}^m$ between a multi-vector and a vector.
\hfill $\Box$
\end{Th}
The {\it Conformal Willmore equation} also called {\it Constrained Willmore} is obtained by considering critical points of the Willmore functional among immersions realizing a fixed conformal class and assuming the critical point is not {\it isothermic} - we shall see this notion a bit below - (see \cite{BPP}). $f(z)$ is just the expression in the conformal chart of an holomorphic quadratic
differential $q:=f(z)\ dz\otimes dz$ of the Riemann surface whose conformal structure is generated by $\vec{\Phi}$.

\medskip

The analyticity of the weak limits to local {\it Palais-Smale sequences} described in the theorem~\ref{th-0.3} above was obtained by proving that the {\it Constrained Willmore equation}
(\ref{0.12}) is equivalent to the the existence of $S\in W^{1,2}(D^2,{\R})$ and $\vec{R}\in W^{1,2}(D^2,\wedge^2{\R}^m)$ satisfying
\be
\label{0.13}
\left\{\begin{array}{l}
\ds-\Delta S=-\big(\nabla\star\vec{n}_{\vec{\Phi}}\big)\cdot\nabla^\perp\vec{R}\\[5mm]
\ds-\Delta\vec{R}=(-1)^{m-1}\star\big(\nabla\vec{n}_{\vec{\Phi}}\bullet\nabla^\perp\vec{R}\big)\,+\,(\nabla\star \vec{n}_{\vec{\Phi}})\,\nabla^\perp S\\[5mm]
\ds-\,\Delta\vec{\Phi}\;=\;\nabla\vec{R}\,\bullet\nabla^\perp\vec{\Phi}\;+\,\nabla S\;\nabla^\perp\vec{\Phi}\:.
\end{array}
\right.
\ee
We call this system the {\it Conservative Conformal Willmore System}. Observe that the right-hand-sides of this system is only made of linear combinations of jacobians of functions which
 are at least $W^{1,2}$. Using Wente theorem~\ref{th-0.2}
one easily bootstrap in the equation and obtain the smoothness of $\vec{\Phi}$. This {\it conservative form} of the {\it Conformal Willmore equation} is also the key tool for 
passing to the limit in {\it Palais-Smale sequences} of Willmore Lagrangian (see \cite{BR}).

\medskip

We have then understood how to control isothermal coordinates locally and the local convergence of ''almost Willmore surfaces'' towards analytic conformal Willmore surfaces
Moreover these two operations can be extended to the framework of weak immersions : to the space ${\mathcal E}_\Sigma$ of lipschitz immersions with $L^2-$bounded
second fundamental forms. Our task now is to collect these local procedures in order to be able to deal with the minimization procedure itself.
To that aim we introduce a distance $d$ - or more precisely a family of equivalent distances on ${\mathcal E}_\Sigma$ - for which the space $({\mathcal E}_\Sigma,d)$
will be \underbar{complete}. The details of the construction of this distance are given in section II.  As explained in section III, under the assumption
that there is a minimizing sequence of Willmore functional with conformal class not diverging in the Moduli space of $\Sigma$ we can make use of {\it Ekeland Variational Principle} in order
to produce sequences of immersions that will be Palais Smale - satisfying (\ref{0.11})  -  in the {\it controlled isothermal coordinates} constructed in theorem~\ref{th-0.1}  and that will converge
to an analytic immersion in these charts. A consequence of theorem~\ref{th-0.1} is that an extraction of subsequence is possible in order to cover any compact part of the surface minus finitely many fixed points by {\it controlled isothermal charts}.
Hence we obtain  at the limit an element in ${\mathcal E}_\Sigma$ which is minimizing $W$ and which is analytic away at most from finitely many points. These points are removable
due in one hand to the fact that the total Willmore energy of this immersion has to be  less than $8\pi$ - for minimality reason - and in the other hand to the fact that the with Li-Yau inequality excludes the possibility of having a branched point below $8\pi$ (see lemma~\ref{lm-A.4} and the argument at the end of section III).
Finally we exclude the possibility of the conformal class to degenerate to the ''boundary of the moduli space'' for energetic reasons (this is explained in \cite{Ri1}) and we have given a new  proof
of the following theorem originally due to L.Simon and M.Bauer-E.Kuwert for the space of smooth immersions (see \cite{Si} and \cite{BK}) and that we can extend to the space ${\mathcal E}_\Sigma$ of Lipschitz immersions with $L^2-$bounded second fundamental form.
\begin{Th}
\label{th-0.4}
{\bf [Existence of a minimizer of $W$ in ${\mathcal E}_\Sigma$]}.
Let $m$ be an arbitrary dimension larger than 2, let $\Sigma$ be a smooth  compact orientable surface without boundary. Then there exists a smooth Willmore embedding of $\Sigma$ into ${\R}^m$
minimizing the Willmore energy among all Lipschitz immersions with $L^2-$bounded second fundamental form (i.e. elements in ${\mathcal E}_\Sigma$).  \hfill $\Box$
\end{Th}

In section IV of the paper we explain how to adapt the argument for proving the previous theorem to a minimization problem under constraint.
Before to state the result we recall a definition.

\begin{Dfi}
\label{df-0.5}{\bf [Local isothermic immersions.]}
A $C^2$ immersion of a surface $\Sigma$ is called local isothermic if, away from the umbilic points, the curvature lines define conformal coordinates. \hfill $\Box$
\end{Dfi}
{\it Local isothermic immersions} realizes very particular surfaces that have been studied since the XIXth century. A survey on the classical geometry of {\it local isothermic immersions} as well
as their role in {\it integrable system theory} can be found in \cite{Bu1}, \cite{Bu2} and \cite{Toj}. There are some characterization of {\it local isothermic immersions}
that we recall in the 3 following propositions for which we prove in section V.
\begin{Prop}
\label{pr-0.6}
A $C^2$ immersion $\vec{\Phi}$ of a surface $\Sigma$ is local isothermic if and only if there exists an holomorphic quadratic differential  $q$ (locally in some complex coordinates $q=f(z)\ dz\otimes dz$) of the Riemann
surface $\ti{\Sigma}$, obtained by equiping $\Sigma$ with the complex structure generated by $\vec{\Phi}$ and by removing the umbilic points of the immersion $\vec{\Phi}$, such that
\be
\label{aza}
\lf<q,h_0\rg>_{WP}:=\Im\lf(f(z)\ \vec{H}_0\rg)\ dz\wedge d\ov{z}\equiv 0
\ee
where $\vec{H}_0$ is the conjugate of the {\it Weingarten map} in local coordinates 
$$\vec{H}_0:=-2^{-1}\ e^{-2\la}\ \lf[\p_x\vec{n}_{\vec{\Phi}}\res\p_x\vec{\Phi}-\p_y\vec{n}_{\vec{\Phi}}\res\p_y\vec{\Phi}+2i\ \p_x\vec{n}_{\vec{\Phi}}\res\p_y\vec{\Phi}\rg]$$
and $h_0$ is the Weingarten Operator given locally by 
$$
h_0:=\p_z\vec{n}_{\vec{\Phi}}\res\p_z\vec{\Phi}\ dz\otimes dz=2^{-1}\ e^{2\la}\ \ov{\vec{H}_0}\ dz\otimes dz
$$
 and $\lf<\cdot,\cdot\rg>_{WP}$
 is the Weil-Peterson pointwise product\footnote{ One verifyies easily that the two form $f(z)\ \vec{H}_0\ dz\,d\ov{z}$ is independent of the
local choice of complex coordinates and defines a complex valued 2-form on $\Sigma$ which is the pointwise Hermitian Weil Peterson product.}
\hfill $\Box$
\end{Prop}
Since an element in ${\mathcal E}_\Sigma$ defines a smooth conformal structure and since it defines an $L^2$ Weingarten operator one easily extend
the notion of isothermic immersions to elements in ${\mathcal E}_\Sigma$. There is a last characterisation of isothermic immersions that we recall
and which also permits to define isothermic immersions for an element in ${\mathcal E}_\Sigma$ and which coincide with the previous characterization
also for element in ${\mathcal E}_\Sigma$
\begin{Prop}
\label{pr-0.7}
A Lipschitz immersion $\vec{\Phi}$ in ${\mathcal E}_\Sigma$ is local isothermic if and only if, away from umbilic points, there exists local complex coordinates
for the structure defined by $\vec{\Phi}$ and a Lipschitz map $\vec{L}$ on this coordinate domain such that
\be
\label{bzb}
e^{-2\la}\,\p_z\vec{\Phi}=\p_{\ov{z}}\vec{L}\quad.
\ee
where $e^{2\la}=|\p_x\vec{\Phi}|^2=|\p_y\vec{\Phi}|^2$ and $z=x+iy$. \hfill$\Box$
\end{Prop}
This proposition is also proved in section V.

\medskip

Finally we give a last characterization of local isothermic immersion which we will meet in our proof, this is the following proposition which is also proved in the section V.

\begin{Prop}
\label{pr-0.7b}
A Lipschitz immersion $\vec{\Phi}$ in ${\mathcal E}_\Sigma$ is local isothermic if and only if, away from umbilic points, there exists in local
conformal coordinates a Lipschitz map $\vec{L}$ from these coordinates into ${\R}^m$ such that
\be
\label{0.013}
\lf\{
\begin{array}{l}
\ds\nabla^\perp\vec{L}\cdot\nabla\vec{\Phi}=0\\[5mm]
\ds\nabla^\perp\vec{L}\wedge\nabla\vec{\Phi}=0
\end{array}
\rg.
\ee
\hfill $\Box$
\end{Prop}
\begin{Rm}
\label{rm-0.9}
The Isothermic equation in the form (\ref{0.013}) has to be compared with the Conformal Willmore equation written in \cite{BR} : there exists $\vec{L}$ such that
\be
\label{0.014}
\lf\{
\begin{array}{l}
\nabla^\perp\vec{L}\cdot\nabla\vec{\Phi}=0\\[5mm]
\nabla^\perp\vec{L}\wedge\nabla\vec{\Phi}=2\ (-1)^m\nabla^\perp\lf(\star(\vec{n}_{\vec{\Phi}}\res\vec{H})\rg)\res\nabla\vec{\Phi}\quad .
\end{array}
\rg.
\ee
There is an apparent strong similarity between the two equations, the isothermic equation (\ref{0.013}) ''corresponds'' to the Conformal Willmore equation (\ref{0.014}) when its right-hand-side of (\ref{0.014}) is replaced by zero. There is however a major difference between these two equations.
Equation (\ref{0.014}) is \underbar{elliptic} and this can be seen by showing that $(S,{\vec{R}})$ given by 
\[
\lf\{
\begin{array}{l}
\ds\nabla S:=\nabla\vec{\Phi}\cdot \vec{L}\quad\quad,\\[5mm]
\ds \nabla {\vec{R}}:=\nabla\vec{\Phi}\wedge\vec{L}+2\nabla^\perp\vec{\Phi}\wedge\vec{H}
\end{array}
\rg.
\]
solves the elliptic system (\ref{0.13}) (see \cite{Ri2}), whereas (\ref{0.013}) is \underbar{hyperbolic} (see \cite{Ri1}), which is confirmed by the next remark.
\hfill $\Box$
\end{Rm}

\begin{Rm}
\label{rm-0.8}
 It is an interesting question to ask how regular lipschitz {\it local isothermic immersions} with $L^2-$ bounded fundamental form can be.
They are not necessarily analytic since axially symmetric surfaces are automatically {\it isothermic immersions}. Hence an arbitrary axially symmetric surface
with $L^2-$bounded second fundamental form is isothermic but not necessarily $C^2$ - it is however $C^{1,1/2}$ and it would be interesting either to try to find 
 less regular examples in ${\mathcal E}_\Sigma$ or to prove that {\it isothermic surfaces } are necessarily $C^{1,1/2}$.\hfill $\Box$
\end{Rm}  
Finally we need the following more restrictive notion of isothermic immersions that we call {\it global isothermic}. Because of the lack of regularity of the most
elementary examples of isothermic immersions such as rotationally invariant surfaces the definition of {\it global isothermic immersions} requires  a framework
that includes non $C^2$ immersions but which, however, define a \underbar{smooth} complex structure on $\Sigma$. The framework of Lipschitz immersions
with $L^2-$bounded second fundamental form ${\mathcal E}_\Sigma$ seems to be the most suitable for that and comes also naturally in the minimization
procedure of Willmore surfaces inside a conformal class as we will describe below. 
\begin{Dfi}
\label{df-0.9}{\bf [Global Isothermic Immersions.]}
An immersion $\vec{\Phi}$ in ${\mathcal E}_\Sigma$ is called global isothermic if there exists an holomorphic quadratic differential $q$ of the Riemann surface defined
by $\vec{\Phi}$ such that
\[
\lf<q,h_0\rg>_{WP}\equiv 0\quad.
\]
\hfill $\Box$
\end{Dfi}
A characteristic of {\it global isothermic immersions} is to be the degenerate
points for the map which to an immersion assigns its conformal class (see lemma~\ref{lm-A.5} in section V). This is why it is not so surprising to see them
appearing as singular points in the minimization process under constrained conformal class - and will appear also as singular points in min-max procedures (see \cite{Ri3}). Adapting the method we used to prove theorem~\ref{th-0.4} to the
constrained case we establish the following result which is the main result of the present work.
\begin{Th}
\label{th-0.9}
{\bf [Existence of a minimizer of $W$ in a conformal class]}.
Let $m$ be an arbitrary integer larger than 2, let $\Sigma$ be a smooth  compact orientable surface without boundary and $c$ a conformal class for this surface. Then there exists an immersion, away from possibly isolated branched points, minimizing the Willmore energy
in the sub-space of lipschitz immersions with $L^2-$bounded second fundamental form, ${\mathcal E}_\Sigma$, realizing the conformal class $c$. Such a minimal immersion is either a smooth Conformal Willmore immersion of $\Sigma$ in ${\R}^m$
satisfying 
\[
\Delta_\perp\vec{H}\,-\,2\,|\vec{H}|^2\vec{H}\,+\,\ti{A}(\vec{H})\:=\Im (q,h_0)_{WP}
\]
where $q$ is an holomorphic quadratic differential of $(\Sigma,c)$, $h_0$ the Weingarten Operator  and $(\cdot,\cdot)$ is the pointwise Weil-Peterson product\footnote{In local complex coordinates $(q,h_0)_{WP}=e^{-2\la}\ f(z)\ \vec{H}_0$ where $q=f(z)\ dz\otimes dz$.} or\footnote{The {\it ''or''} is not exclusive, there are isothermic immersions which are Conformal Willmore.} it is a global isothermic immersion. If the minimal Willmore energy in this conformal class is less
than $8\pi$ then there the immersion has no branched points and it  extends to an embedding of $\Sigma$. \hfill $\Box$
\end{Th}
Partial existence results of minimizers of the Willmore energy in a given conformal class for the dimensions $m=3$ and $m=4$ have been announced in \cite{Scm} and \cite{Sct}.

\medskip

As a byproduct of our analysis we observe that in a compact subset of the moduli space of the surface $\Sigma$ the following holds : weak limit of Palais Smale {\bf Willmore} are {\bf Conformal Willmore}, that Weak Limits of Palais Smale sequences
of {\bf Conformal Willmore} are either {\bf Conformal Willmore} or {\bf Global Isothermic} and finally
we observe also that weakly converging Palais Smale sequences of {\bf Global Isothermic Immersions} are {\bf Global Isothermic}.
This notion of global Palais smale will be presented and used in a forthcoming paper \cite{Ri3} to present the Mountain Pass Lemma for Willmore energy in order to produce saddle points
for this lagrangian with or without constraints.

\medskip

Our paper is organized as follows : in section II we define the metric space of Lipschitz conformal immersions with $L^2-$bounded second fundamental forms. In section III we give
a proof of the existence of a minimizer of the Willmore energy for an arbitrary closed surface $\Sigma$ and an arbitrary codimension (i.e. proof of the theorem~\ref{th-0.4}). In section IV we show how the proof in section III
can be adapted to prove the existence of a minimizer of the Willmore energy in a conformal class (i.e. proof of the theorem~\ref{th-0.9}). In section V we present isothermic
immersion  and explain why they are the degenerate points of the conformal class mapping. In the appendix we give the proof of several
lemmas and propositions used in the previous sections.

\section{ The metric space of lipschitz immersions with $L^2-$bounded second fundamental form.}
\subsection{Definitions and notations.}

Let $\Phi\in \mathcal{E}_\Sigma$ and $\Psi$ be a Lipschitz homeomorphism from $D^2$ into $\Sigma$. For a.e. $(x,y)\in D^2$
we denote $H(D(\vec{\Phi}\circ\Psi))$ the Hopf differential of $\vec{\Phi}\circ\Psi$ :
\[
H(\nabla(\vec{\Phi}\circ\Psi)):=\lf[|\p_x(\vec{\Phi}\circ\Psi)|^2-|\p_y(\vec{\Phi}\circ\Psi)|^2\rg]-2\,i\ \p_x(\vec{\Phi}\circ\Psi)\cdot\p_y(\vec{\Phi}\circ\Psi)
\]
Similarly for a metric $g=\sum_{ij=1}^2g_{ij}\ dx_i\otimes dx_j$ on the disc $D^2$ we define
\[
H(g):=\lf[g_{11}-g_{22}\rg]-2\,i\, g_{12}\quad.
\footnote{Observe that for any immersion $\vec{\Phi}$ of $D^2$ into ${\R}^m$ with our notations $H(\nabla\vec{\Phi})=H(\vec{\Phi}^\ast g_{{\R}^m})$.} 
\]
Remark that if $H(\nabla(\vec{\Phi}\circ\Psi))=0$ then, due to the conformal invariance of the Dirichlet energy one has
\be
\label{I.5}
\int_{\Psi(D^2)}|d\vec{n}_{\vec{\Phi}}|^2_g\ dvol_g=\int_{D^2}|\nabla (\vec{n}_{\vec{\Phi}\circ\Psi})|^2(x,y)\ dx\,dy\quad,
\ee
where $\nabla$ is the usual gradient operator on the disk $D^2$ for the flat metric : $\nabla :=(\p_x\cdot,\p_y\cdot)$. 

\medskip

For $\Psi\in W^{1,\infty}(D^2,D^2)$ such that $\log|\nabla\Psi|\in L^\infty(D^2)$ we denote by $Dis(\vec{\Phi}\circ\Psi)(x,y)$ the complex distortion at $(x,y)\in D^2$ given by
\[
Dis(\vec{\Phi}\circ\Psi)(x,y):=\frac{H(\nabla(\vec{\Phi}\circ\Psi))}{|\nabla \vec{\Phi}\circ\Psi|^2}(x,y)\quad.
\]
Similarly also we define for a metric $g=\sum_{ij=1}^2g_{ij}\ dx_i\otimes dx_j$ on the disc $D^2$ the complex distortion of this metric
to be 
\[
Dis(g):=\frac{H(g)}{tr\, g}=\frac{ g_{11}-g_{22}-2\, i\, g_{12}}{g_{11}+g_{22}}\quad.\footnote{ Once again for any immersion $\vec{\Phi}$ of $D^2$ into ${\R}^m$ $Dis(\nabla\vec{\Phi})=Dis(\vec{\Phi}^\ast g_{{\R}^m})$.}
\]
Observe that
\be
\label{I.7a0}
|Dis(g)|^2=1-4\frac{detg}{(tr\, g)^2}<1\quad.
\ee
\begin{Dfi}
\label{df-I.1}
An admissible measurable complex structure on $\Sigma$ is a measurable section\footnote{ i.e. $J$ is a measurable map
from $\Sigma$ into $End(T\Sigma)$ such that for a.e. $p\in \Sigma$ $J(p)$ is an endomorphism of $T_p\Sigma$, the tangent space
to $\Sigma$ at $p$ such that $J(p)\circ J(p)$ is minus the identity map of $T_p\Sigma$ into itself.}
 $J$ of the endomorphism bundle of $\Sigma$
satisfying $J^2=-Id$ and such that
\be
\label{I.7b0}
\lf\|\ln\frac{|X\wedge JX|_{g_0}}{|X\wedge jX|_{g_0}}\rg\|_{L^\infty((T\Sigma)_0)}<+\infty\quad,
\ee
where $j$ is an arbitrary smooth complex structure on $\Sigma$, the metric $|\cdot|_{g_0}$ on $T\Sigma\wedge T\Sigma$ is induced from an arbitrary reference
metric\footnote{Observe that $\ln\frac{|X\wedge JX|_{g_0}}{|X\wedge J_0X|_{g_0}}$ is independent of the choice of $g_0$}  $g_0$ on $T\Sigma$ and $(T\Sigma)_0$ is the tangent bundle minus the zero section. \hfill $\Box$
\end{Dfi}
\begin{Rm}
\label{rm-I.1}
Observe that for $\vec{\Phi}\in{\mathcal E}_\Sigma$ the complex structure induced by $\vec{\Phi}^\ast g_{{\R}^m}$ is admissible. Indeed
condition(\ref{I.1}) together with the fact that $\vec{\Phi}\in W^{1,\infty}(\Sigma)$ implies that there exists $C_1>0$ such that
\[
\forall p\in{\Sigma}\quad\forall X\in T_p\Sigma\setminus\{0\}\quad\quad C_1^{-1}\ |X|_{g_0}\le |X|_g\le C_1\ |X|_{g_0}\quad.
\]
from which one easily deduce the existence of $C_2$ such that at every point in $\Sigma$
\[
C_2^{-1}\ vol_{g_0}\le vol_{g}\le C_2\ vol_{g_0}
\]
and we deduce (\ref{I.7b0}) by combining the previous equivalences of the metrics and their volume form with the identity
\[
X\wedge [g]\cdot X=|d\vec{\Phi}\cdot X|^2\ vol_g\quad.
\]
where $[g]\cdot$ is the action of the complex structure associated to $g$.
\hfill$\Box$
\end{Rm}
Now given a measurable admissible complex structure $J$ on $\Sigma$, we define the {\it  complex Distortion with respect to $J$} of an immersion
$\Psi\in W^{1,\infty}(D^2,\Sigma)$ to be the function $Dis^J(\Psi)\in L^\infty(D^2,{\C})$ given by
\be
\label{I.7a1}
Dis^J(\Psi):=\frac{|\p_x\Psi|_g^2-|\p_y\Psi|_g^2-2\,i\,\lf<\p_x\Psi,\p_y\Psi\rg>_g}{|\p_x\Psi|^2_g+|\p_y\Psi|^2_g}\quad,
\ee
where $g$ is an arbitrary metric compatible\footnote{$g(J\cdot,J\cdot)=g(\cdot,\cdot)$} with the complex structure $J$.
It is also not difficult to check that $|Dis^J(\Psi)|<1$.

\subsection{The distance functions $d^J$.}

For any measurable admissible complex structure $J$ and non-negative integer $k$, we shall consider the following space of quasi-conformal 
lipschitz parametrization into $\Sigma$ :
\[
{\mathcal Q}^J_k:=\lf\{
\begin{array}{c}
\Psi\in W^{1,\infty}(D^2,\Sigma)\ ,\quad\log|\nabla\Psi|\in L^\infty(D^2)  \\[5mm]
\ds \mbox{ and }\quad \lf\|Dis^J(\Psi)\rg\|_{L^\infty(D^2)}\le 1-2^{-k}\quad\mbox{ a.e. in }D^2
\end{array}
\rg\}
\]
%We also introduce
%\[
%{\mathcal Q}^J_k(\vec{\Phi}):=\lf\{
%\begin{array}{c}
%\ds\Psi\in {\mathcal Q}^J_k\quad\mbox{ s.t. }\quad\quad\int_{D^2}|d\vec{n}_{\vec{\Phi}\circ\Psi}|_g^2\ dvol_g<2\pi/3
%\end{array}
%\rg\}
%\]
%where $g$ denotes on $D^2$ the pull-back metric by $\vec{\Phi}\circ\Psi$ of the standard metric in ${\R^m}$.

On $\mathcal{E}_\Sigma\times \mathcal{E}_\Sigma$ we  introduce the following non negative function

\[
\Delta^J_k(\vec{\Phi}_1,\vec{\Phi}_2):=\sup_{\Psi\in {\mathcal Q}^J_{k}}\ \Lambda(\vec{\Phi}_1,\vec{\Phi}_2,\Psi)\quad,
\]
where
\[
\begin{array}{l}
\ds\Lambda(\vec{\Phi}_1,\vec{\Phi}_2,\Psi):=\lf[\int_{D^2}|\nabla(\vec{n}_1-\vec{n}_2)|^2\ dx\,dy\rg]^\frac{1}{2}\\[5mm]
%   \quad\quad\ds+\lf\|\log_{\C}\lf[1+Dis(\vec{\Phi}_1\circ\Psi)\rg]-\log_{\C}\lf[1+Dis(\vec{\Phi}_2\circ\Psi)\rg]\rg\|_{L^\infty(D^2)}\\[5mm]
   \quad\quad\ds+\lf\|\log|\nabla(\vec{\Phi}_1\circ\Psi)|-\log|\nabla(\vec{\Phi}_2\circ\Psi)|\rg\|_{L^\infty(D^2)}\\[5mm]
   \quad\quad+\lf\|\nabla\lf(\vec{\Phi}_1\circ\Psi-\vec{\Phi}_2\circ\Psi\rg)\rg\|_{L^2(D^2)}\quad,
\end{array}
\]
 where we made use of the following notation
\[
\vec{n}_i:=\vec{n}_{\vec{\Phi}_i\circ\Psi}=\star\frac{\p_x(\vec{\Phi}_i\circ\Psi)\wedge\p_y(\vec{\Phi}_i\circ\Psi)}{|\p_x(\vec{\Phi}_i\circ\Psi)\wedge\p_y(\vec{\Phi}_i\circ\Psi)|}\quad,\quad\quad\mbox{ for }i=1,2\quad .
\]
\begin{Prop}
\label{pr-I.1}
Let $J$ be an admissible measurable complex structure on $\Sigma$, we define  $d^J$ to be the following nonnegative function on $\mathcal{E}_\Sigma\times \mathcal{E}_\Sigma$ 
\[
d^J(\vec{\Phi}_1,\vec{\Phi}_2):=\sum_{k\in {\N}}2^{-k}\Delta^J_k(\vec{\Phi}_1,\vec{\Phi}_2)+\|\ln|X|_{g_1}-\ln|X|_{g_2}\|_{L^\infty((T\Sigma)_0)}
\]
where $g_i:=\vec{\Phi}_i^\ast g_{{\R}^m}$ and $(T\Sigma)_0$ is equal to the tangent bundle to $\Sigma$ minus the zero section.
$d^J$ defines a distance-function on $\mathcal{E}_\Sigma$.\hfill $\Box$
\end{Prop}
{\bf Proof of proposition~\ref{pr-I.1}.}
First we have to prove that $d^J$ is a well defined function. There is indeed a $\sup$ operation and we have to show 
that $d^J(\vec{\Phi}_1,\vec{\Phi}_2)<+\infty$ for any pair $(\vec{\Phi}_1,\vec{\Phi}_2)\in {\mathcal E}_\Sigma\times{\mathcal E}_\Sigma$. 

Since $\vec{\Phi}_1$ and $\vec{\Phi}_2$ are in ${\mathcal E}_\Sigma$, because of (\ref{I.1}) - see also remark~\ref{rm-I.1} - the two metrics
$g_1:=\vec{\Phi}_1^\ast g_{{\R}^m}$ and $g_2:=\vec{\Phi}_2^\ast g_{{\R}^m}$ are equivalent to a reference metric $g_0$ that we assume to be compatible
with $g_0$ : i.e. there exists $C_{g_i,g_0}$ such that
\be
\label{I.7a2}
\forall X\in (T\Sigma)_0\quad\quad C_{g_i,g_0}^{-1}|X|_{g_i}\le |X|_{g_0}\le  C_{g_i,g_0}\ |X|_{g_i}\quad.
\ee
Hence 
\[
\|\ln|X|_{g_1}-\ln|X|_{g_2}\|_{L^\infty((T\Sigma)_0)}<+\infty\quad.
\]
Let $k\in {\N}$ and let $\Psi\in {\mathcal Q}^J_k$. Let $\ov{g}_0$, $\ov{g}_1$ and $\ov{g}_2$ be the 3 following metrics on $D^2$ given by  $\ov{g}_0:=\Psi^\ast g_0$  and $\ov{g}_i:=\Psi^\ast g_i$ where $g_0$ is a reference metric compatible
with $J$. We have
\[
|Dis(\ov{g}_0)|=|Dis^J(\Psi)|<1-2^{-k}\quad.
\]
Hence we deduce from lemma~\ref{lm-A.00} the following inequalities
\be
\label{I.7a3}
\frac{1}{2}\ \inf_{X\in (T\Sigma)_0} \frac{|X|^2_{g_0}}{|X|^2_{g_i}}\le\frac{tr(g_0)}{tr(g_i)}\le 2^k\ \sup_{X\in (T\Sigma)_0} \frac{|X|^2_{g_0}}{|X|^2_{g_i}}\quad .
\ee
From which, together with (\ref{I.7a2}), we deduce that
\be
\label{I.7a4}
2^{-(k+1)}\ C^{-1}_{g_2,g_0}\ C^{-1}_{g_1,g_0}\le\frac{tr(\ov{g}_1)}{tr(\ov{g}_2)}=\frac{|\nabla (\vec{\Phi}_1\circ\Psi)|^2}{|\nabla(\vec{\Phi}_2\circ\Psi)|^2}\le 2^{k+1}\ C_{g_2,g_0}\ C_{g_1,g_0}
\quad.
\ee

Let now $\al=\al_1\, dx+\al_2\, dy$ be a 1-form on $D^2$. Denote $G_0:=(\ov{g}_{0,ij})$ and $A:=(\al_1,\al_2)$. we have
\[
|\al|^2_{\ov{g}_0}\ dvol_{\ov{g}_0}= A\,G_0^{-1}\,A^T\ \sqrt{det(G_0)}\ dx\,dy\quad.
\]
We have also
\[
A\,G_0^{-1}\,A^T\ \sqrt{det(G_0)}\ge\ A\,A^T\ (det(G_0)^{-1})^{-1/2}\inf_{\la\in Spec((G_0)^{-1})}\ \la \quad,
\]
where $Spec((G_0)^{-1})$ denotes the spectrum of the inverse of $G_0$. Hence we have
\[
A\,G_0^{-1}\,A^T\ \sqrt{det(G_0)}\ge\mu\ A\,A^T\quad.
\]
where 
\[
\mu:=\inf\lf\{\sqrt{\frac{\la_1}{\la_2}};\sqrt{\frac{\la_2}{\la_1}}\rg\}\quad,
\]
where $\la_1$ and $\la_2$ are the two eigenvalues of $G_0^{-1}$. Clearly $0<\mu<1$. From (\ref{I.7a0})
we have
\[
\frac{1}{\mu+\frac{1}{\mu}}=\frac{1}{2}\ \sqrt{1-(Dis(\ov{g}_0))^2}\quad.
\]
Hence we deduce, since $|Dis(\ov{g}_0)|=|Dis^J(\Psi)|<1-2^{-k}<1$,
\[
\mu>2^{-1-k/2}\quad.
\]
We deduce from the previous identities that for $i=1,2$
\[
\int_{\Psi(D^2)}\ |d\vec{n}_{\vec{\Phi}_i}|^2_{g_0}\ dvol_{g_0}\ge 2^{-1-k/2}\ \int_{D^2}\ |\nabla \vec{n}_{\vec{\Phi}_i\circ\Psi}|^2\ dx\, dy\quad.
\]
Using now the equivalence of the norms mentioned in remark~\ref{rm-I.1} we obtain that
\be
\label{I.7a5}
2^{1+k/2}\ \int_\Sigma |d\vec{n}_{\vec{\Phi}_i}|^2_{g_i}\ dvol_{g_i}\ge \int_{D^2}\ |\nabla \vec{n}_{\vec{\Phi}_i\circ\Psi}|^2\ dx\, dy\quad.
\ee
In a similar way we deduce also that
\be
\label{I.7a5a}
2^{1+k/2}\ \int_\Sigma |d{\vec{\Phi}_i}|^2_{g_i}\ dvol_{g_i}\ge \int_{D^2}\ |\nabla ({\vec{\Phi}_i\circ\Psi})|^2\ dx\, dy\quad.
\ee
%Observe that 
%\[
%Dis(\vec{\Phi}_i\circ\Psi)=Dis(\Psi^\ast g_i)=Dis^{J_i}(\Psi)\quad,
%\]
%where $J_i$ is a complex structure compatible with $g_i$. Using inequality (\ref{I.7a14}) we know that there exists a constant\footnote{$C_i$ is equal to $4\, e^{\delta(J_i,J)}$  
%where the distance $\delta$ is defined in proposition~\ref{pr-I.2}. } $C_i$
%independent of $\Psi$ such that
%\[
%2^{-k}\le 1-|Dis^{J}(\Psi)|\le C_i\ [1-|Dis^{J_i}(\Psi)]\quad.
%\]
%Using the fact that for a complex number $a$ in the unit disc $|\log_{\C}(1+a)|^2<(\pi/2)^2+(\log|1-|a||^{-1})^2$ one has
%\be
%\label{I.7a5b}
%|\log_{\C}(1+Dis(\vec{\Phi}_i\circ\Psi)|\le C_{\vec{\Phi}_i}\ k
%\ee
We hence deduce from (\ref{I.7a4})   (\ref{I.7a5}) and (\ref{I.7a5a}) that
\be
\label{I.7a6}
\Delta^J_k(\vec{\Phi}_1,\vec{\Phi}_2)=\sup_{\Psi_1\in {\mathcal Q}^J_{k}}\ \Lambda(\vec{\Phi}_1,\vec{\Phi}_2,\Psi)<C_{\vec{\Phi}_1,\vec{\Phi}_2} (2^{k/4}+k)\quad.
\ee
Combining this fact together with (\ref{I.7a3}) we obtain that $d^J(\vec{\Phi}_1,\vec{\Phi}_2)<+\infty$ for any pair $(\vec{\Phi}_1,\vec{\Phi}_2)$ in ${\mathcal E}_\Sigma\times{\mathcal E}_\Sigma$.

\medskip

We prove now that $d^J$ is a distance function.

\medskip

\noindent{\bf Symmetry.} It is clear by definition.

\medskip

\noindent{\bf Triangular inequality.}
 Observe first that $\forall\  (\vec{\Phi}_1,\vec{\Phi}_2)\in({\mathcal E}_\Sigma)^2\quad \forall\ \Psi\in{\mathcal Q}_k^J $
\be
\label{I.7a7}
\forall\  (\vec{\Phi}_1,\vec{\Phi}_2)\in({\mathcal E}_\Sigma)^2\quad \forall\ \Psi\in{\mathcal Q}_k^J \quad\quad\La(\vec{\Phi}_1,\vec{\Phi}_2,\Psi)=\La(\vec{\Phi}_2,\vec{\Phi}_1,\Psi)\quad.
\ee
and moreover $\forall  (\vec{\Phi}_1,\vec{\Phi}_2,\vec{\Phi}_3)\in({\mathcal E}_\Sigma)^3$ and $\Psi\in {\mathcal Q}_k^J$
\be
\label{I.7a8}
\La(\vec{\Phi}_1,\vec{\Phi}_2,\Psi)\le\La(\vec{\Phi}_1,\vec{\Phi}_3,\Psi)+\La(\vec{\Phi}_3,\vec{\Phi}_2,\Psi)\quad.
\ee
Let $\ep>0$, there exists $\Psi\in {\mathcal E}_{\Sigma}$ such that $\forall\ \Psi\in {\mathcal Q}_k^J$
\be
\label{I.7a9}
\Delta_k^J(\vec{\Phi}_1,\vec{\Phi}_2)\le \La(\vec{\Phi}_1,\vec{\Phi}_2,\Psi)+\ep\quad.
\ee
Combining (\ref{I.7a8}) and (\ref{I.7a9}) we obtain for any $\ep>0$
\[
\Delta_k^J(\vec{\Phi}_1,\vec{\Phi}_2)\le\Delta_k^J(\vec{\Phi}_1,\vec{\Phi}_3)+\Delta_k^J(\vec{\Phi}_3,\vec{\Phi}_2)+\ep\quad.
\]
which implies the triangular inequality for $d^J$

\medskip

\noindent{\bf discernibility}. Assume $d^J(\vec{\Phi}_1,\vec{\Phi}_2)=0$. Then for any quasi-conformal map $\Psi\ :\ D^2\rightarrow\Sigma$ we have
$\vec{\Phi}_1\circ\Psi=\vec{\Phi}_2\circ\Psi$. This clearly implies that $\vec{\Phi}_1=\vec{\Phi}_2$.

\medskip

\noindent This concludes the proof of the fact that $d^J$ is a distance function on ${\mathcal E}_\Sigma$ and proposition~\ref{pr-I.1} is
proved. \hfill$\Box$

\begin{Prop}
\label{pr-I.2}
Let $J$ and $J'$ be two admissible measurable complex structures on $\Sigma$, then $d^J$ and $d^{J'}$ are equivalent
distances on ${\mathcal E}_\Sigma$ and there holds
\be
\label{I.7a12}
\begin{array}{l}
\forall\ (\vec{\Phi}_1,\vec{\Phi}_2)\in ({\mathcal E}_\Sigma)^2\\[5mm]
\  8^{-1}\ e^{-\delta(J,J')}\ d^J(\vec{\Phi}_1,\vec{\Phi}_2)\le d^{J'}(\vec{\Phi}_1,\vec{\Phi}_2)
\le 8\ e^{\delta(J,J')}\  d^J(\vec{\Phi}_1,\vec{\Phi}_2)\quad,
\end{array}
\ee
 where 
 $\delta(J,J')$ is the following distance between the two complex structures :
 \[
 \delta(J,J'):=\lf\|\ln\frac{|X\wedge JX|_{g_0}}{|X\wedge J'X|_{g_0}}\rg\|_{L^\infty((T\Sigma)_0)}
 \]
  (for an arbitrary
metric $g_0$).

\hfill $\Box$
\end{Prop}
{\bf Proof of proposition~\ref{pr-I.2}.} Let $\Psi$ be  a quasi-conformal map from the disc $D^2$ into $\Sigma$. (\ref{I.7a1}) implies
\[
1-|Dis^J(\Psi)|^2=4\frac{|\p_x\Psi\wedge\p_y\Psi|_g^2}{\lf[|\p_x\Psi|^2_g+|\p_y\Psi|^2_g\rg]^2}\quad.
\]
where $g$ is an arbitrary metric compatible with $J$. We have hence for instance $|\p_x\Psi|^2_g=|\p_x\Psi\wedge J\p_x\Psi|$
and $|\p_y\Psi|^2_g=|\p_y\Psi\wedge J\p_y\Psi|$. Let ${\mathfrak J}$ be the complex structure such that ${\mathfrak J}\p_x\Psi=\p_y\Psi$.
We have
\[
\begin{array}{rl}
1-|Dis^J(\Psi)|^2&\ds=4\ \lf[\frac{|\p_x\Psi\wedge J\p_x\Psi|_g}{|\p_x\Psi\wedge{\mathfrak J}\p_x\Psi|_g}+\frac{|\p_y\Psi\wedge J\p_y\Psi|_g}{|\p_y\Psi\wedge{\mathfrak J}\p_y\Psi|_g}\rg]^{-1}\\[5mm]
 &\ds =4\ \lf[\frac{|\p_x\Psi\wedge J\p_x\Psi|_{g_0}}{|\p_x\Psi\wedge{\mathfrak J}\p_x\Psi|_{g_0}}+\frac{|\p_y\Psi\wedge J\p_y\Psi|_{g_0}}{|\p_y\Psi\wedge{\mathfrak J}\p_y\Psi|_{g_0}}\rg]^{-1}
\end{array}
\]
where $g_0$ is an arbitrary reference metric on $\Sigma$.
Using an elementary algebraic inequality\footnote{\[
\forall\  a,a',b,b'>0\quad\quad
\frac{1}{a+b}\le\lf[\frac{a'}{a}+\frac{b'}{b}\rg]\ \frac{1}{a'+b'}\quad.
\]} we obtain that for any quasi-conformal map $\Psi$ from $D^2$ into $\Sigma$
\be
\label{I.7a13}
\begin{array}{rl}
\ds\frac{1-|Dis^J(\Psi)|^2}{1-|Dis^{J'}(\Psi)|^2}&\ds\le \lf[\frac{|\p_x\Psi\wedge J'\p_x\Psi|_{g_0}}{|\p_x\Psi\wedge J\p_x\Psi|_{g_0}}+
\frac{|\p_y\Psi\wedge J'\p_y\Psi|_{g_0}}{|\p_y\Psi\wedge J\p_y\Psi|_{g_0}}\rg]\\[5mm]
 &\ds \le 2 e^{\delta(J,J')}
 \end{array}
\ee
Hence we have
\be
\label{I.7a14}
1-|Dis^J(\Psi)|\le 4\ e^{\delta(J,J')}\ [1-|Dis^{J'}(\Psi)]
\ee
let $k_0=[\delta(J,J')/\log 2]+1$. We have that
\[
|Dis^J(\Psi)|<1-2^{-k}\quad\quad\Longrightarrow\quad\quad |Dis^{J'}(\Psi)|<1-2^{-k-k_0-2}\quad.
\]
Hence
\[
\forall k\in {\N}\quad\quad 2^{-k}\ \Delta_k^J(\vec{\Phi}_1,\vec{\Phi}_2)\le 2^{k_0+2}\ 2^{-k-k_0-2}\ \Delta_{k+k_0+2}^{J'}(\vec{\Phi}_1,\vec{\Phi}_2)\quad,
\]
from which we deduce
\[
d^{J}(\vec{\Phi}_1,\vec{\Phi}_2)\le 2^{k_0+2}\ d^{J'}(\vec{\Phi}_1,\vec{\Phi}_2)\quad.
\]
This last inequality implies proposition~\ref{pr-I.2}.\hfill $\Box$

\subsection{ Completeness of the metric spaces $({\mathcal E}_\Sigma,d^J)$.}
In this subsection we prove the following result.
\begin{Prop}
\label{pr-II.1}
For any  admissible measurable complex structure $J$ on $\Sigma$, the metric space $({\mathcal E}_\Sigma,d^J)$ is complete.
\hfill $\Box$
\end{Prop}
{\bf Proof of proposition~\ref{pr-II.1}.}
Because of the equivalence of the $d^J$ we can choose an arbitrary $J$ that we will assume to be smooth. We first choose a finite covering of $\Sigma$
by open sets $U_i$ such that each $U_i$ is diffeomorphic to a disc and we denote by $\Psi_i$ diffeomorphisms from $D^2$ into $U_i$
such that $Dis^J(\Psi_i)=0$.
 Let now $\vec{\Phi}_k\in{\mathcal E}_{\Sigma}$ such that $d(\vec{\Phi}_{k-1},\vec{\Phi}_k)\le 2^{-k-1}$.
 Denotes $g_k:=\vec{\Phi}_k^\ast g_{{\R}^m}$. The assumption implies that
 $\vec{\Phi}_k\circ\Psi_i$ converges strongly in $W^{1,2}$ to a limit $\vec{\xi}_i$ : $D^2\rightarrow {\R}^m$.
 Hence $|\nabla\vec{\Phi}_k\circ\Psi_i|$ converges a.e. to $|\nabla \vec{\xi}_i|$ and hence $\log|\nabla\vec{\Phi}_k\circ\Psi_i|$
 converges a.e. to $\log|\nabla \vec{\xi}_i|\in {\R}\cup\{+\infty\}\cup\{-\infty\}$.
 From the Cauchy sequence assumption for $\vec{\Phi}_k$ w.r.t. $d^J$ we have that $\log|\nabla\vec{\Phi}_k\circ\Psi_i|$ is Cauchy in $L^\infty$
 and this limit can only be $\log|\nabla\vec{\xi}_i|$ which is then in $L^\infty(D^2)$.
 We have moreover
 \[
 \forall\ j,l=1,2\quad\quad \Psi_i^\ast g_{k,jl}=\p_{x_j}(\vec{\Phi}_k\circ\Psi_i)\cdot\p_{x_l}(\vec{\Phi}_k\circ\Psi_i)\longrightarrow \p_{x_j}\vec{\xi}_i\cdot\p_{x_l}\vec{\xi}_i\quad\mbox{a.e.}\quad.
 \]
Since 
\be
 \label{I.7b1}
\forall k\in {\N}\quad\quad\quad \|\log|X|_{g_k}-\log|X|_{g_1}\|_{L^\infty((T\Sigma)_0)}\le 2^{-1}
\ee
we have that
\[
\forall \, X\in{\R}^2\ne 0 \quad\quad \log \lf[\sum_{j,l=1}^2\Psi_i^\ast g_{k,jl}X^j\,X^l\rg]\longrightarrow \log \lf[\sum_{j,l=1}^2\p_{x_j}\vec{\xi}_i\cdot\p_{x_l}\vec{\xi}_i\, X^j\, X^l\rg]\quad\mbox{ a.e. }
\]
and hence
\be
\label{I.7b2}
\forall \, X\in{\R}^2\ne 0 \quad\quad\lf|\log\lf[ \sum_{j,l=1}^2\p_{x_j}\vec{\xi}_i\cdot\p_{x_l}\vec{\xi}_i\ X^j\, X^l \rg]\rg|\le \log|X|_{g_1}+2^{-1}\quad.
\ee
We deduce from the previous inequality that $\vec{\xi}_i$ is an immersion from $D^2$ into ${\R}^m$ and there exists $c_i>0$ such that
\be
\label{I.7b3}
|d\vec{\xi}_i\wedge d\vec{\xi}_i|\ge c_i>0\quad\quad\mbox{ on }U_i\quad.
\ee
Hence the $\vec{\xi}_i$ are lipschitz immersions. We have, for any pair $i,j$,
$\vec{\xi}^{-1}_i\circ\vec{\xi}_j=\Psi_i^{-1}\circ\Psi_j$, hence there exists $\vec{\xi}$ a Lipschitz immersion from $\Sigma$ into ${\R}^m$ such that
$\vec{\xi}_i=\vec{\xi}\circ\Psi_i$ for all $i$. 
Let
\[
\vec{n}_{i,k}:=\vec{n}_{\vec{\Phi}_k\circ\Psi_i}=\star\frac{\p_x(\vec{\Phi}_k\circ\Psi_i)\wedge\p_y(\vec{\Phi}_k\circ\Psi_i)}{|\p_x(\vec{\Phi}_k\circ\Psi_i)\wedge\p_y(\vec{\Phi}_k\circ\Psi_i)|}\quad\quad\quad .
\]
From the Cauchy sequence assumption for $\vec{\Phi}_k$ we deduce that $\vec{n}_{i,k}$ converges strongly in $\dot{W}^{1,2}(D^2)$ 
to a limit that we denote $\vec{n}_i$, moreover, from the above we have that
\[
\star\frac{\p_x(\vec{\Phi}_k\circ\Psi_i)\wedge\p_y(\vec{\Phi}_k\circ\Psi_i)}{|\p_x(\vec{\Phi}_k\circ\Psi_i)\wedge\p_y(\vec{\Phi}_k\circ\Psi_i)|}
\longrightarrow
\star\frac{\p_x\vec{\xi}_i\wedge\p_y\vec{\xi}_i}{| \p_x\vec{\xi}_i\wedge\p_y\vec{\xi}_i|}\quad\quad\quad\mbox{ a.e. }\quad.
\]
Hence we have that
\be
\label{I.7b4}
\star\frac{\p_x\vec{\xi}_i\wedge\p_y\vec{\xi}_i}{| \p_x\vec{\xi}_i\wedge\p_y\vec{\xi}_i|}=\vec{n}_i\quad\in \quad W^{1,2}(D^2)\quad.
\ee 
Hence $\vec{\xi}$ is a Lipschitz immersion inducing a metric comparable to the smooth reference metric $g_0$ whose Gauss map is in $W^{1,2}$
with respect to this metric. This implies that $\vec{\xi}\in {\mathcal E}_\Sigma$ and it can be proved with moderate efforts that $d^J(\vec{\Phi}_k,\vec{\xi})\longrightarrow 0$. This concludes the proof of proposition~\ref{pr-II.1}.\hfill $\Box$

\subsection{Control of $d^g(\vec{\Phi},\vec{\Phi}+t\vec{w})$ for conformal $\vec{\Phi}$ and for $W^{2,2}\cap W^{1,\infty}$ perturbations $\vec{w}$.}

\begin{Lm}
\label{lm-II.1}
There exists $\ep_0>0$ such that for any conformal immersion $\vec{\Phi}$  of the disc $D^2$ into ${\R}^m$ in $W^{2,2}\cap W^{1,\infty}$  satisfying
\be
\label{II.00}
\int_{D^2}|\nabla\vec{n}_{\vec{\Phi}}|^2\ dx\ dy<\frac{4\pi}{3}\quad,
\ee
Let $\vec{w}\in W^{1,\infty}\cap W^{2,2}(D^2,{\R}^m)$ compactly supported in $D^2_{1/2}$ such that
\be
\label{II.1}
\|\nabla \vec{w}\|_{L^\infty(D^2)}+\|\nabla^2\vec{w}\|_{L^2(D^2)}\le 1
\ee
Denote $\vec{\Phi}_t:=\vec{\Phi}+t\vec{w}$. Then there exists $C>0$ independent of $\vec{\Phi}$
and $\vec{w}$ such that, for $|t|<[\inf_{D^2_{1/2}}|\nabla \vec{\Phi}|]/4$,
\be
\label{II.2}
\begin{array}{l}
\ds d^g(\vec{\Phi}_t,\vec{\Phi})\le \frac{C}{\inf_{D^2_{1/2}}|\nabla\vec{\Phi}|}\ |t|\ \|\nabla^2 w\|_2\\[7mm]
 \ds\quad\quad\quad\quad+\frac{C}{\inf_{D^2_{1/2}}|\nabla\vec{\Phi}|}\lf[1+\frac{\|\nabla\vec{\Phi}\|_\infty\ \|\nabla\vec{n}\|_2+\|\nabla\vec{\Phi}\|_2}{\inf_{D^2_{1/2}}|\nabla\vec{\Phi}|}\rg]\ |t|\ \|\nabla w\|_\infty\quad.
 \end{array}
\ee
where $g:=\vec{\Phi}^\ast g_{{\R}^m}$.\hfill$\Box$
\end{Lm}
{\bf Proof of lemma~\ref{lm-II.1}.}

We denote the conformal factor as usual as follows : $e^\la:=|\p_x\vec{\Phi}|=|\p_y\vec{\Phi}|$. 
Denote $e^{\ov{\la}}=\|\nabla\vec{\Phi}\|_\infty$ and $e^{\underline{\la}}:=\inf_{D^2_{1/2}} e^\la$. Consider $t$ such that $4\,|t|< e^{\underline{\la}}$. Since
$w$ is supported in $D^2_{1/2}$ and since $\|\nabla w\|_\infty\le 1$, we have
\[
\begin{array}{rl}
\ds|\p_x\vec{\Phi}_t\wedge\p_y\vec{\Phi}_t|&\ds\ge e^{2\la}-|t|\ |\p_x\vec{w}\wedge\p_y\vec{\Phi}|-|t|\ |\p_x\vec{\Phi}\wedge\p_y\vec{w}|\\[5mm]
 &\ds\ -t^2 |\p_x\vec{w}\wedge\p_y\vec{w}|\\[5mm]
 &\ds\ge e^{2\la}-2\,|t|\,e^{\la}-t^2\ge \frac{7}{16}\ e^{2\la}=\frac{7}{16}\ |\p_x\vec{\Phi}\wedge\p_y\vec{\Phi}|
 \end{array}
 \]
A straightforward but a bit lengthy computation shows that
\be
\label{II.15a}
|\nabla(\vec{n}_{\vec{\Phi}}-\vec{n}_{\vec{\Phi}_t})|\le C\ |t|\ e^{-\la}\ \lf[ |\nabla \vec{w}|\ e^{-\la}\ |\nabla^2\vec{\Phi}|+|\nabla^2\vec{w}|\rg]\quad,
\ee
where $C$ is independent of all the datas $\vec{\Phi}$, $\vec{w}$ and $t$. Since 
\[
\Delta\vec{\Phi}=2e^{2\la}\ \vec{H}
\]
where $\vec{H}$ is the mean curvature vector of the immersion of $D^2$ which is pointwisely controled by $|\nabla \vec{n}|$, standard
elliptic estimates imply
\[
\int_{D^2_{1/2}}|\nabla^2\vec{\Phi}|^2\le C\ e^{2\ov{\la}}\ \int_{D^2}|\nabla\vec{n}|^2+C\ \int_{D^2}|\nabla\vec{\Phi}|^2\quad.
\]
Integrating hence (\ref{II.15a}) on $D^2$, since $|\nabla(\vec{n}_{\vec{\Phi}}-\vec{n}_{\vec{\Phi}_t})|$ is supported on 
$D^2_{1/2}$, we obtain
\be
\label{II.16}
\begin{array}{l}
\ds\|\nabla(\vec{n}_{\vec{\Phi}}-\vec{n}_{\vec{\Phi}_t})\|_2
\le \frac{C}{\inf_{D^2_{1/2}}|\nabla\vec{\Phi}|}\lf[ |t|\ \|\nabla^2 \vec{w}\|_2\rg.\\[5mm]
\quad\quad\quad\quad\quad\quad\quad\quad\ds+\frac{\|\nabla\vec{\Phi}\|_\infty\ \|\nabla\vec{n}\|_2+\|\nabla\vec{\Phi}\|_2}{\inf_{D^2_{1/2}}|\nabla\vec{\Phi}|}\ |t|\ \|\nabla\vec{w}\|_\infty\big]\quad.
\end{array}
\ee
%The conformality of $\vec{\Phi}$ is equivalent to $Dis(\vec{\Phi})=0$. Again for  $4\,|t|< e^{\underline{\la}}$ we have
%\[
%|Dis(\vec{\Phi}_t)|\le\ \frac{C}{|\nabla\vec{\Phi}|}\ |t|\ |\nabla \vec{w}|\quad.
%\]
%Thus we obtain
%\be
%\label{II.17}
%\|Dis(\vec{\Phi})-Dis(\vec{\Phi}_t)\|_\infty\le \frac{C}{\inf_{D^2_{1/2}}|\nabla\vec{\Phi}|}\ |t|\ \|\nabla \vec{w}\|_\infty\quad.
%\ee
We have the pointwise identity
\[
|\nabla\vec{\Phi}_t-\nabla\vec{\Phi}|\le |t|\ |\nabla \vec{w}|
\]
and hence
\be
\label{II.017a}
\int_{D^2}|\nabla\vec{\Phi}_t-\nabla\vec{\Phi}|^2\le C\  |t|^2\ \|\nabla \vec{w}\|_\infty^2\quad.
\ee
We have also 
\[
\lf|\log\frac{|\nabla\vec{\Phi}_t|^2}{|\nabla\vec{\Phi}|^2}\rg|=\lf|
\log\lf[1+2e^{-2\la}\ t\ \nabla\vec{w}\cdot\nabla\vec{\Phi}+e^{-2\la}\ t^2\ |\nabla\vec{w}|^2\rg]\rg|\quad.
\]
Hence for $4\ |t|<e^{\underline{\la}}$ we deduce
\be
\label{II.17a}
\lf\|\log|\nabla\vec{\Phi}_t|-\log|\nabla\vec{\Phi}|\rg\|_{L^\infty(D^2)}\le  \frac{C}{\inf_{D^2_{1/2}}|\nabla\vec{\Phi}|}\ |t|\ \|\nabla \vec{w}\|_\infty\quad.
\ee
Combining (\ref{II.16}), (\ref{II.017a}) and  (\ref{II.17a}), we obtain
\be
\label{II.18}
\begin{array}{l}
\ds\La(\vec{\Phi},\vec{\Phi}_t,id_{D^2})\le \frac{C}{\inf_{D^2_{1/2}}|\nabla\vec{\Phi}|}\ |t|\ \|\nabla^2 \vec{w}\|_2\\[7mm]
 \ds\quad+\frac{C}{\inf_{D^2_{1/2}}|\nabla\vec{\Phi}|}\lf[1+\frac{\|\nabla\vec{\Phi}\|_\infty\ \|\nabla\vec{n}\|_2+\|\nabla\vec{\Phi}\|_2}{\inf_{D^2_{1/2}}|\nabla\vec{\Phi}|}\rg]\ |t|\ \|\nabla w\|_\infty\ .
 \end{array}
\ee
Let now $\Psi$ be an arbitrary map in ${\mathcal Q}_k^g(\vec{\Phi})$. Since $\vec{\Phi}$ is conformal $\Psi$ is a quasiconformal map
satisfying
\[
|Dis^g(\Psi)|=|Dis(\vec{\Phi}\circ\Psi)|\le1-2^{-k}\quad .
\]
This implies
\be
\label{II.18a}
(2^{k+1}-1)^{-1}|\p_y\Psi|^2<|\p_x\Psi|^2\le\,(2^{k+1}-1)\,|\p_y\Psi|^2\quad,
\ee
and 
\be
\label{II.18b}
\lf|\frac{\p_x\Psi}{|\p_x\Psi|}\cdot \frac{\p_y\Psi}{|\p_y\Psi|}\rg|^2<1-2^{-k-1}\quad.
\ee
After some short computation, we deduce from the previous line
\be
\label{II.18c}
2^{-5/2}\,2^{-3k/2}\ |\nabla\Psi|^2\le \det\nabla\Psi=\p_x\Psi\times\p_y\Psi\le|\nabla\Psi|^2/2\quad.
\ee
Hence we have
\be
\label{II.19}
\begin{array}{l}
\ds\int_{D^2}|\nabla(\vec{n}_{\vec{\Phi}\circ\Psi}-\vec{n}_{\vec{\Phi}_t\circ\Psi})|^2\ dx\,dy\\[5mm]
\ds\quad\quad\le\int_{D^2}|\nabla(\vec{n}_{\vec{\Phi}}-\vec{n}_{\vec{\Phi}_t})|^2\circ\Psi\ |\nabla\Psi|^2\ dx\,dy\\[5mm]
\ds\quad\quad\le 2^{5/2}\,2^{3k/2}\ \int_{D^2}|\nabla(\vec{n}_{\vec{\Phi}}-\vec{n}_{\vec{\Phi}_t})|^2\circ\Psi\ \det\nabla\Psi\ dx\,dy\\[5mm]
\ds\quad\quad\le 2^{5/2}\,2^{3k/2}\ \int_{D^2}|\nabla(\vec{n}_{\vec{\Phi}}-\vec{n}_{\vec{\Phi}_t})|^2\ dx\,dy\\[5mm]
\ds\quad\quad\le 2^{5/2}\,2^{3k/2}\ \La^2(\vec{\Phi},\vec{\Phi}_t,id_{D^2})\quad.
\end{array}
\ee
Similarly we have
\be
\label{II.20}
\begin{array}{l}
\ds\int_{D^2}|\nabla(\vec{\Phi}\circ\Psi-\vec{\Phi}_t\circ\Psi)|^2\ dx\,dy\\[5mm]
\ds\quad\quad\le
 2^{5/2}\,2^{3k/2}\ \int_{D^2}|\nabla(\vec{\Phi}-\vec{\Phi}_t)|^2\circ\Psi\ \det\nabla\Psi\ dx\,dy\\[5mm]
 \ds\quad\quad\le 2^{5/2}\,2^{3k/2}\ \La^2(\vec{\Phi},\vec{\Phi}_t,id_{D^2})\quad.
\end{array}
\ee

We have moreover, since $|\nabla(\vec{\Phi}\circ\Psi)|=e^\la\ |\nabla\Psi|$,
\[
\begin{array}{l}
\ds\lf|\log\frac{|\nabla(\vec{\Phi}_t\circ\Psi)|^2}{|\nabla(\vec{\Phi}\circ\Psi)|^2}\rg|=\lf|
\log\lf[1+2e^{-2\la}\ |\nabla\Psi|^{-2} t\ \nabla(\vec{w}\circ\Psi)\cdot\nabla(\vec{\Phi}\circ\Psi)\rg.\rg.\\[5mm]
 \ds\quad\quad\quad\lf.\lf.+e^{-2\la}\ t^2\ |\nabla\Psi|^{-2}\ |\nabla\vec{w}\circ\Psi|^2\rg]\rg|\quad.
 \end{array}
\]
Using the fact that $|\nabla(\vec{w}\circ\Psi)|\le|\nabla\vec{w}|\ |\nabla\Psi|$, we then have for $|t|<[\inf_{D^2_{1/2}}|\nabla \vec{\Phi}|]/4$
\be
\label{II.23a}
\lf\|\log|\nabla(\vec{\Phi}_t\circ\Psi)|-\log|\nabla(\vec{\Phi}\circ\Psi)|\rg\|_{L^\infty(D^2)}\le  \frac{C}{\inf_{D^2_{1/2}}|\nabla\vec{\Phi}|}\ |t|\ \|\nabla \vec{w}\|_\infty\quad.
\ee
Hence combining (\ref{II.19}), (\ref{II.20}) and (\ref{II.23a}) we have obtained the existence of $C>0$ independent of $\vec{\Phi}$
and $\vec{w}$ such that for $|t|<[\inf_{D^2_{1/2}}|\nabla \vec{\Phi}|]/4$, for any $\Psi\in {\mathcal Q}^{[g]}_k$
\be
\label{II.24}
\begin{array}{l}
\ds\La(\vec{\Phi},\vec{\Phi}_t,\Psi)\le \frac{2^{3k/4}\ C}{\inf_{D^2_{1/2}}|\nabla\vec{\Phi}|}\ |t|\ \|\nabla^2 \vec{w}\|_2\\[7mm]
 \ds\quad\quad\quad+\frac{2^{3k/4}\ C}{\inf_{D^2_{1/2}}|\nabla\vec{\Phi}|}\lf[1+\frac{\|\nabla\vec{\Phi}\|_\infty\ \|\nabla\vec{n}\|_2+\|\nabla\vec{\Phi}\|_2}{\inf_{D^2_{1/2}}|\nabla\vec{\Phi}|}\rg]\ |t|\ \|\nabla \vec{w}\|_\infty\quad.
 \end{array}
\ee
where $C$ is independent of $\Psi$. Hence we deduce
\be
\label{II.024}
\begin{array}{l}
\ds\sum_{k\in{\N}}2^{-k}\Delta^{[g]}_k(\vec{\Phi}_t,\vec{\Phi})\le \frac{ C}{\inf_{D^2_{1/2}}|\nabla\vec{\Phi}|}\ |t|\ \|\nabla^2 \vec{w}\|_2\\[7mm]
 \ds\quad\quad\quad+\frac{C}{\inf_{D^2_{1/2}}|\nabla\vec{\Phi}|}\lf[1+\frac{\|\nabla\vec{\Phi}\|_\infty\ \|\nabla\vec{n}\|_2+\|\nabla\vec{\Phi}\|_2}{\inf_{D^2_{1/2}}|\nabla\vec{\Phi}|}\rg]\ |t|\ \|\nabla \vec{w}\|_\infty\quad.
 \end{array}
\ee
at a point $p\in D^2$ we have for $X=X_1\,\p_{x_1}+X_2\,\p_{x_2}$, denoting $|X|^2_0=X_1^2+X_2^2$,
\[
\frac{|X|_{g_{\vec{\Phi}_t}}}{|X|_{g_{\vec{\Phi}}}}=1+2\ e^{-2\la}\ t\ \nabla\vec{w}\frac{X}{|X|_0}\cdot\nabla\vec{\Phi}\frac{X}{|X|_0}+e^{-2\la}\ t^2\ |\nabla\vec{w}\frac{X}{|X|_0}|^2
\]
Hence for $|t|<e^{\underline{\la}}$ we have that
\be
\label{II.025}
\|\log\,|X|_{g_{\vec{\Phi}_t}}-\log\,|X|_{g_{\vec{\Phi}}}\|_{L^\infty((T\Sigma)_0)}\le  \frac{C}{\inf_{D^2_{1/2}}|\nabla\vec{\Phi}|}\ |t|\ \|\nabla \vec{w}\|_\infty
\ee
where we are using the fact that, for such $t$ and $\vec{w}$,  $e^{-\la}\ |t|\ |\nabla\vec{w}|<1$. Inequality (\ref{II.024})
together with inequality (\ref{II.025}) imply (\ref{II.2}) and Lemma~\ref{lm-II.1} is proved.\hfill $\Box$

\section{Existence of Minimizers of the Willmore Energy.}
\reset

In this section we prove the following result

\begin{Th}
\label{th-III.1}
Let $\Sigma$ be an abstract closed two dimensional smooth manifolds. Assume that there exists a minimizing sequence
$\vec{\Phi}_k\in {\mathcal E}_\Sigma$ of the Willmore energy such that the conformal class induced by $\vec{\Phi}_k^\ast g_{{\R}^m}$
stays in the compact subset of the Riemann Moduli Space of $\Sigma$ and assume that
\[
\limsup_{k\rightarrow +\infty}W(\vec{\Phi}_k)= \int_{\Sigma}|\vec{H}_k|^2_{g_k}\ dvol_{g_k}<8\pi\quad .
\]
where $g_k:=\vec{\Phi}_k^\ast g_{{\R}^m}$ and $\vec{H}_k$ is the mean-curvature vector of the immersion $\vec{\Phi}_k$, then
\[
\inf_{\vec{\Phi}\in{\mathcal E}_\Sigma}W(\vec{\Phi})
\]
is achieved by a smooth embedding.\hfill $\Box$
\end{Th}

Before to prove theorem~\ref{th-III.1} we first state the following proposition which is a direct application of Ekeland's Variational
Principle (see theorem 5.1 in \cite{St}).

\begin{Prop}
\label{pr-III.1}
Let $J$ be an arbitrary smooth complex structure on $\Sigma$ and let $\vec{\Phi}_k$ be a minimizing sequence such that
\[
W(\vec{\Phi}_k)\le\inf_{\vec{\Phi}\in{\mathcal E}_\Sigma}W(\vec{\Phi})+2^{-k}\quad,
\]
then there exists $\vec{\xi}_k\in {\mathcal E}_\Sigma$ such that 
\begin{itemize}
\item[i)]
$\vec{\xi}_k$ minimizes in ${\mathcal E}_\Sigma$ the following functional
\be
\label{III.a01}
W(\vec{\xi}_k)=\inf_{\vec{\Phi}\in{\mathcal E}_\Sigma}W(\vec{\Phi})+2^{-k/2}\ d^J(\vec{\Phi},\vec{\xi}_k)\quad ,
\ee
\item[ii)]
\be
\label{III.a02}
W(\vec{\xi}_k)\le W(\vec{\Phi}_k)\quad,
\ee
\item[iii)]
\be
\label{III.a03}
d^J(\vec{\xi}_k,\vec{\Phi}_k)\le 2^{-k/2}\quad.
\ee
\end{itemize}
 \hfill$\Box$
\end{Prop}
{\bf Proof of theorem~\ref{th-III.1}.}
We can assume that $\Sigma$ is not $S^2$ since  a classical result implies that
\[
\inf_{\vec{\Phi}\in {\mathcal E}_{S^2}}W(\vec{\Phi})=4\pi
\]
and is achieved by the unit sphere of ${\R}^3\subset{\R}^m$ (see for instance \cite{Wi}). Let $\vec{\Phi}_k$ be a minimizing sequence of the Willmore energy $W$ in the space ${\mathcal E}_\Sigma$ satisfying
\be
\label{III.1}
W(\vec{\Phi}_k)\le\inf_{\vec{\Phi}\in{\mathcal E}_\Sigma}W(\vec{\Phi})+2^{-k}\quad.
\ee
 A straightforward
molification argument allows to work under the assumption that $\vec{\Phi}_k\in C^\infty(\Sigma,{\R}^m)$. The assumption
of the theorem is telling us that the conformal class of the induced metric $\vec{\Phi}_k^\ast g_{{\R}^m}$ is contained
in a compact subset of the Moduli space of $\Sigma$. Therefore, modulo extraction of a subsequence, we can find a sequence
of complex structure $J_k$ and diffeomorphisms $f_k$ of $\Sigma$ such that 
\[
\vec{\Phi}_k\circ f_k\ :\ (\Sigma,J_k)\longrightarrow {\R}^m\quad\mbox{ is conformal }\quad,
\]
and, if $h_k$ denotes 
\be
\label{III.1a}
J_{k}\longrightarrow J_\infty \quad\mbox{ w.r.t. }\delta\quad.
\ee
$\vec{\Phi}_k\circ f_k$ satisfies of course still (\ref{III.1}). Denote by $h_k$, resp.  $h_\infty$, the smooth constant scalar curvature metric
compatible with $J_k$ resp. $J_\infty$ having a fixed volume 1 on $\Sigma$, we may also ensure that
\be
\label{III.2a}
\|\log\,|X|_{h_k}-\log\,|X|_{h_\infty}\|_{L^\infty((T\Sigma)_0)}\longrightarrow 0\quad .
\ee
To each $x\in\Sigma$ we assign $\rho_x>0$ such that 
\[
\int_{B_{\rho_x}(x)}|d\vec{n}|^2_{h_k}\ dvol_{h_k}=\int_{B_{\rho_x}(x)}|d\vec{n}|^2_{g_k}\ dvol_{g_k}=8\pi/3\quad,
\]
where $B_{\rho_x}(x)$ is the geodesic ball in $(\Sigma,h_k)$ of center $x$ and radius $\rho_x$ and
 $g_k:=f_k^\ast\vec{\Phi}_k^\ast g_{{\R}^m}$.
 We extract a finite Besicovitch covering : each point in $\Sigma$ is covered by at most $N$ of such balls where $N$
 only depends on  $(\Sigma,g_\infty)$. Let $(B_{\rho_k^i}(x_k^i))_{i\in I}$ be this finite covering. We can extract a subsequence
 such that $I$ is independent of $k$, such that each $x_k^i$ converges to a limit $x_\infty^i$ and each $\rho_k^i$
 converges to a limit $\rho_\infty^i$.
 Let
 \[
 I_0:=\{i\in I\quad\mbox{ s. t. }\quad\rho_\infty^i=0\}\quad.
 \] 
 Let $I_1:=I\setminus I_0$. It is clear that the union of the closures of the balls $\cup_{i\in I_1}\ov{B}_{\rho_\infty^i}(x_\infty^i)$ covers $\Sigma$.
 Because of the strict convexity of the balls with respect either to the euclidian distance ($\Sigma=T^2$) or the hyperbolic distance
 (genus$(\Sigma)>1$) the points in $\Sigma$ which are not contained in the union of the \underbar{open} balls
  $\cup_{i\in I}B_{\rho_\infty^i}(x^i_\infty)$ cannot accumulate and therefore are isolated and hence finite. Denote
  \be
  \label{III.2}
  \{a_1\cdots a_N\}:=\Sigma\setminus\cup_{i\in I_1}B_{\rho_\infty^i}(x^i_\infty)\quad.
  \ee
We are now using the following lemma in order to ''normalize'' the embeddings $\vec{\Phi}_k\circ f_k$.
\begin{Lm}
\label{lm-III.1}{\bf[3-points normalization lemma]}
For any $\La>0$ there exists $R, r>0$ such that for any closed two dimensional manifold $\Sigma$, for any choice
of 3 distinct points $P_1, P_2$ and $P_3$ in $\Sigma$ and for any embedding $\vec{\Phi}$ of $\Sigma$ into ${\R}^m$
satisfying
\be
\label{III.3}
\int_{\Sigma} |d\vec{n}_{\vec{\Phi}}|^2_g\ dvol_g<\La\quad,
\ee
where $g:=\vec{\Phi}^\ast g_{{\R}^m}$, then there exists a Moebius transformation $\Xi$ of ${\R}^m$ such that 
\be
\label{III.4}
\Xi\circ\vec{\Phi}(\Sigma)\subset B_R(0)\quad\mbox{ and }\quad\forall i\ne j \quad |\Xi\circ\vec{\Phi}(P_i)-\Xi\circ\vec{\Phi}(P_j)|\ge r\quad.
\ee
Moreover the following control of the total area of $\Xi\circ\vec{\Phi}(\Sigma)$ holds
\be
\label{III.5}
{\mathcal H}^2(\Xi\circ\vec{\Phi}(\Sigma))\le C\ R^2\ \La\quad,
\ee
where $C>0$ is a universal constant.
\hfill$\Box$
\end{Lm}  
{\bf Proof of lemma~\ref{lm-III.1}.}  
 We apply a translation and a dilation in such a way that $P_1=0$ and 
 \[
 |\vec{\Phi}(P_1)-\vec{\Phi}(P_2)|=\min_{i\ne j}|\vec{\Phi}(P_i)-\vec{\Phi}(P_j)|=1\quad.
 \]
 We keep denoting $\vec{\Phi}$ the resulting embedding - observe that due to it's conformal invariance the Willmore energy has
 not been modified.
 From  ({A.6}) in \cite{KS2} there exists a universal constant $C>0$ such that for any $x_0\in{\R}^m$ and $0<\sigma<\rho<+\infty$
 \be
 \label{III.6}
 \sigma^{-2}\ {\mathcal H}^2(\vec{\Phi}(\Sigma)\cap B_\sigma(x_0))\le C\ \lf[\rho^{-2}\ {\mathcal H}^2(\vec{\Phi}(\Sigma)\cap B_\rho(x_0))+\int_{\vec{\Phi}^{-1}(B_\rho(x_0))}|\vec{H}|^2\ dvol_{g}\rg]\ .
 \ee
 \medskip
 
 We claim that there exists $\rho_0$ depending only on $\La>W(\vec{\Phi})$ and $x_1\in B_1(0)$ such that
 $\vec{\Phi}(\Sigma)\cap B_{\rho_0}(x_1)=\emptyset$. For $y\in \vec{\Phi}(\Sigma)$ one has 
 \be
 \label{III.7}
 \lim_{\sigma\rightarrow 0} \sigma^{-2}\ {\mathcal H}^2(\vec{\Phi}(\Sigma)\cap B_\sigma(y))=\pi \quad.
 \ee
 For $0<\rho<1/2$ we consider a regular covering of $B_1(0)$ by balls $B_\rho(z_l)$ in such a way that any point in $B_1(0)$ is contained
 in at most $C(m)$ balls of the form $B_{2\rho}(z_l)$ . The number of $l$ 
 such that $\int_{\vec{\Phi}^{-1}(B_{2\rho}(z_l))}|\vec{H}|^2\ dvol_{g}>C^{-1}\,\pi/2$ is bounded by $2\,\La\,C\, C(m)$. For an $l$ such that
$\int_{\vec{\Phi}^{-1}(B_{2\rho}(z_l))}|\vec{H}|^2\ dvol_{g}<C^{-1}\,\pi/2$ and such that there exists $y\in B_\rho(z_l)\cap\Sigma\ne\emptyset$,
combining (\ref{III.6}) and (\ref{III.7}) one obtains that
\[
(2\rho)^{-2}\ {\mathcal H}^2(\vec{\Phi}(\Sigma)\cap B_{2\rho}(z_l))>C^{-1}\,\pi/2\quad.
\]
the number of such $l$ is then bounded by $\rho^{-2}$ times a number depending only on $m$ and $\La$ - where we are using again
(\ref{III.6}) but for $x_0=0$, $\sigma=1$ and $\rho\rightarrow+\infty$. The total number of ball $B_\rho(z_l)$ is proportional to $\rho^{-m}$.
Since $m>2$, for $\rho=\rho_0$ chosen small enough, depending only on $m$ and $\La$ we deduce the claim.

\medskip

Let $x_1$ and $\rho_0$ given by the claim we choose $\Xi$ to be the inversion with respect to $x_1$ : $\Xi(x):=(x-x_1)/|x-x_1|^2$. We have then
\be
\label{III.7a}
\Xi(\vec{\Phi}(\Sigma))\subset B_{1/\rho_0}(0)\quad,
\ee
moreover,  since none of the $P_i$ is in $B_{\rho_0}(x_1)$, we have that $\forall i=1,2,3$ $\Xi(\vec{\Phi}(P_i))\in B_{1/\rho_0}(0)$. We have also
that $|\vec{\Phi}(P_2)-x_1|+|\vec{\Phi}(P_1)-x_1|<3$ hence $\Xi(\vec{\Phi}(P_1))$ and $\Xi(\vec{\Phi}(P_2))$ are contained in ${\R}^m\setminus B_{1/3}(0)$ thus
\[
|\Xi(\vec{\Phi}(P_1))-\Xi(\vec{\Phi}(P_2))|\ \|\nabla \Xi^{-1}\|_{L^\infty(B_{1/\rho_0}(0)\setminus B_{1/3}(0))}|\ge |\vec{\Phi}(P_1)-\vec{\Phi}(P_2)|=1\quad,
\]
which implies that
\be
\label{III.8}
|\Xi(\vec{\Phi}(P_1))-\Xi(\vec{\Phi}(P_2))|\ge 9\quad.
\ee
Either $\vec{\Phi}(P_3)\in B_{10}(0)$ or $\vec{\Phi}(P_3)\in {\R}^m\setminus B_{10}(0)$. In the first case one has that all the $\Xi(\vec{\Phi}(P_i))$
are included in $B_{1/\rho_0}(0)\setminus B_{1/11}(0)$ hence we have
\[
\begin{array}{l}
\ds\forall\ i\ne j\quad |\Xi(\vec{\Phi}(P_i))-\Xi(\vec{\Phi}(P_j))|\ \|\nabla \Xi^{-1}\|_{L^\infty(B_{1/\rho_0}(0)\setminus B_{1/11}(0))}\\[5mm]
\ds\quad\quad\quad\ge |\vec{\Phi}(P_i)-\vec{\Phi}(P_j)|\ge  1
\end{array}
\]
which implies that
\be
\label{III.9}
\forall\ i\ne j\quad |\Xi(\vec{\Phi}(P_i))-\Xi(\vec{\Phi}(P_j))|\ge 11^2\quad.
\ee
This implies (\ref{III.4}) in this case. In the case when $\vec{\Phi}(P_3)\in {\R}^m\setminus B_{10}(0)$ we deduce that $\Xi(\vec{\Phi}(P_3))\in B_{1/9}(0)$ and since $\Xi(\vec{\Phi}(P_1))$ and $\Xi(\vec{\Phi}(P_2))$ are contained in ${\R}^m\setminus B_{1/3}(0)$, we obtain that
\be
\label{III.10}
\forall\ i=1,2\quad |\Xi(\vec{\Phi}(P_i))-\Xi(\vec{\Phi}(P_3))|\ge {2}/{9}\quad.
\ee
This lower bound combined with (\ref{III.8}) gives (\ref{III.4}) in this case too.

Regarding the proof of estimate (\ref{III.5}), we first observe that inequality (\ref{III.6}) (which holds also for $\vec{\Phi}$ replaced by $\Xi\circ\vec{\Phi}$) implies that
, for any $\rho\ge \rho_0^{-1}$
\[
{\mathcal H}^2(\Xi\circ\vec{\Phi}(\Sigma))\le C\rho_0^{-2}\rho^{-2}\  {\mathcal H}^2(\Xi\circ\vec{\Phi}(\Sigma))+\rho_0^{-2}\ \La\quad.
\]
Letting $\rho$ converge to $+\infty$ yields the desired estimate (\ref{III.5}). Hence lemma~\ref{lm-III.1} is proved.\hfill $\Box$

\medskip

\noindent{\bf Proof of theorem~\ref{th-III.1} continued.}

For the minimizing sequence $\vec{\Phi}_k\circ f_k$ we consider the ''Normalization Moebius Transformations'' $\Xi_k$ given by
lemma~\ref{lm-III.1} : we fix 3 distinct points $P_1$, $P_2$ and $P_3$ in the interior of one of the  balls ${B}_{\rho_\infty^i}(x_\infty^i)$  for ${\rho_\infty^i}\ne 0$
(i.e. $i\in I_1$) and we replace our minimizing sequence $ \vec{\Phi}_k\circ f_k$ by $\Xi_k\circ\vec{\Phi}_k\circ f_k$.

In order to simplify the notations we then write $\vec{\Phi}_k$ instead of $\Xi_k\circ\vec{\Phi}_k\circ f_k$. Hence from now on we have a sequence of complex structures
$J_k$ on $\Sigma$ such that
\[
J_{k}\longrightarrow J_\infty \quad\mbox{ w.r.t. }\delta\quad.
\]
with associated constant scalar curvature metrics $h_k$ of volume 1 and satisfying
\[
\|\log\,|X|_{h_k}-\log\,|X|_{h_\infty}\|_{L^\infty((T\Sigma)_0)}\longrightarrow 0\quad .
\]
and we have a sequence of smooth immersions $\vec{\Phi}_k$ of $\Sigma$ into ${\R}^m$ satisfying (\ref{III.1}) and the following five conditions
\begin{itemize}
\item[i)]
\be
\label{III.11}
\vec{\Phi}_k\quad\quad\mbox{ is conformal from } (\Sigma,J_k)\quad\mbox{ into }{\R}^m\quad.
\ee
\item[ii)] There exists finitely many points $a_1\cdots a_N$ in $\Sigma$ and a fixed finite covering $({B}_{\rho_\infty^i}(x_\infty^i))_{i\in I_1}$
of $\Sigma\setminus\{a_1\cdots a_N\}$ such that for any $i\in I_1$, $0<\rho<\rho_\infty$ and $k$ large enough
\be
\label{III.12}
\int_{B_\rho(x_\infty^i)}|d\vec{n}_{\vec{\Phi}_k}|^2_{g_k}\ dvol_{g_k}<8\pi/3\quad .
\ee
where $g_k:=\vec{\Phi}_k^\ast g_{{\R}^m}$.
\item[iii)] There exists a positive real $R>0$ such that
\be
\label{III.13}
\vec{\Phi}_k(\Sigma)\subset B_R(0)\quad.
\ee
\item[iv)] There exists a constant $C>0$ such that
\be
\label{III.14}
{\mathcal H}^2(\vec{\Phi}_k(\Sigma))\le C\quad .
\ee
\item[v)] There exist a positive  real number $r>0$, independent of $k$ and three distinct points $P_1$, $P_2$ and $P_3$, independent of $k$ too, in the interior of one ball $B_{\rho_\infty^i}(x_\infty^i)$
such that
\be
\label{III.15}
\forall\ i\ne j\quad\quad\quad|\vec{\Phi}_k(P_i)-\vec{\Phi}_k(P_j)|\ge r>0\quad\quad.
\ee
\end{itemize}

\noindent Using now proposition~\ref{pr-III.1} we construct $\vec{\xi}_k$ satisfying (\ref{III.a01}), (\ref{III.a02}) and (\ref{III.a03}).

\medskip

We claim now the following 

\begin{Lm}
\label{lm-III.0}
For any compact $K\subset \Sigma\setminus\{a_1\cdots a_N\}$ there exists $C_K>0$ and $k_K\in {\N}$ such that
\be
\label{III.16}
\sup_{k\ge k_K}\|\log|d \vec{\Phi}_k|_{h_k}\|_{L^\infty(K)}\le C_K<+\infty\quad.
\ee
\hfill $\Box$
\end{Lm}
\noindent{\bf Proof of the lemma~\ref{lm-III.0}.} For any compact subset $K$ of $\Sigma\setminus\{a_1\cdots a_N\}$ there exists $\delta>0$ such that
$K\subset \Sigma\setminus\cup_{i=1}^N \ov{B}_\delta(a_i)$. Since $\Sigma\setminus\cup_{i=1}^N {B}_\delta(a_i)\subset \cup_{i\in I_1} B_{\rho_\infty^i}(x_\infty^i)$,
there exist $\rho^i_\infty>r_i>0$ such that
\be
\label{III.17}
\Sigma\setminus\cup_{i=1}^N {B}_\delta(a_i)\subset \cup_{i\in I_1} B_{r^i}(x_\infty^i)
\ee
and for $k$ large enough one has for any $i\in I_1$ $B_{r^i}(x_\infty^i)\subset B_{\rho_k^i}(x_k^i)$. We may also have chosen $r^i$ large enough in such a way that
$P_1$, $P_2$ and $P_3$ belong to the same $B_{r^i}(x^i_\infty)$ for $i=i_P$ say.

 Let $s^i=(r^i+\rho^i_\infty)/2$. We consider $k$ large enough in such a way that 
$B_{s^i}(x_\infty^i)\subset B_{\rho_k^i}(x_k^i)$ for any $i\in I_1$. Since, on $B_{s^i}(x^i_\infty)$ $\vec{\Phi}_k$ is conformal and (\ref{III.12}) holds for $\rho=s^i$, using 
lemma 5.1.4 of \cite{Hel}, we deduce the existence of
 $\vec{e}_1$ and $\vec{e}_2$ in $W^{1,2}(B_{s^i}(x_\infty^i),S^{m-1})$ such that 
\be
\label{III.18}
\vec{e}_1\cdot \vec{e}_2=0\quad\quad,\quad\quad n_{\vec{\Phi}_k}=\vec{e}_1\wedge \vec{e}_2 \quad,
\ee
\be
\label{III.19}
\int_{B_{s^i}(x_\infty^i)}\lf[|\nabla \vec{e}_1|_{g_k}^2+|\nabla \vec{e}_2|_{g_k}^2\rg]\ dvol_{g_k}\le 2\int_{B_{s^i}(x_\infty^i)}|\nabla\vec{n}_{\vec{\Phi}_k}|^2_{g_k}\ dvol_{g_k}
\ee
where $g_k:=\vec{\Phi}_k^\ast\, g_{{\R}^m}$ and
\be
\label{III.20}
\lf\{
\begin{array}{l}
\ds d(\ast_{g_k}(\vec{e}_1,d \vec{e}_2))=0\quad,\\[5mm]
\ds\iota_{\p B_{s^i}(x_\infty^i)}^\ast\ast_{g_k}(\vec{e}_1,d \vec{e}_2)=0\quad .
\end{array}
\rg.
\ee
where $\iota_{\p B_{s^i}(x_\infty^i)}$ is the canonical inclusion of $\p B_{s^i}(x_\infty^i)$ in ${\R}^2$. 

$(\Sigma, h_k)$ identifies to a fundamental domain $(\Xi_k,g_{hyp})$ of the Poincar\'e upper half-plane $H$ where $\Xi_k$ converges
to a  limiting domain $\Xi_\infty$. On each geodesic ball $B_{s^i}(x_\infty^i)$ (We can always assume that the geodesic  ball $B_{s^i}(x_\infty^i)$ is connected in the fundamental, otherwise we change the domain for each $i$) we consider the canonical coordinates $(x,y)$ of $H$.
We shall denote also $\Psi^i$ the coordinate map, that modulo a dilation we can always assume to go from $D^2$ into $B_{s^i}(x_\infty^i)$.
The metric $h_k$ on such a ball is conformally equivalent to the flat metric $g^i_{flat}$ given by the canonical coordinates $(x,y)$: $h_k= e^{\sigma}\ g^i_{flat}$ with $\|\sigma\|_{L^\infty(B_{s^i}(x_\infty^i))}$ uniformly bounded independently of $k$. It suffices then to prove
that $\|\log|\nabla( \vec{\Phi}_k\circ\Psi^i)|\|_\infty=\|\log|d\vec{\Phi}_k|_{g^i_{flat}}\|$ is uniformly bounded w.r.t. $k$ for every $i$ in order to show (\ref{III.16}).

\medskip

Since $\vec{\Phi}_k$ is conformal, the logarithm of the conformal factor $\la_k=\log e^{\la_k}=\log|d \vec{\Phi}_k|_{g^i_{flat}}$ - i.e. $g_k=e^{2\la_k}\, (dx^2+dy^2)$ - satisfies on $B_{s^i}(x_\infty^i)$ the following equation
\be
\label{III.21}
-\Delta\la_k={\mathcal K}\,e^{2{\la_k}}=(\nabla^\perp \vec{f}_1,\nabla\vec{f}_2)=\p_x\vec{f}_1\cdot\p_y\vec{f}_2-\p_x\vec{f}_2\cdot\p_y\vec{f}_1\quad,
\ee
where ${\mathcal K}$ is the scalar curvature and $(\vec{f}_1,\vec{f}_2)$ is the positive orthonormal Basis of $T_{\vec{\Phi}(x,y)}\vec{\Phi}_k(B_{s^i}(x_\infty^i))$ given by $(\vec{f}_1,\vec{f}_2)=e^{-\la_k}(\p_x\vec{\Phi}_k,\p_y\vec{\Phi}_k)$. A classical computation yields
\be
\label{III.22}
div((\vec{f}_1,\nabla \vec{f}_2))=0\quad .
\ee
There exists $\theta$ such that $\vec{f}=e^{i\theta}\vec{e}$ (where we are identifying $\vec{e}$ (resp. $\vec{f}$) with $\vec{e}_1+i\vec{e}_2\in {\C}\otimes{\R}^m$ (resp. $\vec{f}_1+i\vec{f}_2$)). $\theta$ satisfies 
\[
(\vec{f}_1,\nabla\vec{f}_2)=\nabla\theta+(\vec{e}_1,\nabla\vec{e}_2)\quad
\]
This gives that $\theta$ is harmonic but also that $(\nabla^\perp \vec{e}_1,\nabla\vec{e}_2)=(\nabla^\perp \vec{f}_1,\nabla\vec{f}_2)$. Hence we have
\be
\label{III.23}
-\Delta\la_k=(\nabla^\perp \vec{e}_1,\nabla\vec{e}_2)\quad.
\ee
Because of (\ref{III.14}) we have that
Let $\mu_k$ be the function on $B_{s^i}(x_\infty^i)$ satisfying
\be
\label{III.24}
\lf\{
\begin{array}{l}
\ds -\Delta\mu^i_k=(\nabla^\perp \vec{e}_1,\nabla\vec{e}_2)\quad\mbox{ on }B_{s^i}(x_\infty^i)\\[5mm]
\mu^i_k=0\quad\quad\quad\mbox{ on }\p B_{s^i}(x_\infty^i)\quad.
\end{array}
\rg.
\ee
Using Wente's inequality (\cite{We} or \cite{Hel} theorem 3.1.2) we obtain
\be
\label{III.25}
\|\mu^i_k\|_{L^\infty(B_{s^i}(x_\infty^i))}\le C\int_{B_{s^i}(x_\infty^i)}|\nabla\vec{n}_{\vec{\Phi}_k}|^2\ dx\,dy\quad.
\ee
Let $\nu^i_k:=\la_k-\mu_k^i$. $\nu^i_k$ is harmonic and, hence, $e^{2\nu_k^i}$ is subharmonic ($\Delta e^{2\nu_k^i}\ge 0$). 
Let $r^i<t^i<s^i$ be independent of $k$.
The subharmonicity of $e^{2\nu_k^i}$  implies  the existence of constant $C^i$, independent of $k$, such that
\be
\label{III.26}
\begin{array}{l}
\ \ds\sup_{x\in B_{t^i}(x_\infty^i)} e^{\nu_k^i}\le C^i\ \lf[\int_{B_{s^i}(x_\infty^i)} e^{2\nu_k^i}\rg]^{1/2}\\[5mm]
\quad\quad\le C^i\lf[{\mathcal H}^2(\vec{\Phi}_k(B_{s^i}(x_\infty^i)))\rg]^{1/2}\ e^{\|\mu_k^i\|_\infty}\quad.
\end{array}
\ee
Combining this inequality with (\ref{III.14}) and (\ref{III.25}) we obtain
\be
\label{III.27}
\sup_{x\in B_{t^i}(x_\infty^i)}\la_k\le \log\lf[C^i\ {\mathcal H}^2(\vec{\Phi}_k(\Sigma))\ W(\vec{\Phi}_k)\rg]
\ee
and hence $\la_k$ is bounded from above on $B_{r^i}(x_\infty^i)$ independently of $k$. 

Consider now the index $i=i_P$ corresponding to the ball $B_{r^{i_P}}(x^{i_P}_\infty)$ which contains the
3 points $P_1$, $P_2$ and $P_3$ satisfying (\ref{III.15}). Since $\nu_k^i$ is harmonic, the Poisson representation formula
reads
\[
\forall p\in B_{t^{i_P}}(x^{i_P}_\infty)\quad\quad \nu_k^{i_P}(p)=\frac{1}{2\pi}\int_{\p B_{{t^{i_P}}}}\nu_k^{i_P}(z)\ \frac{|t^{i_P}|^2-|p|^2}{|z-p|^2}\ d\sigma(z)\quad.
\]
Combining this identity together with (\ref{III.26}) implies the existence of a constant $C^{i_P}$ independent of $k$ such that
\be
\label{III.28}
\begin{array}{l}
\ds\forall p\in B_{r^{i_P}}(x^{i_P}_\infty)\\[5mm]
\ds\quad\quad \nu_k^{i_P}(p)\le C^{i_P}\ \int_{\p B_{{t^{i_P}}}}(\nu_k^{i_P})^-(z)\ \ d\sigma(z)+ C^{i_P}W(\vec{\Phi}_k)\ .
\end{array}
\ee
where $(\nu_k^{i_P})^-:=\min\{\nu_k^{i_P},0\}$. From (\ref{III.15}) there exists $r>0$ independent of $k$ such that
\[
\begin{array}{rl}
\ds r<\int_{[P_1,P_2]} e^{\la_k}&\ds\le e^{\|\mu_k^{i_P}\|_\infty}\ \int_{[P_1,P_2]}e^{\nu_k^{i_P}}\\[5mm]
 &\ds\le e^{C^{i_P}W(\vec{\Phi}_k)}\ \exp\lf[C^{i_P}\,\int_{\p B_{{t^{i_P}}}}(\nu_k^{i_P})^-(z)\ \ d\sigma(z)\rg]\quad .
 \end{array}
\]
Hence $\int_{\p B_{{t^{i_P}}}}(\nu_k^{i_P})^-(z)\ \ d\sigma(z)$ is bounded from below independently of $k$. This implies
\[
\limsup_{k\rightarrow +\infty}\int_{\p B_{{t^{i_P}}}}|\nu_k^{i_P}|(z)\ \ d\sigma(z)<+\infty
\]
Using again the representation formula (\ref{III.28}) and the bound (\ref{III.25}) we obtain finally that
\be
\label{III.29}
\limsup_{k\rightarrow +\infty}\|\la_k\|_{L^\infty(B_{r^{i_P}}(x^{i_P}_\infty))}<+\infty\quad.
\ee
This is exactly the claim (\ref{III.16}) but on the ball $B_{r^{i_P}}(x^{i_P}_\infty)$. Considering now any other ball
$B_{r^i}(x^i_\infty)$ which intersection with $B_{r^{i_P}}(x^{i_P}_\infty)$ is non empty. Therefore one can find a 
non trivial segment $[Q_1,Q_2]\subset B_{r^i}(x^i_\infty)\cap B_{r^{i_P}}(x^{i_P}_\infty) $ for which then
\[
\liminf_{k\rightarrow +\infty}\int_{[Q_1,Q_2]}e^{\la_k}>0
\]
Hence we have all the ingredients in order to repeat the argument above and deduce (\ref{III.29}) for this $B_{r^i}(x^i_\infty)$ too.
We iterate this procedure until having reached every ball $B_{r^i}(x^i_\infty)$ for $i\in I_1$ since $\Sigma$ is assumed to be connected.
Hence the claim (\ref{III.16}) is proved and this finishes the proof of lemma~\ref{lm-III.0}. \hfill $\Box$

\medskip

\noindent{\bf Proof of theorem~\ref{th-III.1} continued .}  Since , from (\ref{III.a03}), $d^{J}(\vec{\Phi}_k,\vec{\xi}_k)\le 2^{-k/2}$, by taking on each $B_{r^i}(x^i_\infty)$, $\Psi^i$ to be the canonical
coordinate map $(x,y)$ of the Poincar\'e half plane ${\mathbb H}$ once $B_{r^i}(x^i_\infty)$ has been identified with a connected part of a fundamental domain associated to $(\Sigma,J_k)$ we have 
\be
\label{III.29a}
\begin{array}{l}
\|\log|d \vec{\Phi}_k|_{h_k}-\log|d \vec{\xi}_k|_{h_k}\|_{L^\infty(B_{r^i}(x^i_\infty))}\\[5mm]
\quad\quad\le\|\log|\nabla(\vec{\Phi}_k\circ\Psi^i)|-\log|\nabla(\vec{\xi}_k\circ\Psi^i)|\|_{L^\infty(D^2)}\\[5mm]
\quad\quad\le d^{J_k}(\vec{\Phi}_k,\vec{\xi}_k)\le C\ 2^{-k/2}
\end{array}
\ee
Where we have used proposition~\ref{pr-I.2} and the fact that $\delta(J_k,J)$ is uniformly bounded since $J_k$ converges to a limit $J_\infty$.
Hence we deduce that for any compact $K\subset \Sigma\setminus\{a_1\cdots a_N\}$ there exists $C_K>0$ and $k_K\in {\N}$ such that
\be
\label{III.30}
\sup_{k\ge k_K}\|\log|d \vec{\xi}_k|_{h_k}\|_{L^\infty(K)}\le C_K<+\infty\quad.
\ee
Moreover, since for the same $\Psi^i$ on $B_{r^i}(x^i_\infty)$, $\vec{\Phi}_k\circ \Psi^i$ is conformal from $D^2$ into ${\R}^m$, combining the fact that
$$
\|\log\,|X|_{g_{\vec{\xi}_k}}-\log\,|X|_{g_{\vec{\Phi}_k}}\|_{L^\infty((T\Sigma)_0)}\le d^{J_k}(\vec{\xi}_k,\vec{\Phi}_k)\le\ C\  2^{-k/2}
$$
together with inequality (\ref{III.29a}) and identitiy (\ref{A.02}), we obtain
\be
\label{III.31}
\lf\|\log_{\C}\lf[1+Dis(\vec{\xi}_k\circ\Psi^i)\rg]\rg\|_{L^\infty(D^2)}\le 2^{-k/2}\quad.
\ee
This implies that for any compact $K\subset \Sigma\setminus\{a_1\cdots a_N\}$ there exists $C_K>0$ for any $k\in {\N}$
\be
\label{III.32}
\lf\|\log\frac{|d\vec{\xi}_k\cdot X\wedge d\vec{\xi}_k\cdot J_kX|}{|d\vec{\xi}_k\cdot X|^2}\rg\|_{L^\infty(TK)}\le C_K\ 2^{-k/2}\quad.
\ee
We have proved that
\be
\label{III.33}
\sup_{k\ge k_K}\|\log|\nabla \vec{\Phi}_k\circ\Psi^i|\|_{L^\infty(B_{r^i}(x^i_\infty))}\le C^i<+\infty\quad.
\ee
Moreover 
\be
\label{III.34}
\begin{array}{rl}
\ds 4^{-1}\sup_{k\ge k_K}\int_{D^2} |\Delta(\vec{\Phi}_k\circ\Psi^i)|^2 e^{-2\la_k}\ dx\,dy&\ds=\sup_{k\ge k_K}\int_{B_{r^i}(x^i_\infty)}|\vec{H}|^2_{g_k}\ dvol_{g_k}
\\[5mm]
 &\ds\le 8\pi/3\quad.
 \end{array}
\ee
Hence, combining (\ref{III.33}) and (\ref{III.34}) we deduce that 
\be
\label{III.35a}
\sup_{k\ge k_K}\int_{\cup_{i\in I_1}B_{r^i}(x^i_\infty)}|\Delta_{h_k}\vec{\Phi}_k|^2\ dvol_{h_k}<+\infty\quad.
\ee
Combining this fact and the fact that
\be
\label{III.35}
\begin{array}{l}
\ds\sup_{k\ge k_K}\int_{\cup_{i\in I_1}B_{r^i}(x^i_\infty)} |d\vec{\Phi}_k|_{h_k}^2\ dvol_{h_k}=\sup_{k\ge k_K}\int_{\cup_{i\in I_1}B_{r^i}(x^i_\infty)} |d\vec{\Phi}_k|_{g_k}^2\ dvol_{g_k}\\[7mm]
\ds\quad\quad\quad\quad\quad\le \sup_{k\ge k_K}{\mathcal H}^2(\vec{\Phi}_k(\Sigma))<+\infty\quad,
\end{array}
\ee
we have that, modulo extraction of a subsequence, $d\vec{\Phi}_k$ converges strongly in $L^p(K)$ ($\forall p<+\infty$) w.r.t. $h_k$ (which itself converges to $h_\infty$ in every norm).
This implies that $\nabla(\vec{\Phi}_k\circ\Psi^i)$ converges strongly in $L^p(K\cap B_{r^i}(x^i_\infty))$ for all $p<+\infty$ . From (\ref{III.29a}) and (\ref{III.31}) we deduce
that $d\vec{\xi}_k$ converges also strongly in $L^p(K)$ to a limit $d\vec{\xi}_\infty$ where $\vec{\xi}_\infty$ is a Lipschitz conformal immersion
of $(K,J_\infty)$ into ${\R}^m$. We have then, using also (\ref{III.30}), for all $p<+\infty$
\be
\label{III.36}
\ds\vec{n}_{\vec{\xi}_k}=\ast_{h_k}\frac{d\vec{\xi}_k\wedge d\vec{\xi}_k}{|d\vec{\xi}_k\wedge d\vec{\xi}_k|_{h_k}}\longrightarrow
\ast_{h_\infty}\frac{d\vec{\xi}_\infty\wedge d\vec{\xi}_\infty}{|d\vec{\xi}_\infty\wedge d\vec{\xi}_\infty|_{h_\infty}}=\vec{n}_{\vec{\xi}_\infty}\quad\mbox{ in }L^p(K)
\ee
From the definition of $d^{J}$ have for any $i\in I_1$
\be
\label{III.37}
\begin{array}{rl}
\ds\int_{D^2}|\nabla(\vec{n}_{\vec{\Psi}_k}\circ\Psi^i-\vec{n}_{\vec{\xi}_k}\circ\Psi^i)|^2\ dx\,dy&\ds\le d^{J_k}(\vec{\xi}_k,\vec{\Phi}_k)\\[5mm]
 &\ds\le C\ d^{J}(\vec{\xi}_k,\vec{\Phi}_k)\le C\ 2^{-k/2}
 \end{array}
\ee
Hence we have
\be
\label{III.38}
\begin{array}{l}
\ds\limsup_{k\rightarrow+\infty}\int_{K}|d\vec{n}_{\vec{\xi}_k}|^2_{h_k}\ dvol_{h_k}\le
 \limsup_{k\rightarrow+\infty}\int_{K}|d\vec{n}_{\vec{\Phi}_k}|^2_{h_k}\ dvol_{h_k}\\[5mm]
\ds \quad\quad\quad\quad= \limsup_{k\rightarrow+\infty}\int_{K}|d\vec{n}_{\vec{\Phi}_k}|^2_{g_k}\ dvol_{g_k}<+\infty
 \end{array}
 \ee
Combining (\ref{III.36}) and (\ref{III.38}) we deduce that
\be 
\label{III.39}
\begin{array}{l}
\ds\int_{K}|d\vec{n}_{\vec{\xi}_\infty}|^2_{h_\infty}\ dvol_{h_\infty}=\int_{K}|d\vec{n}_{\vec{\xi}_\infty}|^2_{g_\infty}\ dvol_{g_\infty}\\[5mm]
\ds\quad\quad\le\limsup_{k\rightarrow+\infty}\int_{K}|d\vec{n}_{\vec{\Phi}_k}|^2_{g_k}\ dvol_{g_k}\le \inf_{\vec{\Phi}\in{\mathcal E}_\Sigma} 4W(\vec{\Phi})-4\pi\chi(\Sigma)\quad,
\end{array}
\ee
where $g_\infty:=\vec{\xi}^\ast_\infty g_{{\R}^m}$. Hence, by iterating the previous facts for a sequence of compacts $K_l$ of $\Sigma\setminus\{a_1\cdots a_N\}$ $\vec{\xi}_\infty$ such that $\cup_lK_l=\Sigma\setminus\{a_1\cdots a_N\}$ one obtains that $\vec{\xi}_\infty$ realizes a conformal, locally lipschitz, immersion of $\Sigma\setminus\{a_1\cdots a_N\}$ such that
\be
\label{III.40}
\int_{\Sigma\setminus\{a_1\cdots a_N\}}|d\vec{n}_{\vec{\xi}_\infty}|^2_{g_\infty}\ dvol_{g_\infty}\le \inf_{\vec{\Phi}\in{\mathcal E}_\Sigma} 4W(\vec{\Phi})-4\pi\chi(\Sigma)\quad.
\ee
We claim now that 
\begin{Lm}
\label{lm-III.2}
\be
\label{III.41}
\vec{\xi}_\infty\mbox{ is a {\it Conformal Willmore immersion} of }\Sigma\setminus\{a_1\cdots a_N\}
\ee
and hence is analytic on $\Sigma\setminus\{a_1\cdots a_N\}$ (see \cite{BR}).\hfill $\Box$
\end{Lm}

\medskip

\noindent{\bf Proof of the lemma~\ref{lm-III.2}.} 
 Denote $\gamma_k$ the metric $\gamma_k=\vec{\xi}_k^\ast g_{{\R}^m}$. Because of (\ref{III.31})
 We have that
 \[
 \|Dis^{J_k}(\vec{\xi_k})\|_{L^\infty(B_{r^i}(x^i_\infty))}\le C\ 2^{-k/2}\quad\longrightarrow\quad 0
 \]
 We can then apply lemma~\ref{lm-A.3} in order to obtain the existence of a diffeomorphism
 $\zeta_k$ from $D^2$ into $D^2$ such that $\vec{\Xi}_k^i:=\vec{\xi}_k\circ\Psi^i\circ\zeta_k$ is conformal
 and satisfy
\be
\label{III.41a}
\limsup_{k\rightarrow+\infty}\|\zeta_k\|_{C^{0,\al}(D^2)}+\|\zeta_k^{-1}\|_{C^{0,\al}(D^2)}<+\infty
\ee
and for any $\rho<1$,  (\ref{A.11}) implies 
\be
\label{III.41b}
\limsup_{k\rightarrow+\infty}\|\log|\nabla\vec{\Xi}_k|\|_{L^\infty(D^2_\rho)}+\|\vec{\Xi}_k^j\|_{W^{2,2}(D^2_\rho)}<+\infty\quad.
\ee
Combining this with (\ref{III.29a}) and (\ref{III.31}) we also obtain
 \be
 \label{III.41c}
\limsup_{k\rightarrow+\infty}\|\log|\nabla\vec{\zeta}_k|\|_{L^\infty(D^2_\rho)}<+\infty \quad.
 \ee
 Let $\rho<1$ and $\vec{w}\in W^{1,\infty}\cap W^{2,2}(D^2,{\R}^m)$ such that $\vec{w}\in C^\infty_0(D^2_\rho,{\R}^m)$ and 
$$
\|\nabla\vec{w}\|_{L^\infty(D^2)}+\|\nabla^2\vec{w}\|_{L^2(D^2)}\le 1\quad.
$$
 Denote $\gamma_k$ the
metric $\gamma_k=\vec{\xi}_k^\ast g_{{\R}^m}$. Because of (\ref{III.30}) and (\ref{III.31}) we have, because of proposition~\ref{pr-I.2}, 
\be
\label{III.42}
d^{\gamma_k}\ \simeq\ d^{g_k}\ \simeq\ d^J\quad\quad\mbox{ indep. of }k\quad.
\ee
Since $\vec{\xi}_k$ minimizes $W(\cdot)+2^{-k/2}\ d^J(\cdot,\vec{\xi}_k)$, we have that for any such $\vec{w}$ and for $|t|$ small enough, independent of $k$, say $|t|<t_0$, denoting $$\vec{\xi}_k^t:=\vec{\xi}_k+t\ \chi_{(B_{r^i}(x^i_\infty))}\ \vec{w}\circ\zeta_k^{-1}\circ(\Psi^i)^{-1}\in{\mathcal E}_\Sigma$$
where $\chi_{(B_{r^i}(x^i_\infty))}$ is the characteristic function of the ball $B_{r^i}(x^i_\infty)$,
\be
\label{III.43}
\begin{array}{rl}
\ds W(\vec{\xi}_k)&\ds\le W(\vec{\xi}_k^t)+2^{-k/2}\ d^J(\vec{\xi}_k^t,\vec{\xi}_k)\\[5mm]
 &\ds\le W(\vec{\xi}_k^t)
+C\ 2^{-k/2}\ d^{\gamma_k}(\vec{\xi}_k^t,\vec{\xi}_k)
\end{array}
\ee
Using now lemma~\ref{lm-II.1}, we deduce the existence of a constant $C>0$ independent of $k$, $\vec{w}$ and $t$ such that
\be
\label{III.44}
W(\vec{\xi}_k)\le W(\vec{\xi}_k^t)+C\ 2^{-k/2}\ |t|\ \lf[\|\nabla \vec{w}\|_{\infty}+\|\nabla^2\vec{w}\|_{2}\rg]
\ee
We have
\[
W(\vec{\xi}_k)-W(\vec{\xi}_k^t)=4^{-1}\int_{B_{r^i}(x^i_\infty)}|d\vec{n}_{\vec{\xi}_k}|^2_{\gamma_k}\ dvol_{\gamma_k}-|d\vec{n}_{\vec{\xi}^t_k}|^2_{\gamma^t_k}\ dvol_{\gamma^t_k}
\]
where $\gamma^t_k:=(\vec{\xi}_k^t)^\ast g_{{\R}^m}$. A straightforward but a bit lengthy argument shows that
\be
\label{III.45}
\begin{array}{rl}
\ds|d\vec{n}_{\vec{\xi}^t_k}|^2_{\gamma^t_k}\ dvol_{\gamma^t_k}=|d\vec{n}_{\vec{\xi}_k}|^2_{\gamma_k}\ dvol_{\gamma_k}&\ds+t\ \vec{w}\circ\zeta_k^{-1}\circ(\Psi^i)^{-1}\cdot\vec{F}_{\vec{\xi}_k}\ dvol_{\vec{\xi}_k}\\[5mm]
 \\+t^2G(\vec{\xi}_k,\vec{w},t)\ dvol_{\vec{\xi}_k}
 \end{array}
\ee
where 
\be
\label{III.46}
\limsup_{k\rightarrow+\infty}\sup_{|t|<t_0}\int_{B_{r^i}(x^i_\infty)}|G(\vec{\xi}_k,\vec{w},t)|\ dvol_{\vec{\xi}_k}<+\infty\quad.
\ee
Now a classical computation from Blashke \cite{Bla} for $m=3$ and \cite{Wei} for arbitrary $m$  gives for a regular immersion $\vec{\xi}$
from $B_{r^i}(x^i_\infty)$ into ${\R}^m$ that
\be
\label{III.47}
4^{-1}\vec{F}_{\vec{\xi}}=\Delta_\perp\vec{H}_{\vec{\xi}}+\ti{A}(\vec{H}_{\vec{\xi}})-2|\vec{H}_{\vec{\xi}}|^2\vec{H}_{\vec{\xi}}
\ee
where $\Delta_\perp$ is the negative covariant laplacian on the normal bundle of the immersion $\vec{\xi}$, moreover for any $\vec{L}\in {\R}^m$, $\ti{A}(\vec{L}):=\sum_{i,j=1}^2\vec{B}(\vec{e}_i,\vec{e}_j)\ \vec{B}(\vec{e}_i,\vec{e}_j)\cdot\vec{L}$
where $\vec{B}$ is the second fundamental form of the immersion $\vec{\xi}$. At this stage it is very important to observe that we are computing
$\vec{F}_{\vec{\xi}}$ for a smooth immersion $\xi$. It does not make sense for an immersion in ${\mathcal E}_\Sigma$ such as $\vec{\xi}_k$.
One of the main computation in \cite{Ri2} establishes that in conformal coordinates $\Psi$ from $D^2$
into $B_{r^i}(x^i_\infty)$ one has
\be
\label{III.48}
\begin{array}{l}
\ds div\lf[\nabla\vec{H}_{\vec{\xi}}-3\pi_{\vec{n}_{\vec{\xi}}}(\nabla \vec{H}_{\vec{\xi}})+\star(\nabla^\perp\vec{n}_{\vec{\xi}}\wedge\vec{H}_{\vec{\xi}})\rg]\\[5mm]
\ds=-2\,e^{2\la}\lf[\Delta_\perp\vec{H}_{\vec{\xi}}+\ti{A}(\vec{H}_{\vec{\xi}})-2|\vec{H}_{\vec{\xi}}|^2\vec{H}_{\vec{\xi}}\rg]\quad .
\end{array}
\ee
where $e^\la$ is the conformal factor of the immersion in conformal coordinates w.r.t. these coordinates $(x,y)$ and $\pi_{\vec{n}_{\vec{\xi}}}$
is the orthogonal projection onto the normal space to the immersion $\vec{\xi}$. One observe that 
$e^{2\la}\ [dx^2+dy^2]=dvol_{\vec{\xi}}$. Hence (\ref{III.48}) implies that for any function $\vec{f}$ in $C^\infty_0(D^2,{\R}^m)$, for any
smooth immersion $\vec{\xi}$ from $B_{r^i}(x^i_\infty)$ into ${\R}^m$ and for any conformal coordinates $\Psi$ one has
\be
\label{III.49}
\begin{array}{l}
\ds\int_{D^2}\nabla \vec{f}\cdot\lf[\nabla\vec{H}_{\vec{\xi}}-3\pi_{\vec{n}_{\vec{\xi}}}(\nabla \vec{H}_{\vec{\xi}})+\star(\nabla^\perp\vec{n}_{\vec{\xi}}\wedge\vec{H}_{\vec{\xi}})\rg]\ dx\,dy\\[5mm]
\ds=2\int_{B_{r^i}(x^i_\infty)} \vec{f}\circ\Psi^{-1}\cdot \lf[\Delta_\perp\vec{H}_{\vec{\xi}}+\ti{A}(\vec{H}_{\vec{\xi}})-2|\vec{H}_{\vec{\xi}}|^2\vec{H}_{\vec{\xi}}\rg]\ 
dvol_{\gamma}
\end{array}
\ee
where $\gamma=\xi^\ast g_{{\R}^m}$. As observed in \cite{Ri2}, The projection $\pi_{\vec{n}_{\vec{\xi}}}$ can be expressed using the Gauss
$m-2$-vector $\vec{n}_{\vec{\xi}}$ and the interior multiplication $\res$ between multivectors
\[
\pi_{\vec{n}_{\vec{\xi}}}(\vec{v}):=\vec{n}_{\vec{\xi}}\res(\vec{n}_{\vec{\xi}}\res\vec{v})\quad.
\]
Hence we have in particular 
\be
\label{III.50}
\pi_{\vec{n}_{\vec{\xi}}}(\nabla \vec{H}_{\vec{\xi}})=(\nabla\vec{n}_{\vec{\xi}})\res(\vec{n}_{\vec{\xi}}\res\vec{H}_{\vec{\xi}})+ \vec{n}_{\vec{\xi}}\res((\nabla\vec{n}_{\vec{\xi}})\res\vec{H}_{\vec{\xi}})
+\vec{n}_{\vec{\xi}}\res(\vec{n}_{\vec{\xi}}\res\nabla\vec{H}_{\vec{\xi}})
\ee
Taking now $\varphi\in C^\infty_0(D^2)$ such that $\int_{D^2}\varphi=1$ and denote $\varphi_\ep(x):=\ep^{-2}\varphi(\ep^{-1}\, x)$. We also denote
$\vec{\xi}_{\ep,k}:=\varphi_\ep\star\vec{\xi}_k$ and $$\vec{\xi}_{\ep,k}^t:=\varphi_\ep\star\vec{\xi}_k+t\ \chi_{(B_{r^i}(x^i_\infty))}\ \vec{w}\circ\zeta_k^{-1}\circ(\Psi^i)^{-1}$$
We have that $\vec{\xi}_{\ep,k}\longrightarrow\vec{\xi}$ strongly in $W^{2,2}(B_{r^i}(x^i_\infty),{\R}^m)$ and for $\ep$ small
enough $\|\log|\nabla \vec{\xi}_{\ep,k}|\|_\infty$ remains uniformly bounded. Hence we deduce that, as $\ep$ goes to zero
\be
\label{III.51}
\vec{n}_{\vec{\xi}_{\ep,k}}\longrightarrow\vec{n}_{\vec{\xi}_{k}}\quad\mbox{ strongly in }W^{1,2}(B_{r^i}(x^i_\infty),{\R}^m)\quad,
\ee
 and 
 \be
 \label{III.52}
\vec{H}_{\vec{\xi}_{\ep,k}}\longrightarrow\vec{H}_{\vec{\xi}_{k}}\quad\mbox{ strongly in }L^2(B_{r^i}(x^i_\infty),{\R}^m)\quad.
\ee
Hence, combining (\ref{III.50}), (\ref{III.51}) and (\ref{III.52}) we obtain that, for any $\vec{f}\in W^{1,\infty}\cap W^{2,2}(D^2,{\R}^m)$, as $\ep$ goes to zero 
\be
\label{III.53}
\begin{array}{c}
\ds\int_{D^2}\nabla \vec{f}\cdot\lf[\nabla\vec{H}_{\vec{\xi}_{\ep,k}}-3\pi_{\vec{n}_{\vec{\xi}_{\ep,k}}}(\nabla \vec{H}_{\vec{\xi}_{\ep,k}})+\star(\nabla^\perp\vec{n}_{\vec{\xi}_{\ep,k}}\wedge\vec{H}_{\vec{\xi}_{\ep,k}})\rg]\ dx\,dy\\[5mm]
 \ds \quad\quad\downarrow\\[5mm]
\ds\int_{D^2}\nabla \vec{f}\cdot\lf[\nabla\vec{H}_{\vec{\xi}_{k}}-3\pi_{\vec{n}_{\vec{\xi}_{k}}}(\nabla \vec{H}_{\vec{\xi}_{k}})+\star(\nabla^\perp\vec{n}_{\vec{\xi}_{k}}\wedge\vec{H}_{\vec{\xi}_{k}})\rg]\ dx\,dy\quad.
\end{array}
\ee
One verifies easily moreover that $G(\vec{\xi}_{\ep,k},\vec{w},t)\longrightarrow G(\vec{\xi}_{k},\vec{w},t)$ in $L^1$. Hence applying
(\ref{III.45}) and (\ref{III.48}) to $\vec{\xi}:=\vec{\xi}_{\ep,k}$ and passing to the limit as $\ep$ goes to zero using again (\ref{III.51})
\be
\label{III.54}
\begin{array}{l}
W(\vec{\xi}_k)-W(\vec{\xi}_k^t)\\[5mm]
\ds=-2^{-1}\ t\ \int_{D^2}\nabla \vec{w}\cdot\lf[\nabla\vec{H}_{\vec{\xi}_{k}}-3\pi_{\vec{n}_{\vec{\xi}_{k}}}(\nabla \vec{H}_{\vec{\xi}_{k}})+\star(\nabla^\perp\vec{n}_{\vec{\xi}_{k}}\wedge\vec{H}_{\vec{\xi}_{k}})\rg]\ dx\,dy\\[5mm]
\ds\quad-8^{-1}\ t^2\ \int_{B_{r^i}(x^i_\infty)}G(\vec{\xi}_{k},\vec{w},t)
\end{array}
\ee
Thus, combining (\ref{III.44}) and (\ref{III.54}), dividing by $|t|$ taking respectively the limit as $t\rightarrow 0^+$ and $t\rightarrow 0^-$ one obtains
for any $\vec{w}$ in $W^{1,\infty}\cap W^{2,2}(D^2,{\R}^m)$ supported in a strict open subset to $D^2$
\be
\label{III.55}
\begin{array}{l}
\ds\lf| \int_{D^2}\nabla \vec{w}\cdot\lf[\nabla\vec{H}_{\vec{\xi}_{k}}-3\pi_{\vec{n}_{\vec{\xi}_{k}}}(\nabla \vec{H}_{\vec{\xi}_{k}})+\star(\nabla^\perp\vec{n}_{\vec{\xi}_{k}}\wedge\vec{H}_{\vec{\xi}_{k}})\rg]\ dx\,dy\rg|\\[5mm]
\ds\quad\quad\le C\ 2^{-k/2}\ \lf[\|\nabla \vec{w}\|_{\infty}+\|\nabla^2\vec{w}\|_{2}\rg]
\end{array}
\ee
This implies that
\be
\label{III.56}
\begin{array}{c}
\ds div\lf[\nabla\vec{H}_{\vec{\xi}_{k}}-3\pi_{\vec{n}_{\vec{\xi}_{k}}}(\nabla \vec{H}_{\vec{\xi}_{k}})+\star(\nabla^\perp\vec{n}_{\vec{\xi}_{k}}\wedge\vec{H}_{\vec{\xi}_{k}})\rg]\\[5mm]
\ds\longrightarrow 0\quad\mbox{ in }
\quad(W^{1,\infty}\cap W^{2,2})^\ast\quad.
\end{array}
\ee
Using theorem II.1 and theorem II.2 of \cite{BR} we deduce that $\vec{\xi}_\infty$ is
conformal willmore on $B_{r^i}(x^i_\infty)$ but since $\vec{\xi}_\infty$ is conformal from $(B_{r^i}(x^i_\infty),J_\infty)$ into ${\R}^m$, using 
again theorem II.2 of \cite{BR} we obtain that $\vec{\xi}_\infty$ is analytic on $B_{r^i}(x^i_\infty)$. This holds for any $i$ in $I_1$
and hence we have proved lemma~\ref{lm-III.2}.\hfill $\Box$

\medskip

\noindent{\bf Proof of theorem~\ref{th-III.1} continued.}

The goal now is to extend $\vec{\xi}_\infty$ as a smooth embedding through the points $a_j$.

\medskip

 Let $a_j$ be such a point. Let $\Psi_j$ be 
a positive conformal diffeomorphism from $D^2$, equipped with the canonical complex structure, into a neighborhood $U_j$ of $a_j$ in $(\Sigma,J_\infty)$
and such that $\Psi_j(0)=a_j$
We keep denoting $\vec{\xi}_\infty$ the composition $\vec{\xi}_\infty\circ\Psi_j$.

\medskip

 $\vec{\xi}_\infty$ is conformal from $D^2\setminus\{0\}$ into ${\R}^m$
and we have that ${\mathcal H}^2(\vec{\xi}_\infty(D^2\setminus\{0\}))<+\infty$ moreover $\vec{\xi}_\infty(D^2\setminus\{0\})\subset B_R(0)$. Hence $\vec{\xi}_\infty$ is
in $L^\infty\cap W^{1,2}(D^2\setminus\{0\},{\R}^m)$.
Since the 2-capacity of a point in 2 dimension is zero we deduce that $$\vec{\xi}_\infty\in L^\infty\cap W^{1,2}(D^2,{\R}^m)\quad .$$

Similarly,
$\vec{n}_{\vec{\xi}_\infty}$ realizes a map in $W^{1,2}(D^2\setminus\{0\},Gr_{m-2}({\R}^m))$. For the same reason as before,
$\vec{n}_{\vec{\xi}_\infty}$ extends to a map in $W^{1,2}(D^2,Gr_{m-2}({\R}^m))$.We now use the lemma~\ref{lm-A.4} which is already implicitly present in \cite{Hub}, \cite{MS} and \cite{Hel} but for which we thought that it
could have been useful for the reader to have the details of a proof of it presented in the appendix. We can deduce from this lemma that $\vec{\xi}_\infty$ extends to a Lipshitz map
through $0$ and that there exists an integer $m$ such that
\be
\label{III.57}
(C-o(1))\ |z|^{m-1}\le\lf|\frac{\p\vec{\xi}_\infty}{\p z}\rg|\le (C+o(1))\ |z|^{m-1}\quad.
\ee
We claim that $m=1$. Because of this estimate, for any $\delta>0$ there exists $r_\delta>0$ such that, for any $r<r_\delta$,
 $\vec{\xi}_\infty(B_r(0))\subset B_{\rho}(\vec{\xi}_\infty(0))$ and $|\p_x\vec{\xi}_\infty|=|\p_x\vec{\xi}_\infty|=e^\la_\infty\ge C(1-\delta)/\sqrt{2}\ |z|^{m-1}$.
 where $\rho=C(\sqrt{2})^{-1}\ m^{-1}\, (1+\delta)\ r^m$. We have then that the mass of $\vec{\xi}_\infty(\Sigma)$ present in $B_\rho(\vec{\xi}_\infty(0))$ can be estimated from below
 as follows
 \be
 \label{III.58}
 \begin{array}{rl}
 \ds M\lf(\vec{\xi}_\infty(\Sigma)\res B_\rho(\vec{\xi}_\infty(0))\rg)&\ds\ge C^2\ \frac{(1-\delta)^2}{2}\ \int_{B_r(0)}|z|^{2m-2}\\[5mm]
  &\ds\ge \frac{\pi}{m}\ C^2\ \frac{(1-\delta)^2}{2}\ r^{2m}\\[5mm]
  &\ds\ge m\ \lf(\frac{1-\delta}{1+\delta}\rg)^2\ \pi\ \rho^2
 \end{array}
 \ee
This implies that the lower 2-density $\theta_\ast^2((\vec{\xi}_\infty)_\ast[\Sigma],\vec{\xi}_\infty(0))$ of $(\vec{\xi}_\infty)_\ast[\Sigma]$ at $\vec{\xi}_\infty(0)$ is larger
or equal to $m$. The Li-Yau inequality (see \cite{LY}) which also holds for varifolds with weak $L^2-$bounded mean curvature which are smooth outside
one point (as proved in \cite{KS2}) implies
\be
\label{III.59}
m\le\theta_\ast^2((\vec{\xi}_\infty)_\ast[\Sigma],\vec{\xi}_\infty(0))\le\frac{W(\vec{\xi}_\infty(\Sigma))}{4\pi}\quad.
\ee
Because of the lower semi-continuity of $W$ and the assumption that $W(\vec{\xi}_k(\Sigma))<8\pi-\delta$ for some $\delta>0$, we deduce that $m=1$.

\medskip 

We have then proved that $\vec{\xi}_\infty$ is a $W^{2,2}-lipschitz$ immersion, that is an element from ${\mathcal E}_\Sigma$, which is smooth
on $\Sigma\setminus\{a_1\cdots a_N\}$ and satisfying
\[
W(\vec{\xi}_\infty)\le \inf_{\vec{\Phi}\in{\mathcal E}_\Sigma}W(\vec{\Phi})\quad.
\]
Because of the minimality of $\vec{\xi}$, for any $i$ and any $\vec{w}\in C^\infty_0(B_r(a_i))$ for some $r$, we have that, for $t$ small enough
\be
\label{III.60}
W(\vec{\xi}_\infty)\le W(\vec{\xi}_\infty+t\vec{w})\quad\quad.
\ee
Arguing like above we have the asymptotic expansion 
\be
\label{III.61}
\begin{array}{l}
W(\vec{\xi}_\infty)-W(\vec{\xi}_\infty+t\vec{w})\\[5mm]
\ds=-2^{-1}\ t\ \int_{D^2}\nabla \vec{w}\cdot\lf[\nabla\vec{H}_{\vec{\xi}_{\infty}}-3\pi_{\vec{n}_{\vec{\xi}_{\infty}}}(\nabla \vec{H}_{\vec{\xi}_{\infty}})+\star(\nabla^\perp\vec{n}_{\vec{\xi}_{\infty}}\wedge\vec{H}_{\vec{\xi}_{\infty}})\rg]\ dx\,dy\\[5mm]
\ds\quad-8^{-1}\ t^2\ \int_{B_{r}(a_i)}G(\vec{\xi}_{\infty},\vec{w},t)
\end{array}
\ee
where we are using some holomorphic chart on $B_r(a_i)$ which identifies $B_r(a_i)$ with $D^2$ and where $\int_{B_{r}(a_i)}G(\vec{\xi}_{\infty},\vec{w},t)$ is uniformly bounded w.r.t. $t$ as before.
Combining (\ref{III.60}) and (\ref{III.61}) we deduce that $\vec{\xi}_\infty$ realizes a weak Willmore $W^{2,2}\cap W^{1,\infty}$ immersion and
from \cite{Ri2} we deduce that $\vec{\xi}_\infty$ is an analytic immersion. Since $W(\vec{\xi}_\infty)<8\pi$ we deduce
that $\vec{\xi}_\infty$ realizes an embedding which concludes the proof of theorem~\ref{th-III.1}. \hfill $\Box$

\section{Existence of Minimizers of the Willmore Energy in a Conformal Class.}
\reset

\subsection{The completeness of the metric space of $W^{2,2}$ lipschitz immersions of $\Sigma$ in a given conformal class.}

We assume in this section that $\Sigma$ is a connected closed smooth two dimensional manifold of genus larger or equal
to 1.
Let $c$ be a conformal class $\Sigma$ which is represented by a smooth complex structure $J$ on $\Sigma$.  Denote $g$
be an arbitrary smooth metric on $g$ that we can choose to be compatible with $J$.
We introduce the subspace space of ${\mathcal E}_\Sigma$ of lipschitz immersions realizing a complex structure
equivalent to $J$ :
\[
{\mathcal E}_\Sigma^c:=\lf\{
\begin{array}{c}
\vec{\Phi}\in{\mathcal E}_\Sigma\quad\mbox{ s.t. } \exists\ \Psi\in W^{2,2}(\Sigma,\Sigma)\quad\mbox{ s.t. }\\[5mm]
\Psi\mbox{ is a bilipschitz diffeomorphism }\\[5mm]
\vec{\Phi}\circ\Psi\quad :\quad(\Sigma,J) \longrightarrow {\R}^m\quad\mbox{ is conformal }
\end{array}
\rg\}
\]
We are now proving the following proposition
\begin{Prop}
\label{pr-IV.1}
The metric space $({\mathcal E}_\Sigma^c,d^J)$ is complete.\hfill $\Box$
\end{Prop}
{\bf Proof of proposition~\ref{pr-IV.1}.}
Let $\vec{\Phi}_k$ be a Cauchy sequence for $d^J$. From proposition~\ref{pr-II.1} there exists a limit $\vec{\Phi}_\infty$ in ${\mathcal E}_\Sigma$. 
Denote by $\Psi_k$ a lipschitz diffeomorphism such that $\vec{\Phi}_k\circ\Psi_k$ is conformal. Denote by $g_k:=\vec{\Phi}_k^\ast g_{{\R}^m}$
and $J_k$ the associated complex structure.
$\Psi_k$ realizes then a bilipschitz conformal diffeomorphism between $(\Sigma,J)$ and $(\Sigma,J_k)$.
Because of the $d^J$ convergence $J_k$ converges in $L^\infty\cap W^{1,2}$ norm to $J_\infty$ the complex structure associated to $g_\infty:=\vec{\Phi}_\infty^\ast g_{{\R}^m}$.
Denote by $(U^i)_{i\in I}$ a finite covering by balls of $\Sigma$ chosen in such a way that $\int_{U_i}|d\vec{n}_{\vec{\Phi}_\infty}|^2\ dvol_{g_\infty}<4\pi/3$. Denote
by $\zeta^i_\infty\ D^2\ :\ \rightarrow (U^i,J_\infty)$ the conformal parametrization given by lemma 5.1.4 of \cite{Hel} combined with the moving frame technic of the proof of theorem 5.4.3
that we exposed also in the proofs of lemma~\ref{lm-A.3} and lemma~\ref{lm-A.4} below. For each $i$ in $I$ we use lemma~\ref{lm-A.3} in order to construct 
$\varphi_k$ : $D^2\rightarrow D^2$ such that $\zeta^i_k:=\zeta^i_\infty\circ\varphi_k\ D^2\ :\ \rightarrow\ (U_i,J_k)$ is conformal and $\zeta^i_k$ is uniformly bounded in $W^{2,2}$
and $\log|\nabla \zeta^i_k|$ is also uniformly bounded in $L^\infty$ (by taking possibly $U^i$ a bit smaller but still realizing a covering of $\Sigma$).
Denote by $f_k^i(z):=(\Psi_k^i)^{-1}\circ\zeta^i_k$ the maps from $D^2$ into $(\Sigma,J)$. These sequences realize sequences of conformal maps which are harmonic if one equips
$(\Sigma,J)$ with a corresponding constant scalar curvature metric $h$ and denote the corresponding volume form. Observe that since $\Psi_k$ is a conformal diffeomorphism
one has
\[
\begin{array}{l}
\ds\int_\Sigma\om=\int_{\Sigma}(\Psi_k^{-1})^\ast\om=\frac{1}{2}\int_{\Sigma}|d\Psi_k^{-1}|_{h,g_k}^2\ dvol_{g_k}\\[5mm]
\ds\quad\quad\ge\frac{1}{card\, I}\sum_{i\in I}\int_{U_i}|\nabla f_k^i|^2_h\ dx_1\, dx_2\quad.
\end{array}
\]
Hence $\Psi_k^{-1}$ is a uniformly bounded sequence in $W^{1,2}(\Sigma,\Sigma)$ (the metric $g_k$ remains comparable to an arbitrary smooth fixed 
metric on $\Sigma$ because of the $d^J$ convergence) and the $f_k^i$ are uniformly bounded sequences of conformal maps in $W^{1,2}(D^2,\Sigma)$. Hence the $f_k^i$ are uniformly bounded energy harmonic maps. Since the constant  scalar curvature of the metric $h$ is non-positive (genus($\Sigma$)>0)
 the sequences converge strongly in $C^l-$norm in the interiors of $U_i$ (see for instance \cite{Jo}). Since the $\zeta^i_k$ are uniformly bounded
 in $W^{2,2}$, since $\log|\nabla \zeta^i_k|$ is also uniformly bounded in $L^\infty$ and since $\|Dis^{J_\infty}(\zeta^i_k)\|_\infty\longrightarrow 0$,
 we deduce that $\Psi_k^{-1}$ converges to a bilipschitz diffeomorphism $\Psi_\infty^{-1}$ which is conformal between $(\Sigma,J_\infty)$ and
 $(\Sigma,J)$. This implies that the $d^J$ limit $\vec{\Phi}_\infty$ is in ${\mathcal E}_\Sigma^c$ and this concludes the proof of proposition~\ref{pr-IV.1}.\hfill $\Box$
 
\medskip

\subsection{Minimizing Willmore energy in a conformal class.}

In this section we prove the following theorem

\begin{Th}
\label{th-IV.2}
Let $\Sigma$ be a closed surface let $c$ be a conformal class on $\Sigma$ and $m$ an integer larger or equal to 3. 
Assume that
\[
\inf_{\vec{\Phi}\in {\mathcal E}_\Sigma^c}W(\vec{\Phi})\le 8\pi
\]
Then the infimum is achieved by either
\begin{itemize}
\item[i)]  a $C^\infty$ Conformally Willmore embedding of $\Sigma$ into ${\R}^m$
\item[ii)] or a global isothermic embedding of $(\Sigma,c)$.
\end{itemize}
If $\inf_{\vec{\Phi}\in {\mathcal E}_\Sigma^c}W(\vec{\Phi})> 8\pi$ the results is the same modulo the possible existence of isolated branched points.
\hfill $\Box$
\end{Th}
\noindent{\bf Proof of theorem~\ref{th-IV.2}.}
Let $\vec{\Phi}_k$ be a minimizing sequence of $W$ in ${\mathcal E}_\Sigma^c$. Applying the {\it 3-points renormalization lemma} and arguing exactly like in the beginning of the 
proof of theorem~\ref{th-III.1} we can assume that
\begin{itemize}
\item[i)]
\be
\label{IV.11}
\vec{\Phi}_k\quad\quad\mbox{ is conformal from } (\Sigma,J)\quad\mbox{ into }{\R}^m\quad.
\ee
\item[ii)] There exists finitely many points $a_1\cdots a_N$ in $\Sigma$ and a fixed finite covering $({B}_{\rho_\infty^i}(x_\infty^i))_{i\in I_1}$
of $\Sigma\setminus\{a_1\cdots a_N\}$ such that for any $i\in I_1$, $0<\rho<\rho_\infty$ and $k$ large enough
\be
\label{IV.12}
\int_{B_\rho(x_\infty^i)}|d\vec{n}_{\vec{\Phi}_k}|^2_{g_k}\ dvol_{g_k}<8\pi/3\quad .
\ee
where $g_k:=\vec{\Phi}_k^\ast g_{{\R}^m}$.
\item[iii)] There exists a positive real $R>0$ such that
\be
\label{IV.13}
\vec{\Phi}_k(\Sigma)\subset B_R(0)\quad.
\ee
\item[iv)] There exists a constant $C>0$ such that
\be
\label{IV.14}
{\mathcal H}^2(\vec{\Phi}_k(\Sigma))\le C\quad .
\ee
\item[v)] There exist a positive  real number $r>0$, independent of $k$ and three distinct points $P_1$, $P_2$ and $P_3$, independent of $k$ too, in the interior of one ball $B_{\rho_\infty^i}(x_\infty^i)$
such that
\be
\label{IV.15}
\forall\ i\ne j\quad\quad\quad|\vec{\Phi}_k(P_i)-\vec{\Phi}_k(P_j)|\ge r>0\quad\quad.
\ee
\end{itemize}
The following proposition is a direct application of Ekeland's Variational Principle since $({\mathcal E}_\Sigma^c,d^J)$ is a complete metric space as we showed in the previous subsection.
\begin{Prop}
\label{pr-IV.3}
Let $J$ be an arbitrary smooth complex structure on $\Sigma$ and  $c$ be the conformal class of $(\Sigma,J)$. Let $\vec{\Phi}_k$ be a minimizing sequence for $W$
in ${\mathcal E}_\Sigma^c$ such that
\[
W(\vec{\Phi}_k)\le\inf_{\vec{\Phi}\in{\mathcal E}^c_\Sigma}W(\vec{\Phi})+2^{-k}\quad,
\]
then there exists $\vec{\xi}_k\in {\mathcal E}^c_\Sigma$ such that 
\begin{itemize}
\item[i)]
$\vec{\xi}_k$ minimizes in ${\mathcal E}^c_\Sigma$ the following functional
\be
\label{III.av1}
W(\vec{\xi}_k)=\inf_{\vec{\Phi}\in{\mathcal E}^c_\Sigma}W(\vec{\Phi})+2^{-k/2}\ d^J(\vec{\Phi},\vec{\xi}_k)\quad ,
\ee
\item[ii)]
\be
\label{III.av2}
W(\vec{\xi}_k)\le W(\vec{\Phi}_k)\quad,
\ee
\item[iii)]
\be
\label{III.av3}
d^J(\vec{\xi}_k,\vec{\Phi}_k)\le 2^{-k/2}\quad.
\ee
\end{itemize}
 \hfill$\Box$
\end{Prop}
As in the previous section we prove that $d\vec{\xi}_k$ converges strongly in $L^p_{loc}({\Sigma}\setminus\{a_1\cdots a_N\})$ to a limiting $W^{2,2}$ immersion $\vec{\xi}_\infty$
of ${\Sigma}\setminus\{a_1\cdots a_N\}$ and moreover we have that
\be
\label{IV.16}
\int_{\Sigma\setminus\{a_1\cdots a_N\}}|d\vec{n}_{\vec{\xi}_\infty}|^2_{g_\infty}\ dvol_{g_\infty}\le \inf_{\vec{\Phi}\in{\mathcal E}^c_\Sigma} 4W(\vec{\Phi})-4\pi\chi(\Sigma)\quad.
\ee
We claim now that
\begin{Lm}
\label{lm-IV.5} Under the previous notations we have that either
\be
\label{IV.17}
\vec{\xi}_\infty\mbox{ is a {\it Conformal Willmore immersion} of }\Sigma\setminus\{a_1\cdots a_N\}
\ee
and hence is analytic on $\Sigma\setminus\{a_1\cdots a_N\}$ (see \cite{BR}) or
\[
\vec{\xi}_\infty\mbox{ is an isothermic immersion of }(\Sigma,J)\quad.
\]
\hfill $\Box$
\end{Lm}

\medskip

\noindent{\bf Proof of the lemma~\ref{lm-IV.5}.}

\noindent{\bf First case : there exists a subsequence $\vec{\xi}_k$ which is not made of isothermic surfaces.}

\medskip

Let $\Psi^i$ be a local conformal chart on $B_{r^i}(x^i_\infty)$ for the complex structure $J$. We have that $\vec{\xi}_k\circ\Psi^i$ is conformal on $D^2$.
From lemma~\ref{lm-A.5} we know that a perturbation $\vec{w}\in W^{1,\infty}\cap W^{2,2}$ of $\vec{\xi}_k$ keeps infinitesimally the conformal class $c$ if and only if for any
holomorphic quadratic differential $q$ of $(\Sigma,J)$ which is an holomorphic section of $K \otimes K$, where $K$ is the canonical bundle $T^{(0,1)}\Sigma$ of (1-0)-forms over $(\Sigma,J)$,
one has
\be
\label{IV.18}
(\p_z\vec{w}\cdot\p_z\vec{\xi}_k\ dz\otimes dz,q)_{WP}=0
\ee
where $(\cdot,\cdot)_{WP}$ is the Weil-Petersson Hermitian product given locally (assuming $\vec{w}$ is supported in a ball $B_{r^i}(x^i_\infty)$ on which we have holomorphic chart
given by $\Psi^i$ that we simply denote by $z$), writing $q=f(z)\ dz\otimes dz$,
\[
(\p_z\vec{w}\cdot\p_z\vec{\xi}_k\ dz\otimes dz,q)_{WP}:=\frac{i}{2}\int_{D^2}e^{-2\la_k}\ \p_z\vec{w}\cdot\p_z\vec{\xi}_k\ \ \ov{f(z)}\ dz\wedge d\ov{z}\quad.
\]
where $e^{2\la_k}=|\p_x\vec{\xi}_k|^2=|\p_y\vec{\xi}_k|^2$.

Consider hence $\vec{w}\in W^{1,\infty}\cap W^{2,2}(D^2,{\R}^m)$ supported in the interior of $D^2$, satisfying (\ref{IV.18})
and such that
\[
\|\nabla\vec{w}\|_{L^\infty(D^2)}+\|\nabla^2\vec{w}\|_{L^2(D^2)}\le 1\quad .
\]
Using Lemma~\ref{lm-A.5} and the implicit function theorem there exists a family $\vec{\xi}_k^t$ in ${\mathcal E}_\Sigma^c$ such that
\[
\vec{\xi}_k^t:=\vec{\xi}_k+t\ \chi_{(B_{r^i}(x^i_\infty))}\ \vec{w}\circ(\Psi^i)^{-1}+o(t)\quad.
\]
From now on we shall omit to write explicitly the composition with $(\Psi^i)^{-1}$ and write simply $\vec{w}$
instead of $\vec{w}\circ(\Psi^i)^{-1}$.

Arguing exactly like in the previous section this implies that there exists a constant $C>0$  such that, for all $\vec{w}$ supported in a strict
open subset of $D^2$ and
satisfying (\ref{IV.18}) one has 
\be
\label{IV.19}
\begin{array}{l}
\ds\lf| \int_{D^2}\nabla \vec{w}\cdot\lf[\nabla\vec{H}_{{k}}-3\pi_{\vec{n}_{{k}}}(\nabla \vec{H}_{{k}})+\star(\nabla^\perp\vec{n}_{{k}}\wedge\vec{H}_{{k}})\rg]\ dx\,dy\rg|\\[5mm]
\ds\quad\quad\le C\ 2^{-k/2}\ \lf[\|\nabla \vec{w}\|_{\infty}+\|\nabla^2\vec{w}\|_{2}\rg]
\end{array}
\ee
  where $\vec{H}_k:=\vec{H}_{\vec{\xi}_{k}}$ and $\vec{n}_k:=\vec{n}_{\vec{\xi}_{k}}$. Using the notations and computations in appendix~\ref{lm-A.5} the constraints (\ref{IV.18}) on $\vec{w}$ becomes
\be
\label{IV.20}
\forall j=1\cdots Q\quad\quad\int_{D^2} f^j(z)\ \vec{{H}}_{0,k}\cdot \vec{w}\ \frac{i}{2}\ dz\wedge d\ov{z}= 0\quad\quad.
\ee
where we recall that $f^j(z)\ dz\otimes dz$ is the expression in the $\Psi^i$ conformal chart of the different element $q^j$ of a \underbar{fixed} basis
of the $Q-$dimensional complex space of holomorphic quadratic differentials $Q(J)$ of $(\Sigma, J)$ and where 
$\vec{H}_{0,k}$ is the Weingarten operator associated to the immersion $\vec{\xi}_k$. .
Combining (\ref{IV.19}) and (\ref{IV.20}) we then obtain the existence of a sequence $\mu_k=(\mu^j_k)_{j=1\cdots N}\in{\C}^Q$ such that
\be
\label{IV.21}
\begin{array}{c}
\ds div\lf[\nabla\vec{H}_{{k}}-3\pi_{\vec{n}_{{k}}}(\nabla \vec{H}_{{k}})+\star(\nabla^\perp\vec{n}_{{k}}\wedge\vec{H}_{{k}})\rg]+
\Im\lf[f_k(z)\ \vec{H}_{0,k}\rg]\\[5mm]
\ds\longrightarrow 0\quad\mbox{ in }
\quad(W^{1,\infty}\cap W^{2,2})^\ast\quad.
\end{array}
\ee
where 
\[
f_k(z):=\sum_{j=1}^Q\ \mu^j_k\ f^j(z)\quad.
\]
Using the computations in \cite{BR} section III.2.2  we have
\be
\label{IV.23}
\begin{array}{c}
\ds e^{2\la_k}\lf[\Delta_\perp\vec{H}_{{k}}\;+\;2\,\Re\big((\overline{\vec{H}_{0,k}}\cdot\vec{H}_{{k}})\,\vec{H}_{0,k}\big)\rg]-\Im\big(f_k(z)\,\vec{H}_{k,0}\big)\:.\\[5mm]
\ds\longrightarrow 0\quad\mbox{ in }
\quad(W^{1,\infty}\cap W^{2,2})^\ast\quad.
\end{array}
\ee
The covariant laplacian in conformal coordinates is given by
\be
\label{IV.24}
e^{2\la_k}\,\Delta_\perp\vec{H}_{{k}}\:=\:\pi_{\vec{n}_{k}}\Big({div}\,\pi_{\vec{n}_{k}}\big(\nabla\vec{H}_{{k}}\big)\Big)\:\equiv\:4\,\Im\Big[i\,
\pi_{\vec{n}_{k}}
\,\p_{\ov{z}}\pi_{\vec{n}_{k}}\,\p_z \vec{H}_{{k}}   \Big]\quad.
\ee
Combining (\ref{IV.23}) and (\ref{IV.24}) we obtain that
\be
\label{IV.24a}
\begin{array}{c}
\ds\,\Im\Big[4i\,\pi_{\vec{n}_{k}}\,\p_{\ov{z}}\pi_{\vec{n}_{k}}\,\p_z \vec{H}_{{k}} +2i\  e^{2\la_k}\ (\overline{\vec{H}_{0,k}}\cdot\vec{H}_{{k}})\,\vec{H}_{0,k}-f_k(z)\,\vec{H}_{k,0} \Big]\\[5mm]
\ds\longrightarrow 0\quad\mbox{ in }
\quad(W^{1,\infty}\cap W^{2,2})^\ast\quad.
\end{array}
\ee
It is convenient to introduce $A_k\in {\C}$ given by
\be
\label{IV.25}
A_k\;=\;{e}^{-\la_k}f_k(z)\,-\,2\,i\,{e}^{\la_k}\,\overline{\vec{H}_{0,k}}\cdot\vec{H}_{{k}} \quad.
\ee
With this notation we have in particular
\be
\label{IV.25a}
\p_{\ov{z}}(e^{\la_k}\ A_k)=-\,2\,i\,\p_{\ov{z}}\lf({e}^{2\la_k}\,\overline{\vec{H}_{0,k}}\cdot\vec{H}_{{k}}\rg)\quad.
\ee
Using the general equation $\p_{\ov{z}}(e^{-\la_k}e_{\ov{z}})=2^{-1}\vec{H}_0$ (see again \cite{BR} section III.2.2), we have
\be
\label{IV.25b}
\p_{\ov{z}}(A_k\ e_{\ov{z}})=-\,2\,i\,e^{-\la_k}\ \p_{\ov{z}}\lf({e}^{2\la_k}\,\overline{\vec{H}_{0,k}}\cdot\vec{H}_{{k}}\rg)\,e_{\ov{z}}+\frac{e^{\la_k}}{2}\ A_k\ \vec{H}_{0,k}\quad.
\ee
We  recall at this stage the Codazzi-Mainardi equation \footnote{See \cite{BR} lemma A.3 for a proof.}
\be
\label{IV.26}
{e}^{-2\la_k}\,\p_{\ov{z}}\big({e}^{2\la_k}\,\overline{\vec{H}_{0,k}}\cdot\vec{H}_{{k}}\big)\;=\;\vec{H}_{{k}}\cdot\p_z\vec{H}_{{k}}\,+\,\overline{\vec{H}_{0,k}}\cdot\p_{\ov{z}}\vec{H}_{{k}}\:.
\ee
Combining (\ref{IV.25}), (\ref{IV.25b}) and (\ref{IV.26}) we obtain that
\be
\label{IV.26b}
\begin{array}{l}
\ds\p_{\ov{z}}(A_k\ \vec{e}_{\ov{z}})=-\,2\,i\,e^{\la_k}\,\lf[\vec{H}_{{k}}\cdot\p_z\vec{H}_{{k}}\,+\,\overline{\vec{H}_{0,k}}\cdot\p_{\ov{z}}\vec{H}_{{k}}\rg]\,\vec{e}_{\ov{z}}+\frac{e^{\la_k}}{2}\,A_k\ \vec{H}_{0,k}\\[5mm]
\ds\quad\quad\quad =-\,2\,i\,e^{\la_k}\,\lf[\vec{H}_{{k}}\cdot\p_z\vec{H}_{{k}}\,+\,\overline{\vec{H}_{0,k}}\cdot\p_{\ov{z}}\vec{H}_{{k}}\rg]\,\vec{e}_{\ov{z}}+\frac{1}{2}\,f_k(z)\ \vec{H}_{0,k}\\[5mm]
\ds\quad\quad\quad\quad-i\,e^{2\la_k}\ \lf(\overline{\vec{H}_{0,k}}\cdot\vec{H}_{{k}}\rg)\ \vec{H}_{0,k}
\end{array}
\ee
Another computation in section III.2.2 of \cite{BR} gives
\be
\label{IV.27}
\begin{array}{l}
\ds-2i\p_{\ov{z}}\pi_{\vec{n}_k}\p_z\vec{H}_k=-2i\,\pi_{\vec{n}_k}\p_{\ov{z}}\pi_{\vec{n}_k}\p_z\vec{H}_k\\[5mm]
\ds\quad\quad\quad+2i\,e^{\la_k}\, \lf[(\vec{H}_k\cdot\p_z\vec{H}_k)\,\vec{e}_{\ov{z}}+
(\vec{H}_{0,k}\cdot\p_z\vec{H}_k)\,\vec{e}_z\rg]
\end{array}
\ee
Combining (\ref{IV.26b}) and (\ref{IV.27}) we obtain
\be
\label{IV.28}
\begin{array}{l}
\p_{\ov{z}}(A_k\ \vec{e}_{\ov{z}})-2i\p_{\ov{z}}\pi_{\vec{n}_k}\p_z\vec{H}_k=-2i\,\pi_{\vec{n}_k}\p_{\ov{z}}\pi_{\vec{n}_k}\p_z\vec{H}_k\\[5mm]
\ds\quad\quad\quad-2i\,e^{\la_k}\, \lf[(\overline{\vec{H}_{0,k}}\cdot\p_{\ov{z}}\vec{H}_{{k}})\,\vec{e}_{\ov{z}}-
(\vec{H}_{0,k}\cdot\p_z\vec{H}_k)\,\vec{e}_z\rg]\\[5mm]
\ds\quad\quad+\frac{1}{2}\,f_k(z)\ \vec{H}_{0,k}-i\,e^{2\la_k}\ \lf(\overline{\vec{H}_{0,k}}\cdot\vec{H}_{{k}}\rg)\ \vec{H}_{0,k}
\end{array}
\ee
Observe that
\be
\label{IV.29}
\Im\lf(-2i\,e^{\la_k}\, \lf[(\overline{\vec{H}_{0,k}}\cdot\p_{\ov{z}}\vec{H}_{{k}})\,\vec{e}_{\ov{z}}-
(\vec{H}_{0,k}\cdot\p_z\vec{H}_k)\,\vec{e}_z\rg]\rg)=0
\ee
Combining (\ref{IV.24a}), (\ref{IV.28}) and (\ref{IV.29}) we obtain
\be
\label{IV.30z}
%\begin{array}{c}
\ds\Im\lf(\p_{\ov{z}}\lf[A_k\ \vec{e}_{\ov{z}}-2i\ \pi_{\vec{n}_k}\p_z\vec{H}_k\rg]\rg)\quad
\ds\longrightarrow 0\quad\mbox{ in }
\quad(W^{1,\infty}\cap W^{2,2})^\ast\quad.
%\end{array}
\ee
or in other words
\be
\label{IV.30}
\begin{array}{l}
\ds\Im\lf(\p_{\ov{z}}\lf[e^{-\la_k}\,f_k(z)\,\vec{e}_{\ov{z}}-2\,i\,e^{\la_k}\ \vec{H}_{0,k}\cdot\vec{H}_k\ \vec{e}_z-2\,i\ \pi_{\vec{n}_k}\p_z\vec{H}_k\rg]\rg)\\[5mm]
\quad\quad
\ds\longrightarrow 0\quad\mbox{ in }
\quad(W^{1,\infty}\cap W^{2,2})^\ast\quad.
\end{array}
\ee
Let $\vec{F}_k=\vec{F}^\Re_k+i\vec{F}^\Im_k\in L^{2,\infty}(D^2,{\R}^m\otimes{\C})$ be the unique solution of
\be
\label{IV.31}
\lf\{
\begin{array}{l}
\ds \p_{z}\vec{F}_k=\,e^{-\la_k}\,f_k(z)\,\vec{e}_{\ov{z}}-2\,i\,e^{\la_k}\,\vec{H}_{0,k}\cdot\vec{H}_k\,\vec{e}_z-2\,i\,\pi_{\vec{n}_k}\p_z\vec{H}_k\ \mbox{ in }D^2\\[5mm]
\ds \p_{\nu}\vec{F}^{\Im}_k=0\quad\quad\quad\mbox{ on }\p D^2
\end{array}
\rg.
\ee
Hence combining (\ref{IV.30}) and (\ref{IV.31}) we have
\be
\label{IV.32}
\lf\{
\begin{array}{l}
\ds\Delta\vec{F}_k^\Im\longrightarrow 0\quad\mbox{ in }
\quad(W^{1,\infty}\cap W^{2,2})^\ast\\[5mm]
\ds \p_\nu\vec{F}^{\Im}_k=0\quad\quad\quad\mbox{ on }\p D^2
\end{array}
\rg.
\ee
This implies in particular that
\be
\label{IV.33}
\nabla\vec{F}^\Im_k\longrightarrow 0 \quad\quad\mbox{strongly in }(W^{1,q}(D^2))^\ast\quad\quad\forall\ q>2\quad.
\ee
Let $\vec{Q}^{\C}_k:=\vec{Q}^\Re_k+i\vec{Q}_k^\Im=-4\,\,e^{\la_k}\,\vec{H}_{0,k}\cdot\vec{H}_k\,\vec{e}_z-4\,\,\pi_{\vec{n}_k}\p_z\vec{H}_k\in {\R}^m\otimes{\C}$. It is proved in \cite{Ri2} that for any conformal immersion $\vec{\xi}_k$,  $\vec{Q}:=(\vec{Q}^\Re,\vec{Q}^\Im)$ satisfies
\be
\label{IV.34}
\lf\{
\begin{array}{l}
\vec{Q}_k\cdot\nabla\vec{\xi}_k:= \vec{Q}_k^\Re\cdot\p_x\vec{\xi}_k+\vec{Q}_k^\Im\cdot\p_y\vec{\xi}_k=0\\[5mm]
\vec{Q}_k\wedge\nabla\vec{\xi}_k:= \vec{Q}_k^\Re\wedge\p_x\vec{\xi}_k+\vec{Q}_k^\Im\wedge\p_y\vec{\xi}_k =2\,(-1)^{m}\ \nabla^\perp\lf(\star(\vec{n}_k\res\vec{H}_k)\rg)\res\nabla\vec{\xi}_k\quad.
\end{array}
\rg.
\ee
We rewrite (\ref{IV.31}) in the form
\be
\label{IV.35}
\nabla^\perp\vec{F}^\Re_k+\nabla\vec{F}^\Im_k=\vec{Q}_k+e^{2\la_k}
\lf(
\begin{array}{c}
f^\Im_k(z)\ \p_x\vec{\xi}_k-f^\Re_k(z)\ \p_y\vec{\xi}_k\\[5mm]
-f^\Re_k(z)\ \p_x\vec{\xi}_k-f^\Im_k(z)\ \p_y\vec{\xi}_k
\end{array}
\rg)
\ee
Combining (\ref{IV.34} and (\ref{IV.35}) gives finally
\be
\label{IV.36}
\lf\{
\begin{array}{l}
\lf(\nabla^\perp\vec{F}^\Re_k+\nabla\vec{F}^\Im_k\rg)\cdot\nabla\vec{\xi}_k=0\\[5mm]
\lf(\nabla^\perp\vec{F}^\Re_k+\nabla\vec{F}^\Im_k\rg)\wedge\nabla\vec{\xi}_k:=2\,(-1)^{m}\ \nabla^\perp\lf(\star(\vec{n}_k\res\vec{H}_k)\rg)\res\nabla\vec{\xi}_k\quad.
\end{array}
\rg.
\ee
As in the proof of theorem~\ref{th-III.1} we can extract a subsequence to $\vec{\xi}_k$ that weakly converges to a limiting conformal immersion $\vec{\xi}_\infty$ in $W^{1,\infty}\cap W^{2,2}$ and $\vec{n}_k$ weakly converges in $W^{1,2}$
to the Gauss map $\vec{n}_\infty$ of $\xi_\infty$.
Because of (\ref{IV.32}) and (\ref{IV.33}) we have 
\be
\label{IV.37}
\nabla\vec{F}^\Im_k\cdot\nabla\vec{\xi}_k=div(\nabla\vec{F}^\Im_k\cdot\vec{\xi}_k)-\Delta\vec{F}^\Im_k\cdot\vec{\xi}_k\longrightarrow 0\quad\quad\mbox{ in }{\mathcal D}'(D^2),
\ee
and
\be
\label{IV.38}
\nabla\vec{F}^\Im_k\wedge\nabla\vec{\xi}_k=div(\nabla\vec{F}^\Im_k\wedge\vec{\xi}_k)-\Delta\vec{F}^\Im_k\wedge\vec{\xi}_k\longrightarrow 0\quad\quad\mbox{ in }{\mathcal D}'(D^2),.
\ee
Assume first there exists a subsequence - that we still denote $\vec{\xi}_k$ - such that $|\mu_k|$ is uniformly bounded
and hence a subsequence such that $\mu_k\longrightarrow\mu_\infty=(\mu^j_\infty)_{j=1\cdots Q}$.
This implies that 
\be
\label{IV.22}
f_k(z)\longrightarrow f_\infty(z):=\sum_{j=1}^Q\ \mu^j_\infty\ f^j(z)\quad\quad\mbox{ in }C^l(D^2) \quad\forall l\in {\N}\quad.
\ee
Standard elliptic estimates applied to the system (\ref{IV.31}) imply that, modulo extraction of a subsequence
$\vec{F}_k$ converges weakly in $L^p$ for any $p<2$ \footnote{ Also weakly$^\ast$ in $L^{2,\infty}$.}to a map $\vec{F}_\infty$ which is \underbar{real} because
of (\ref{IV.33}). By Rellich Kondrachov compact embedding $\nabla\vec{\xi}_k$ strongly converges 
to $\nabla\vec{\xi}_\infty$ in $L^q$ for any $q<+\infty$. Hence using the Jacobian structures
 we have
 \be
 \label{IV.39}
 \begin{array}{l}
\ds \nabla^\perp\vec{F}^\Re_k\cdot\nabla\vec{\xi}_k=-div\lf[\vec{F}^\Re_k\cdot\nabla^\perp\vec{\xi}_k\rg]\\[5mm]
\ds\longrightarrow
-div\lf[\vec{F}_\infty\cdot\nabla^\perp\vec{\xi}_\infty\rg]= \nabla^\perp\vec{F}^\Re_\infty\cdot\nabla\vec{\xi}_\infty\\[5mm]
\ds \nabla^\perp\vec{F}^\Re_k\wedge\nabla\vec{\xi}_k=-div\lf[\vec{F}^\Re_k\wedge\nabla^\perp\vec{\xi}_k\rg]\\[5mm]
\ds\longrightarrow
-div\lf[\vec{F}_\infty\wedge\nabla^\perp\vec{\xi}_\infty\rg]= \nabla^\perp\vec{F}_\infty\wedge\nabla\vec{\xi}_\infty\\[5mm]
\ds\nabla^\perp\lf(\star(\vec{n}_k\res\vec{H}_k)\rg)\res\nabla\vec{\xi}_k=-div\lf[ \lf(\star(\vec{n}_k\res\vec{H}_k)\rg)\res\nabla^\perp\vec{\xi}_k \rg]\\[5mm]
\ds \longrightarrow -div\lf[ \lf(\star(\vec{n}_\infty\res\vec{H}_\infty)\rg)\res\nabla^\perp\vec{\xi}_\infty \rg]=
\nabla^\perp\lf(\star(\vec{n}_\infty\res\vec{H}_\infty)\rg)\res\nabla\vec{\xi}_\infty\ .          
\end{array}
\ee
Hence we have proved that $\vec{\xi}_\infty$ satisfies the following system : $\exists \vec{F}_\infty\in L^{2,\infty}(D^2,{\R}^m)$
such that
\[
\lf\{
\begin{array}{l}
\nabla^\perp\vec{F}_\infty\cdot\nabla\vec{\xi}_\infty=0\\[5mm]
\nabla^\perp\vec{F}_\infty\wedge\nabla\vec{\xi}_\infty=2\ (-1)^m\nabla^\perp\lf(\star(\vec{n}_\infty\res\vec{H}_\infty)\rg)\res\nabla\vec{\xi}_\infty\ .
\end{array}
\rg.
\]
This is equivalent to the fact that $\vec{\xi}_\infty$ satisfies the {\it Conformal Willmore } equation
and $\vec{\xi}_\infty$ is analytic.

\medskip

Assume now that
\[
\lim_{k\rightarrow +\infty}|\mu_k|=+\infty\quad.
\]
Then we consider $\vec{F}_k/|\mu_k|$ and dividing by $|\mu_k|$ equation (\ref{IV.31}) and passing to the limit in the Jacobian expressions as above we get the existence of a map $\vec{L}_\infty$ and a non zero holomorphic function $g_\infty(z)$ contained
in the span of $f^j$ such that
\be
\label{IV.40}
\p_z\vec{L}_{\infty}=e^{-2\la_\infty}\ g_\infty(z)\ \p_{\ov{z}}\vec{\xi}_\infty\quad.
\ee
This could have been done on all the balls $B_{r^i}(x_\infty^i)$ simultaneously \footnote{ Either $q_k:=\sum_{j=1}^N\mu_k^j f^j(z)\ dz\otimes dz$ is bounded in
the space $Q(J)$ of holomorphic quadratic forms of $(\Sigma,J)$ or goes to infinity in norm (for the Weil-Peterson hermitian product).} and hence, like in the proof of  lemma~\ref{lm-A.5},
since $\p_{\ov{z}}(e^{-2\la_\infty}\ \p_{\ov{z}}\vec{\xi}_\infty)=2^{-1}\vec{H}_{0,\infty}$,      (\ref{IV.39}) implies the existence
of a non  trivial holomorphic form $q$ of $Q(J)$ such that
\[
\Im(<q, h_{0,\infty}>_{WP})\equiv 0 \quad,
\]
where locally in holomophic coordinates $h_0:=\ov{\vec{H}_0}\ dz\otimes dz$. This is equivalent to the fact that $\vec{\xi}_\infty$ is isothermic. We have then
proved lemma~\ref{lm-IV.5} in the first case : when there exists a subsequence $\vec{\xi}_k$ which is not made of isothermic surfaces.

\medskip

\noindent{\bf Second case : all the $\vec{\xi}_k$ are isothermic conformal immersion of $(\Sigma,J)$.}

\medskip

This would mean that there exists a sequence of holomorphic quadratic differentials $q_k\ne 0$ such that
\be
\label{IV.41}
\Im(<q_k, h_{0,k}>_{WP})\equiv 0 \quad,
\ee
We can normalize $q_k$ in such a way that $<q_k,q_k>_{WP}=1$ and since $Q(J)$ is finite dimensional we can extract a subsequence such that $q_k$
converges strongly in any $C^l$ norm (for any $l\in {\N}$) to a non zero limiting holomorphic quadratic differential $q_\infty$. We have seen that
$d\vec{\xi}_k$ converges strongly in $L^p_{loc}({\Sigma}\setminus\{a_1\cdots a_N\})$ to a limiting $W^{2,2}$ immersion $\vec{\xi}_\infty$
of ${\Sigma}\setminus\{a_1\cdots a_N\}$ and since the second fundamental form of $\vec{\xi}_k$ is uniformly bounded in $L^2$, $h_{0,k}$ converges weakly in $L^2$
to the Weingarten Operator $h_{0,\infty}$ of $\vec{\xi}_\infty$. We can then pass in the limit in the identity (\ref{IV.41}). This implies that
\be
\label{IV.42}
\Im(<q_\infty, h_{0,\infty}>_{WP})\equiv 0 \quad,
\ee
from which we deduce that $\vec{\xi}_\infty$ is an isothermic immersion of $\Sigma\setminus\{a_1,\cdots a_N\}$ into ${\R}^m$. This concludes the proof of lemma~\ref{lm-IV.5} in all cases. \hfill $\Box$

\medskip

The proof theorem~\ref{th-IV.2} can be finished exactly like in the proof of theorem~\ref{th-III.1} in order to exploit the assumption $\inf_{\vec{\Phi}\in{\mathcal E}_\Sigma^c}W(\vec{\Phi})\le 8\pi$ and ''remove'' the singularity points $a_i$.\hfill $\Box$

\section{The conformal constraint and isothermic immersions.}

We first describe the immersions which are the singular points for the map which assigns to an immersion it's conformal class as we will
prove in lemma~\ref{lm-A.5} : the isothermic
immersions.

\medskip

\noindent{\bf Proof of proposition~\ref{pr-0.6}.}

Isothermic immersions have been defined in definition~\ref{df-0.5}. First let consider an isothermic immersion. There are locally, away from umbilic points, conformal coordinates in which
the second fundamental form is diagonal. Hence this means that the Weingarten map is real in such charts. Take two such complex charts $z=x_1+ix_2$ and $\xi=\xi_1+i\xi_2$ overlaping
on some open set. The Weingarten operator is independent of the complex chart and we have
\[
h_0=  \p_z\vec{n}_{\vec{\Phi}}\res\p_{z}\vec{\Phi}\ dz\otimes dz=\p_{\xi}\vec{n}_{\vec{\Phi}}\res \p_{\xi}\vec{\Phi}\ d\xi\otimes d\xi\quad.
\]
Our assumption reads $\Im(\p_z\vec{n}_{\vec{\Phi}}\res\p_{z}\vec{\Phi})=\Im(\p_{\xi}\vec{n}_{\vec{\Phi}}\res \p_{\xi}\vec{\Phi})=0$, moreover, since we are away from umbilic points $h_0\ne 0$
that is then $\Re(\p_z\vec{n}_{\vec{\Phi}}\res\p_{z}\vec{\Phi})\ne 0$. Thus we have that
\[
\frac{\Re(\p_{\xi}\vec{n}_{\vec{\Phi}}\res \p_{\xi}\vec{\Phi})}{\Re(\p_z\vec{n}_{\vec{\Phi}}\res\p_{z}\vec{\Phi})}= (z'(\xi))^2
\]
This implies that the imaginary part of the holomorphic function $(z'(\xi))^2$ is zero which implies that $z'(\xi)$ is constant and $(z'(\xi))^2$ is a real constant. Thus
$dz\otimes dz$ and $d\xi\otimes d\xi$ are proportional to eachother by a real non zero constant and this implies that the form $dz\otimes dz$ extends to an holomorpfhic quadratic differential $q$ of the riemann surface $\ti{\Sigma}$ obtained by withdrawing to $\Sigma$ the umbilic points of the immersion $\vec{\Phi}$ and we have by construction $<q,h_0>_{WP}=0$. 

We are now proving the reciproque. Let $q$ be an holomorphic quadratic differential of $\ti{\Sigma}$. Away from the isolated zeros of $q$ we can choose complex coordinate $z$ 
such that $q(z)=dz\otimes dz$ (indeed in arbitrary complex coordinates $q(\xi)=f(\xi)\ d\xi\otimes d\xi$ where $f$ is holomorphic and just choose $z(\xi)=\sqrt{f(\xi)}$.
In these coordinates  the condition (\ref{aza}) implies $\Im(\vec{H}_0)=0$ which means that the second fundamental form is diagonal in these complex coordinates
and hence $\vec{\Phi}$ is local isothermic.
\hfill $\Box$

\medskip

\noindent{\bf Proof of proposition~\ref{pr-0.7}.}

Let $\vec{\Phi}$ be an immersion. From computations in section III.2.2 of \cite{BR} we have in complex coordinales
\be
\label{V.1}
\p_{\ov{z}}\lf(e^{-2\la}\ \p_{\ov{z}}\vec{\Phi}\rg)=2^{-1}\ \vec{H}_0\quad.
\ee
Assume $\vec{\Phi}$ is local isothermic, because of the previous proposition there exists complex coordinates in which $\Im[\vec{H}_0]=0$. Hence
in these coordinates
\be
\label{V.2}
\Im\lf[\p_{\ov{z}}\lf(e^{-2\la}\ \p_{\ov{z}}\vec{\Phi}\rg)\rg]=0\quad\quad\mbox{ in }D^2\ .
\ee
Let $\vec{L}=\vec{L}_\Re+i\vec{L}_\Im\in{\R}^m\otimes{\C}$ be the unique solution to the following elliptic system
\[
\lf\{
\begin{array}{l}
\ds\p_{z}{\vec{L}}=e^{-2\la}\ \p_{\ov{z}}\vec{\Phi}\quad\quad\mbox{ in }D^2\\[5mm]
\ds\vec{L}_{\Im}=0\quad\quad\quad\mbox{ on }\p D^2
\end{array}
\rg.
\]
Then, because of (\ref{V.2}), $\vec{L}_\Im$ solves
\[
\lf\{
\begin{array}{l}
\ds\Delta\vec{L}_\Im=0\quad\quad\mbox{ in }D^2\\[5mm]
\ds\vec{L}_\Im=0\quad\quad\quad\mbox{ on }\p D^2\quad.
\end{array}
\rg.
\]
This implies that $\vec{L}=\vec{L}_\Re\in{\R}^m$. Hence we have proved (\ref{bzb}). Assuming now that (\ref{bzb}) holds, we obtain the existence of complex coordinates
such that  (\ref{V.2}) is satisfied which implies from (\ref{V.1}) that $\Im(\vec{H}_0)=0$ and from which we deduce that $\vec{\Phi}$ is isothermic. This finishes the proof of 
proposition~\ref{pr-0.7}.\hfill $\Box$

\medskip

\noindent{\bf Proof of proposition~\ref{pr-0.7b}.}

An elementary computation gives for any pair $\vec{\Phi}$ and $\vec{L}$ from $D^2$ into ${\R}^m$
\[
4\, \Im\lf(\p_{\ov{z}}\vec{L}\cdot\p_z\vec{\Phi}\rg)=\nabla \vec{L}\cdot\nabla^{\perp}\vec{\Phi}\quad\quad\mbox{ and }\quad\quad
4\, \Im\lf(\p_{\ov{z}}\vec{L}\wedge\p_z\vec{\Phi}\rg)=\nabla \vec{L}\wedge\nabla^{\perp}\vec{\Phi}
\]
Hence (\ref{bzb}) clearly implies (\ref{0.013}).

Assuming now that (\ref{0.013}) holds, we have then the existence of $\vec{L}\in{\R}^m$ such that
\[
\lf\{
\begin{array}{l}
\vec{e}_1\wedge\p_y\vec{L}=\vec{e}_2\wedge\p_x\vec{L}\\[5mm]
\vec{e}_1\cdot\p_y\vec{L}=\vec{e}_2\cdot\p_x\vec{L}\quad.
\end{array}
\rg.
\]
A short computation shows that this implies the existence of  $(a,b)\in{\R}^2$ such that
\[
\nabla^\perp\vec{L}=\lf(
\begin{array}{cc}
b&a\\[3mm]
a&-b
\end{array}
\rg)\ 
\lf(
\begin{array}{c}
\vec{e}_1\\[3mm]
\vec{e}_2
\end{array}
\rg)
\]
or in other words, introducing $f:=e^{\la} (a+ib)$, one has
\be
\label{V.3}
\p_z\vec{L}=f\ e^{-2\la}\p_{\ov{z}}\vec{\Phi}\quad.
\ee
Since the components of $\vec{L}$ are real, and by consequence the components of $\Delta\vec{L}$ are real as well, we have that
\[
\Im\lf[\p_{\ov{z}}\lf(f\ e^{-2\la}\p_{\ov{z}}\vec{\Phi}\rg)\rg]=0\quad.
\] 
Using (\ref{V.1}) this gives
\[
\p_{\ov{z}}f\ \vec{e}_{\ov{z}}-\p_z\ov{f}\ \vec{e}_{z}=-\frac{e^{\la}}{2}\ \lf[f\ \vec{H}_0-\ov{f}\ \ov{\vec{H}_0}\rg]\quad.
\]
Since $\vec{H}_0$ is orthogonal to the tangent plane of the immersion and since $\vec{e}_z$ and $\vec{e}_{\ov{z}}$ are in the complexified space
to the tangent space and are linearly independent we deduce
\[
\p_{\ov{z}}f=0\quad.
\]
Take now $w=\sqrt{f}$ equation (\ref{V.3}) becomes
\[
\p_w\vec{L}=e^{-2\la}\ |f|\ \p_{\ov{w}}\vec{\Phi}
\]
and one observes that $e^\la\ |f|^{-1/2}$ is the new conformal factor of $\vec{\Phi}$ in the coordinate $w$, which means that $\vec{\Phi}$ satisfies (\ref{bzb}) in these coordinates
and hence, from the previous proposition, $\vec{\Phi}$ is an isothermic immersion.\hfill $\Box$

\medskip

Finally we prove that the {\it global isothermic immersions} are the degenerate points for the conformal class mapping. Precisely we prove the following
result.

\medskip

\begin{Lma}
\label{lm-A.5}
Let $\vec{\Phi}$ be a conformal $W^{2,2}\cap W^{1,\infty}$ immersion of a closed Riemann surface $(\Sigma,J)$ of genus larger or equal to one. Consider in a neighborhood
of $0$ the map 
\[
{\mathcal C} \ :\ \vec{w} \in  W^{2,2}\cap W^{1,\infty}\ \longrightarrow {\mathcal C}(\vec{w})\in {\mathcal T}_\Sigma\quad,
\]
where ${\mathcal T}_\Sigma$ is the Teichm\"uller Space associated to the surface $\Sigma$ and ${\mathcal C}(\vec{w})$ is the
Teichm\"uller class issued from the immersion $\vec{\Phi}+\vec{w}$ with fixed generators of the $\pi_1$ on $\Sigma$.
The map ${\mathcal  C}$ is $C^1$ in a neighborhood of $0$. Identifying ${\mathcal T}_\Sigma$ with the space $Q(J)$
of holomorphic quadratic differentials\footnote{See for instance theorem 4.2.2 in \cite{Jo}.} on $(\Sigma,J)$, its differential at $0$ is given by
\be
\label{A.42b}
d{\mathcal C}(0)\cdot\vec{\nu}=8\,\sum_{j=1}^Qq_j\ \lf<q_j, \p_z\vec{\nu}\cdot\p_z\vec{\Phi}\ dz\otimes dz\rg>_{WP}\quad,
\ee
where $<\cdot,\cdot>_{WP}$ is the Weil Peterson Hermitian product and $(q_j)_{j=1\cdots Q}$ is an orthonormal basis of $Q(J)$ for this
product. Moreover, if $\vec{\Phi}$ is not an isothermic surface, $d{\mathcal C}(0)$ is a submersion onto the space of holomorphic quadratic differentials of $(\Sigma,J)$. \hfill $\Box$
\end{Lma}
{\bf Proof of lemma~\ref{lm-A.5}.}
Let $\vec{\Phi}$ be a conformal $W^{2,2}\cap W^{1,\infty}$ immersion of $(\Sigma,J)$ and $\vec{w}$ be a map in $W^{2,2}\cap W^{1,\infty}(\Sigma,{\R}^m)$,
small enough in this space, in such a way that $\vec{\Phi}+\vec{w}$ still defines an immersion. Denote by $J^{\vec{w}}$ the $W^{1,2}$ complex structure
defined by $(\vec{\Phi}+\vec{w})^\ast g_{{\mathbb{R}}^m}$. Using lemma~\ref{lm-A.3} there exists a covering of $\Sigma$
by disks $(U_i)_{i\in I}$ and $W^{1,\infty}\cap W^{2,2}$ diffeomorphisms $\psi_i^{\vec{w}}$ from $D^2$ into $U_i$ such that $(\vec{\Phi}+\vec{w})\circ\psi_i$
is conformal. Considering now $\Sigma$ together with the covering $U_i$ and the holomorphic transition maps 
\[
h_{ij}^{\vec{w}}(z):=(\psi_j^{\vec{w}})^{-1}\circ\psi_i^{\vec{w}}\quad,
\] 
which satisfy of course the cocycle condition $h_{ij}^{\vec{w}}\circ h_{jk}^{\vec{w}}\circ h_{ki}^{\vec{w}}(z)=z$, we have defined a new smooth complex
structure on $\Sigma$, $\ti{J}^{\vec{w}}$ which is equivalent to $(\Sigma, J^{\vec{w}})$ : there exist smooth conformal diffeomorphisms, $\varphi^{\vec{w}}$ from $D^2$
into $(U_i,\ti{J}^{\vec{w}})$ and an homeomorphism\footnote{$\Sigma$ together with the charts $(U_i,\psi_i)$ defines a smooth complex manifold since $\psi_j^{-1}\circ\psi_i$, are holomorphic,  the smooth complex structure being given by the multiplication by $i$ in the charts. It admits then a constant scalar curvature
metric  $h^{\vec{w}}$ and $\Psi^{\vec{w}}$ is the harmonic diffeomorphism from $(\Sigma, g^{\vec{w}})$ into $(\Sigma, h^{\vec{w}})$ isotopic to the identity, see \cite{Jo}.}
$\Psi^{\vec{w}}$ of $\Sigma$, isotopic to the identity which is conformal from $(\Sigma, J^{\vec{w}})$ into $(\Sigma, \ti{J}^{\vec{w}})$ and such that $\Psi^0=id_{\Sigma}$ and 
$\ti{J}^0=J$. Hence $\Psi^{\vec{w}}$ is bilipschitz and 
both $\Psi^{\vec{w}}$ and $(\Psi^{\vec{w}})^{-1}$ are $W^{2,2}$. By replacing now $\vec{\Phi}$ by $(\vec{\Phi}+\vec{w})\circ (\Psi^{\vec{w}})^{-1}$, if one shows that
that ${\mathcal C}$ is $C^1$ at $0$ one has shown that ${\mathcal C}$ is $C^1$ in a neighborhood of the origin.

In order to show that ${\mathcal C}$ is $C^1$ at $0$ it suffices to show that the mappings wich to $\vec{w}$ assigns the family of holomorphic transition
functions $h_{ij}^{\vec{w}}$ is $C^1$ from $W^{1,\infty}\cap W^{2,2}$ into $C^0$ (which implies that it is $C^1$ from $W^{1,\infty}\cap W^{2,2}$
into $C^l$ for an arbitrary $l$ on a slightly small covering). In order to show that it suffices to show that the mappings which to $\vec{w}$ assigns
$\psi_i^{\vec{w}}$ and $(\psi_i^{\vec{w}})^{-1}$ are $C^1$ from $W^{1,\infty}\cap W^{2,2}$ into $W^{2,p}$ for some $p>1$. This can be done
following carefully the construction of conformal coordinates in lemma~\ref{lm-A.3}. We leave the details to the reader.

We compute now the differential of ${\mathcal C}$ at the origin. As above $h^{\vec{w}}$ denotes the metric of constant scalar curvature
compatible with $(\Sigma,\ti{J}^{\vec{w}})$ and we denote simply by $h$ the constant scalar curvature compatible with $(\Sigma,J)$.
Let $u^{\vec{w}}$ be the harmonic map from $(\Sigma,J)$ into $(\Sigma,\ti{J}^{\vec{w}})$ isotopic to the identity given by corollary 3.10.1 in \cite{Jo}.
The map $C(\vec{w})$ is given explicitly by
\[
C(\vec{w})=\sum_{j=1}^QA_j^{\vec{w}}\ q_j\quad,
\]
where
\[ 
A_j^{\vec{w}}:=\lf<q_j, \lf[\lf<\p_xu^{\vec{w}},\p_xu^{\vec{w}}\rg>_{h^{\vec{w}}}-\lf<\p_yu^{\vec{w}},\p_yu^{\vec{w}}\rg>_{h^{\vec{w}}}-2i\ \lf<\p_xu^{\vec{w}},\p_yu^{\vec{w}}\rg>_{h^{\vec{w}}}\rg]\ (dz)^2\rg>_{WP}
\]
Denote $v^{\vec{w}}:=(\Psi^{\vec{w}})^{-1}\circ u^{\vec{w}}$ and $e^{\mu^{\vec{w}}} g^{\vec{w}}=(\Psi^{\vec{w}})^\ast h^{\vec{w}}$. Hence we have in particular
\[ 
A_j^{\vec{w}}:=\lf<q_j, e^{2\mu^{\vec{w}}}\lf[\lf<\p_xv^{\vec{w}},\p_xv^{\vec{w}}\rg>_{g^{\vec{w}}}-\lf<\p_yv^{\vec{w}},\p_yv^{\vec{w}}\rg>_{g^{\vec{w}}}-2i\,\lf<\p_xv^{\vec{w}},\p_yv^{\vec{w}}\rg>_{g^{\vec{w}}}\rg]\,(dz)^2\rg>_{WP}
\]
Let $X^{\vec{\nu}}:=dv^{\vec{w}}(0)\cdot\nu$. Since $v^{0}=id_\Sigma$, we have, writing locally $q_j=q^z_j(z)\ dz\otimes dz$ in complex coordinates \footnote{As usual
the Weil-Peterson metric is expressed using local complex coordinates bearing in mind that the expression of the integrand is independent of this local choice} satisfying in particular
$|\p_x\vec{\Phi}|=|\p_y\vec{\Phi}|=e^\la$ and then in which $h=e^{2\mu(0)}\ g=e^{2\mu(0)+2\la}\ [dx^2+dy^2]$
\[
\begin{array}{l}
\ds dA_j^{\vec{w}}(0)\cdot\nu=\int_{\Sigma}e^{-2\la-2\mu(0)}\ q^z_j(z)\, e^{2\mu(0)}\,\lf[2\lf<\p_xX^{\vec{\nu}},\p_xid_{\Sigma}\rg>_g-2\lf<\p_yX^{\vec{\nu}},\p_yid_{\Sigma}\rg>_g\rg.
\\[5mm]
\ds\quad\quad\quad+2i\,\lf<\p_xX^{\vec{\nu}},\p_yid_{\Sigma}\rg>_g+2i\,\lf<\p_yX^{\vec{\nu}},\p_xid_{\Sigma}\rg>_g\\[5mm]
\ds\quad\quad\quad\lf.+dg_{11}^{\vec{w}}(0)\cdot\nu-dg_{22}^{\vec{w}}(0)\cdot\nu+2i\,dg_{12}^{\vec{w}}(0)\cdot\nu\rg]\ \frac{i}{2}\, dz\wedge d\ov{z}
\end{array}
\]
Decomposing the vector-field $X^{\vec{\nu}}$ as follows : $X^{\vec{\nu}}=X^{\vec{\nu}}_z\ \p_z+X^{\vec{\nu}}_{\ov{z}}\ \p_{\ov{z}}$ (where $\p_z:=2^{-1}(\p_x-i\p_y)$ ) and observing that $dg^{\vec{w}}(0)\cdot\nu=(\p_{x_j}\vec{\Phi}\cdot\p_{x_i}\vec{\nu}+\p_{x_i}\vec{\Phi}\cdot\p_{x_j}\vec{\nu})_{ij}$ gives
\be
\label{A.43}
\begin{array}{l}
\ds dA_j^{\vec{w}}(0)\cdot\nu=2i\,\int_{\Sigma}q^z_j(z)\ \p_{\ov{z}}X_z^{\vec{\nu}}\  dz\wedge d\ov{z}\\[5mm]
\ds\quad\quad\quad+4i\ \int_{\Sigma} e^{-2\la}\ q^z_j(z)\ \p_{\ov{z}}\vec{\Phi}\cdot\p_{\ov{z}}\vec{\nu}\ dz\wedge d\ov{z}
\end{array}
\ee
Observe that
\[
q^z_j(z)\ \p_{\ov{z}}X_z^{\vec{\nu}}\  dz\wedge d\ov{z}=q^z_j(z)\ dz\wedge d(X_z^{\vec{\nu}})=d\lf[X_z^{\vec{\nu}}\ q^z_j(z)\ dz\rg]\quad.
\]
Let $\xi$ be another complex coordinates. We have that $q_j=q_j^z(z)\ dz\otimes dz=(\xi')^{-2}\ q^z_j(z)\ d\xi\otimes d\xi$. Then
$q_j^\xi(\xi)=(\xi')^{-2}\ q^z_j(z)$. We also have
$X=X^{\vec{\nu}}_z\ \p_z+X^{\vec{\nu}}_{\ov{z}}\ \p_{\ov{z}}=\xi'\,X^{\vec{\nu}}_z\ \p_\xi+\ov{\xi}'\, X^{\vec{\nu}}_{\ov{z}}\ \p_{\ov{\xi}}$. Hence $\xi'\,X^{\vec{\nu}}_z=X^{\vec{\nu}}_\xi$
and
\[
\al=X_z^{\vec{\nu}}\ q^z_j(z)\ dz=X^{\vec{\nu}}_\xi (\xi')^{-1}\ q_j^\xi(\xi)\ (\xi')^{2}\ \xi'\ d\xi=X^{\vec{\nu}}_\xi\ q_j^\xi(\xi)\ d\xi
\]
is an intrinsic 1-form globally defined on $\Sigma$. Thus
\[
\int_{\Sigma}q^z_j(z)\ \p_{\ov{z}}X_z^{\vec{\nu}}\  dz\wedge d\ov{z}=\int_{\Sigma}d\al=0
\]
and (\ref{A.43}) implies (\ref{A.42b}). It remains to prove that $d{\mathcal C}(0)$ is a submersion onto the space of holomorphic quadratic differentials of $(\Sigma,J)$.

In local conformal coordinates for $\vec{\Phi}$ we denote $\vec{e}_i=e^{-\la}\p_{x_i}\vec{\Phi}$ where $e^\la=|\p_{x_i}\vec{\Phi}|$. Let $(\vec{n}_\al)_{\al=1\cdots m-2}$ be a local orthonormal frame of the normal bundle to $\vec{\Phi}(\Sigma)$. We denote $h^\al_{ij}:=-e^{-\la}\ \vec{e}_i\cdot\p_{x_j}\vec{n}_\al$. The {\it Weingarten map}  is given by
\[
\vec{H}_0=\sum_{\al=1}^{m-2} H^\al_0\ \vec{n}_\al:=\frac{1}{2}\sum_{\al=1}^{m-2}(h^\al_{11}-h^\al_{22}+2\, i\ h^\al_{12})\ \vec{n}_\al\quad.
\]
Denote $\vec{e}_z:=e^{-\la}\p_z\vec{\Phi}=2^{-1}(\vec{e}_1-i\,\vec{e}_2)$ and $\vec{e}_{\ov{z}}:=e^{-\la}\p_{\ov{z}}\vec{\Phi}=2^{-1}(\vec{e}_1+i\,\vec{e}_2)$. Some elementary
computations give (see for instance \cite{BR} section III.2.2)
\be
\label{A.44}
\p_{\ov{z}}(e^{-\la}\vec{e}_{\ov{z}})=2^{-1}\ \vec{H}_0\quad.
\ee
Let $\vec{\nu}$ be a map supported in the domain of definition for the local conformal charts that we identify with $D^2$. We assume $\vec{\nu}$ to be in $W^{2,2}_0\cap W^{1,\infty}(D^2)$.
Denote by $f_j(z)\ dz\otimes dz$ the expression of the basis $q_j$ in this conformal charts in such a way that $f_j(z)$ are holomorphic functions on $D^2$.
The expression (\ref{A.42b}) of $d{\mathcal C}(0)\cdot\nu$ gives
\be
\label{A.45}
\begin{array}{rl}
\ds d{\mathcal C}(0)\cdot\vec{\nu}&\ds=8\,\sum_{j=1}^Qq_j\ \int_{D^2} f_j(z)\ e^{-\la}\vec{e}_{\ov{z}}\cdot\p_{\ov{z}}\vec{\nu}\ \frac{i}{2}dz\wedge d\ov{z}\\[5mm]
\ds &\ds=-8\,\sum_{j=1}^Qq_j\ \int_{D^2} f_j(z)\ \p_{\ov{z}}(e^{-\la}\vec{e}_{\ov{z}})\cdot\vec{\nu}\ \frac{i}{2}dz\wedge d\ov{z}\\[5mm]
 &\ds=-8\,\sum_{j=1}^Qq_j\ \int_{D^2} f_j(z)\ \vec{{H}}_0\cdot\vec{\nu}\ \frac{i}{2}dz\wedge d\ov{z}
\end{array}
\ee
%A classical computation gives
%\[
%\begin{array}{rl}
%\ds W(\vec{\Phi})&\ds=\int_{\Sigma}|\vec{H}|^2\ dvol_{g}= 4^{-1}\int_{\Sigma}|d\vec{n}_{\vec{\Phi}}|^2_g\ dvol_g+\pi\ \chi(\Sigma)\\[5mm]
 % &\ds=2^{-1}\int_{\Sigma}|\vec{H}|^2\ dvol_g+ 2^{-1}\int_{\Sigma}|\vec{H}_0|^2\ dvol_g+\pi\ \chi(\Sigma)
  %\end{array}
%\]
%Since we are considering the case $g(\Sigma)\ge 1$, $\chi(\Sigma)=2-2\,g(\Sigma)\le 0$ hence we have
%\[
%\int_{\Sigma}|\vec{H}_0|^2\ge\int_{\Sigma}|\vec{H}|^2> 4\pi>0\quad.
%\]
%the last inequality being a consequence of the fact that $\Sigma$ is  a closed surface and hence the universal
%Willmore energy lower bound \footnote{ See for instance \cite{Wi}.}  by $4\pi$ applies. Thus $\vec{H}_0$ cannot be
%uniformly zero and we can find a disk on which there exist conformal coordinates and in which $\int_{D^2} |\vec{H}_0|^2\ne 0$.
If 
\[
\nu\longrightarrow \lf(\int_{D^2} f_j(z)\ \vec{{H}}_0\cdot\vec{\nu}\ \frac{i}{2}dz\wedge d\ov{z}\rg)_{j=1\cdots Q}
\] 
does not have a complex $N$ dimensional Range then it would mean that the real $2N$ linear forms on $W^{2,2}\cap W^{1,\infty}$ given by
\[
\nu\longrightarrow \lf(\Re\lf[\int_{D^2} f_j(z)\ \vec{{H}}_0\cdot\vec{\nu}\ \frac{i}{2}dz\wedge d\ov{z}\rg],\Im\lf[\int_{D^2} f_j(z)\ \vec{{H}}_0\cdot\vec{\nu}\ \frac{i}{2}dz\wedge d\ov{z}\rg]\rg)_{j=1\cdots Q}
\] 
are linearly dependent. This is equivalent to the existence of a non trivial family of real numbers $(\mu_j,\delta_j)_{j=1\cdots Q}\in {\R}^{2Q}$ such that

$\forall \vec{\nu}\in W^{2,2}\cap W^{1,\infty}$
\[
\sum_{j=1}^Q\mu_j\ \Re\lf[\int_{D^2} f_j(z)\ \vec{{H}}_0\cdot\vec{\nu}\ \frac{i}{2}dz\wedge d\ov{z}\rg]+\delta_j\ \Im\lf[\int_{D^2} f_j(z)\ \vec{{H}}_0\cdot\vec{\nu}\ \frac{i}{2}dz\wedge d\ov{z}\rg]=0
\]
or in other words there would exist $(\delta_j+i\mu_j)\in {\C}^Q$, being not all equal to zero, such that
\[
\Im\lf[\sum_{j=1}^N(\delta_j+i\mu_j)\ f_j(z)\ \vec{H}_0\rg]\equiv 0
\]
In other words again, this would mean that there exist a non zero holomorphic quadratic form $q$ in $Q_J$ such that
\[
\Im(<q, h_0>_{WP})\equiv 0 \quad,
\]
where locally  $h_0:=e^{2\la}\ \ov{\vec{H}_0}\ dz\otimes dz$. This is equivalent to the fact that $\vec{\Phi}$ is isothermic. Hence if we make the assumption that $\vec{\Phi}$ is not
isothermic the dimension of the range of $d{\mathcal C}(0)$ is $Q=dim(Q(J))$, 
 which concludes the proof of lemma~\ref{lm-A.5}.\hfill$\Box$

\appendix
\section{Appendix}
\reset
\begin{Lma}
\label{lm-A.00}
Let  $g$ and $h$ be two metrics at a point $p\in D^2$ such that 
\[
|Dis(g)|<1-2^{-k}\quad.
\]
where $k\in {\N}$. Then the following inequality holds
\be
\label{+A.1}
\frac{1}{2}\ \inf_{X\in {\R}^2\setminus\{0\}} \frac{|X|^2_g}{|X|^2_h}\le\frac{tr(g)}{tr(h)}\le 2^k\ \sup_{X\in {\R}^2\setminus\{0\}} \frac{|X|^2_g}{|X|^2_h}\quad.
\ee
\hfill $\Box$
\end{Lma}
{\bf Proof of Lemma~\ref{lm-A.00}.}
Observe that 
\[
\begin{array}{rl}
\ds |X|^2_g&\ds=\lf[\frac{tr(g)+\Re[H(g)]}{2}\rg]\ X_1^2+\lf[\frac{tr(g)-\Re[H(g)]}{2}\rg]\ X_2^2\\[5mm]
 &\ds-\Im[H(g)]\ X_1\ X_2\quad.
 \end{array}
\]
Hence, denoting also $X=X_1+iX_2$ 
\[
|X|_g^2=\Re\lf[\lf(\frac{H(g)}{2}\ X+\frac{tr(g)}{2}\ \ov{X}\rg)\ X\rg]\quad.
\] 
We deduce that
\be
\label{A.02}
\frac{|X|^2_g}{|X|^2_h}=\frac{tr(g)}{tr(h)}\ \frac{1+\Re\lf[Dis(g)\ e^{2i\,\theta}\rg]}{1+\Re\lf[Dis(h)\ e^{2i\,\theta}\rg]}\quad,
\ee
where $X=|X|\, e^{2i\,\theta}$. 

Since $|Dis(g)|<1-2^{-k}$, using (\ref{A.02}) we have that 
\[
\frac{tr(g)}{tr(h)}\frac{2^{-k}}{\inf_\theta\lf| 1+\Re\lf[Dis(h)\ e^{2i\,\theta}\rg]\rg|}\le \sup_{X\in {\R}^2\setminus\{0\}} \frac{|X|^2_g}{|X|^2_h}\quad.
\]
Hence in particular $\inf_\theta\lf| 1+\Re\lf[Dis(h)\ e^{2i\,\theta}\rg]\rg|>0$ and by taking $e^{2i\theta}:={Dis(h)}/|{Dis(h)}|$ (in the case when $Dis(h)\ne 0$) we see that
there is $\theta$ such that $ 1+\Re\lf[Dis(h)\ e^{2i\,\theta}\rg]>0$, by continuity this implies that
\[
\forall \theta\in {\R}\quad\quad1+\Re\lf[Dis(h)\ e^{2i\,\theta}\rg]>0
\]
This implies that $\inf_\theta\lf| 1+\Re\lf[Dis(h)\ e^{2i\,\theta}\rg]\rg|=1-|Dis(h)|<1$ from which we deduce the upper bound in (\ref{+A.1}).
Take now again $X=|X|\, e^{i\theta}$ in such a way that $e^{2i\theta}:=Dis(h)/|Dis(h)|$ (still in the cas when $Dis(h)\ne 0$). For this $X$ we have
\[
\frac{tr(h)}{tr(g)}= \frac{|X|^2_h}{|X|^2_g}\ \frac{1+\Re\lf[Dis(g)\ e^{2i\,\theta}\rg]}{1+|Dis(h)|}\le 2\ \sup_{X\ne 0}\frac{|X|^2_h}{|X|^2_g}
\]
which gives the lower bound in (\ref{+A.1}). \hfill $\Box$

\begin{Lma}
\label{lm-A.2}
There exists $\ep_0>0$ such that for any $\vec{\xi}\in W^{2,2}\cap W^{1,\infty}(D^2,{\R}^m)$ satisfying
\be
\label{A.5a}
|Dis(\vec{\xi})|<\ep_0
\ee
and
\[
\int_{D^2}|\nabla\vec{n}_{\vec{\xi}}|^2_g\ dvol_g<\frac{4\pi}{3}\quad,
\]
where $g:=\vec{\xi}^\ast g_{{\R}^m}$, there exists $\vec{e}_1$ and $\vec{e}_2$ in $W^{1,2}(D^2,S^{m-1})$ such that 
\be
\label{A.6a}
\vec{e}_1\cdot \vec{e}_2=0\quad\quad,\quad\quad n_{\vec{\xi}}=\vec{e}_1\wedge \vec{e}_2 \quad,
\ee
\be
\label{A.7}
\int_{D^2}\lf[|\nabla \vec{e}_1|_g^2+|\nabla \vec{e}_2|_g^2\rg]\ dvol_g\le 2\int_{D^2}|\nabla\vec{n}_{\vec{\xi}}|^2_g\ dvol_g
\ee
and
\be
\label{A.8}
\lf\{
\begin{array}{l}
\ds d(\ast_g(\vec{e}_1,d \vec{e}_2))=0\quad,\\[5mm]
\ds\iota_{\p D^2}^\ast\ast_g(\vec{e}_1,d \vec{e}_2)=0\quad .
\end{array}
\rg.
\ee
where $\iota_{\p D^2}$ is the canonical inclusion of $\p D^2$ in ${\R}^2$. 
\hfill$\Box$
\end{Lma}
{\bf Proof of lemma~\ref{lm-A.2}.}
This lemma is proved in \cite{Hel} lemma 5.1.4 for $\vec{\xi}$ conformal - which implies that $g=\vec{\xi}^\ast g_{{\R}^m}=e^\la\ dx^2+dy^2$
and hence, in that case, for any function $f$ 
$$
\int_{D^2}|df|^2_g\ dvol_g=\int_{D^2}|\nabla f|^2\ dx\,dy\quad .
$$
We now explain how the strategy in \cite{Hel} adapts to the case when $Dis(\vec{\xi})=0$ is replaced by (\ref{A.5a}).
The assumption (\ref{A.5a}) implies that
\be
\label{A.8a}
\begin{array}{l}
\forall (x,y)\in D^2\quad\forall X\in T_{(x,y)}D^2\quad\quad \\[5mm]
\quad\quad (1-\ep_0^2)\ |X|^2_{g_0}\le {(det g)^{-1}}\ |X|^2_g\le\ (1+\ep_0^2)\ |X|^2_{g_0}\quad,
\end{array}
\ee
where $g_0$ is the flat metric $dx^2+dy^2$. We can first assume that $n_{\vec{\xi}}$ is a smooth map from $D^2$ into the Grassman Space of oriented 2-planes
 in $Gr_2({\R}^m)$ - which are dense in $W^{1,2}(D^2,Gr_2({\R}^m))$ see step 6 of the proof of lemma 5.1.4 of \cite{Hel}.
Let $\ti{e}:=(\ti{e}_1,\ti{e}_2)$ be a smooth orthonormal 2-frame\footnote{this trivialization  exists since we are now working
with a smooth $\vec{n}_{\xi}$ and the pull-back over $D^2$ by   $\vec{n}_{\xi}$ of the tautological bundle $SO(m)/SO(m-2)$ over $Gr_2({\R}^m)$ 
is trivial since $D^2$ is contractible} in ${\R}^m$ realizing (\ref{A.6a}).
For each $r\in (0,1]$ we minimize
\[
F_r(\theta)=\int_{D_r^2} |(e_1,d e_2)|^2_g\ dvol_g\quad,
\]
among $\theta\in W^{1,2}(D_r^2,{\R})$ and $e_1+ie_2=e^{i\theta}(\ti{e}_1+i\ti{e}_2)$. Since $(e^r_1,d e^r_2)= d\theta+(\ti{e}_1,d\ti{e}_2)$ $F_r$ is convex
and the minimum is achieved by a unique $e^r$ satisfying 
\be
\label{A.9a}
\lf\{
\begin{array}{l}
d(\ast_g(e^r_1,d e^r_2))=0\quad\quad\mbox{ in }D^2_r\\[5mm]
\iota_{\partial D^2_r}^\ast\ast_g(e^r_1,de^r_2)=0
\end{array}
\rg.
\ee
where $\iota_{\p D^2_r}$ is the canonical embedding of $\p D^2_r$ in ${\R}^2$. Hence there exists a unique function $f^r\in W^{1,2}_0(D^2_r,{\R})$ such that $\ast_g(e^r_1,de^r_2)=df^r$ and $f^r$ satisfies
\be
\label{A.9b}
\lf\{
\begin{array}{l}
\ds\Delta_{g_0} f^r=\p_xe^r_1\cdot\p_ye^r_2-\p_xe^r_2\cdot\p_y e^r_1+\ast_{g_0}d( (\ast_{g_0}-\ast_g)\ df^r)\quad\quad\mbox{ in }D^2_r\\[5mm]
\ds f^r=0\quad\quad\quad\mbox{ on }\quad\p D^2_r
\end{array}
\rg.
\ee
where $\Delta_{g_0}$ is the Laplace operator for the flat metric $ \Delta_{g_0}=-[\p^2_{x}+\p^2_{y}]$. 
Because of (\ref{A.8a}) we have that $\forall \al$ one-form on $D^2$ $|(\ast_g-\ast_{g_0})\,\al|\le\ep_0\ |\al|$. Hence
Wente estimates together
with more standard elliptic estimates gives the bound
\be
\label{A.9c}
\begin{array}{rl}
\ds\int_{D^2_r}|\nabla f^r|^2\ dx\,dy&\ds\le \frac{3}{16\pi}\int_{D^2_r}|\nabla e_1^r|^2\ dx\,dy\ \int_{D^2_r}|\nabla e_2^r|^2\ dx\,dy\\[5mm]
 &\ds\quad+C\ep_0\ \int_{D^2_r}|\nabla f^r|^2\ dx\,dy
\end{array}
\ee
for some universal $C>0$. Thus for $\ep_0$ chosen small enough we obtain the existence of $C>0$ independent of $r$ and the data of the lemma such that
\be
\label{A.9d}
\int_{D^2_r}|\nabla f^r|^2\ dx\,dy\le \frac{3}{16\pi}(1+2C\ep_0)\ \int_{D^2_r}|\nabla e_1^r|^2\ dx\,dy\ \int_{D^2_r}|\nabla e_2^r|^2\ dx\,dy\
\ee
Once this estimate is established, the rest of the arguments of F. H\'elein carries over and we obtain lemma~\ref{lm-A.2}. \hfill $\Box$

\begin{Lma}
\label{lm-A.3}
There exists $\ep_0>0$ and $0<\al<1$ such that for any $\vec{\xi}\in W^{2,2}\cap W^{1,\infty}(D^2,{\R}^m)$ satisfying
\[
|Dis(\vec{\xi})|<\ep_0
\]
and
\[
\int_{D^2}|\nabla\vec{n}_{\vec{\xi}}|^2_g\ dvol_g<\frac{4\pi}{3}\quad,
\]
where $g:=\vec{\xi}^\ast g_{{\R}^m}$ there exists $\zeta\in W^{1,\infty}_{loc}\cap W^{2,2}_{loc}(D^2,D^2)$ such that
$\vec{\xi}\circ\zeta$ is conformal, 
\be
\label{A.10}
\|\zeta\|_{C^{0,\al}(D^2)}+\|\zeta^{-1}\|_{C^{0,\al}(D^2)}\le\ C\ \lf[\exp(C\ \|\log|\nabla\vec{\xi}|\|_{L^\infty(D^2)})\rg]\quad,
\ee
where $C>0$ is independent of $\vec{\xi}$ and
\be
\label{A.11}
\begin{array}{l}
\forall\ 1>r>0\quad\quad\exists\ C_r>0\quad\mbox{ s. t. }\\[5mm]
\|\log|\nabla(\vec{\xi}\circ\zeta)|\|_{L^\infty(D^2_r)}+\|\vec{\xi}\circ\zeta\|_{W^{2,2}(D^2_r)}\le C_r \lf[1+\exp(\|\log|\nabla\vec{\xi}|\|_\infty)\rg]
\end{array}
\ee
where $C_r>0$ only depends on $r$ and  not on $\vec{\xi}$. \hfill $\Box$
\end{Lma}
{\bf Proof of lemma~\ref{lm-A.3}.}
Let $(\vec{e}_1,\vec{e}_2)$ be the orthonormal 2-framing given by lemma~\ref{lm-A.2} such that $\vec{n}_{\vec{\xi}}=\vec{e}_1\wedge \vec{e}_2$
and let $f$ such that $f\equiv 0$ on $\p D^2$ and $df=\ast_g(\vec{e}_1,\ d\vec{e}_2)$. Hence it solves
\[
\lf\{ 
\begin{array}{l}
\ds\Delta_{g}f=(\nabla^\perp\vec{e_1},\nabla\vec{e}_2)\quad\quad\mbox{ on }D^2\\[5mm]
f=0\quad\quad\quad\mbox{ on }\p D^2\quad.
\end{array}
\rg.
\]
writing as in the proof of lemma~\ref{lm-A.2} 
\be
\label{A.12}
\lf\{
\begin{array}{l}
\ds\Delta_{g_0} f=\p_x\vec{e}_1\cdot\p_y\vec{e}_2-\p_x\vec{e}_2\cdot\p_y\vec{e}_1+\ast_{g_0}d( (\ast_{g_0}-\ast_g)\ df^r)\quad\quad\mbox{ in }D^2\\[5mm]
\ds f=0\quad\quad\quad\mbox{ on }\quad\p D^2
\end{array}
\rg.
\ee
Using integrability compensation result Theorem 3.4.1 we get the a-priori estimate
\be
\label{A.13}
\|\nabla f\|_{L^{2,1}(D^2)}\le C\ \|\nabla \vec{e}_1\|_{L^2(D^2)}\ \|\nabla \vec{e}_2\|_{L^2(D^2)}+C\ \ep_0\ \|\nabla f\|_{L^{2,1}(D^2)}
\ee
where we have used the fact that $\forall \al$ one-form on $D^2$ $|(\ast_g-\ast_{g_0})\,\al|\le\ep_0\ |\al|$. Hence by density of smooth
maps, for $\ep_0$ small enough, we convert the a-priori estimate (\ref{A.13}) into an estimate and by Lorentz-Sobolev embedding
one has
\be
\label{A.14}
\|f\|_{L^\infty(D^2)}\le C\ \|\nabla f\|_{L^{2,1}(D^2)}\le C\ \int_{D^2}|\nabla\vec{n}_{\vec{\xi}}|^2\ dvol_g
\ee
Let $\ep_i:=d\vec{\xi}^{-1}\vec{e}_i$ and let $\ep_i^\ast$ be the dual unit frame for the $g=\vec{\xi}^\ast g_{{\R}^m}$ metric of $(\ep_1,\ep_2)$.
We have that
\be
\label{A.15} \lf\{
\begin{array}{l}
\ds df\wedge \ep_1^\ast=(\ast_g df)\wedge(\ast_g \ep_1^\ast)=-(\vec{e}_1,d\vec{e}_2)\wedge\ep_2^\ast\quad,\\[5mm]
\ds df\wedge \ep_2^\ast=(\ast_g df)\wedge(\ast_g \ep_2^\ast)=(\vec{e}_1,d\vec{e}_2)\wedge\ep_1^\ast\quad,
\end{array}
\rg.
\ee
Moreover Cartan formula gives
\be
\label{A.16}
\begin{array}{rl}
\ds d\ep_i^\ast(\vec{e}_1,\vec{e}_2)&\ds=d(\ep_i^\ast(\vec{e}_2))\cdot \vec{e}_2-d(\ep_i^\ast(\vec{e}_1))\cdot \vec{e}_2-\ep_i^\ast([\vec{e}_1,\vec{e}_2])\\[5mm]
 &=-\ep_i^\ast([\vec{e}_1,\vec{e}_2])=-\ep_i^\ast(D_{\vec{e}_1}\vec{e}_2-D_{\vec{e}_2}\vec{e}_1)
\end{array}
\ee
where $D$ is the Levi-Civita connection associated to the induced metric $g$. Using the immersion $\vec{\xi}$ we have that $D_X\vec{e}_i=P_{\vec{\xi}}(d\vec{e}_i\cdot X)$
where $P_{\vec{\xi}}$ is the orthonormal projection in ${\R}^m$ onto $T\vec{\xi}(D^2)$. Hence $D_X\vec{e}_i=(\vec{e}_1,d\vec{e}_i\cdot X)\ \vec{e}_1+(\vec{e}_2,d\vec{e}_i\cdot X)\ \vec{e}_2$. Combining this later fact with (\ref{A.16}) gives
\be
\label{A.17}
\lf\{
\begin{array}{l}
\ds d\ep_1^\ast=-(\vec{e}_1,d\vec{e}_2)\wedge\ep_2^\ast\quad,\\[5mm]
\ds d\ep_2^\ast=-\ep_1^\ast\wedge(\vec{e}_1,d\vec{e}_2)\quad.
\end{array}
\rg.
\ee
Combining (\ref{A.15}) and (\ref{A.17}) gives, 
\be
\label{A.18}
d(e^{-f}\ep_1^\ast)=0\quad\quad\mbox{ and }\quad\quad d(e^{-f}\ep_2^\ast)=0\quad.
\ee
Hence there exists $\sigma:=(\sigma_1,\sigma_2)\in W^{2,2}\cap W^{1,\infty}(D^2,{\R}^2)$ such that
\be
\label{A.19}
d\si_i=e^{-f}\ \ep^\ast_i\quad\quad\mbox{ and }\quad\quad \si_i(0)=0\quad.
\ee
$\nabla \si$ has maximal rank equal to 2 at every point, therefore it realizes a lipschitz diffeomorphism from $D^2$ into $\Omega:=\si(D^2)$. 
Let $\frac{\p}{\p \si_i}$ be the dual basis to $d\si_i$. Since $d\vec{\xi}=\sum_{i=1}^2d\vec{\xi}\cdot\ep_i\ \ep_i^\ast$,  one has $d\vec{\xi}\cdot\frac{\p}{\p \si_i}=e^f\ \vec{e}_i$
Hence $\vec{\xi}\circ\si^{-1}$ is a conformal immersion. For every one-form $\al$ on $D^2$ one has
\be
\label{A.20}
\quad\quad (1-\ep_0^2)\ |\al|^2_{g_0}\le {(det g)}\ |\al|^2_g\le\ (1+\ep_0^2)\ |\al|^2_{g_0}\quad,
\ee
Hence, since $|d\si|_g^2=e^{2f}$, since $det g=|\p_x\vec{\xi}\times\p_y\vec{\xi}|^2$ and since $1/2- Dis(\vec{\xi})^2/2\le |\nabla\vec{\xi}|^{-2}\  |\p_x\vec{\xi}\times\p_y\vec{\xi}|^2\le1/2+Dis(\vec{\xi})^2/2$, we deduce from (\ref{A.20}) and (\ref{A.14}) the following estimate
\be
\label{A.21}
\|\log|\nabla\si|_{g_0}\|_{L^\infty(D^2)}\le C\ \|\log|\nabla\vec{\xi}|\|_{L^\infty(D^2)}+C\ \int_{D^2}|\nabla\vec{n}_{\vec{\xi}}|^2\ dvol_g\quad.
\ee
Moreover from the above remark\footnote{We use that $|d\si_1\wedge d\si_2|_g=|d\si|^2_g/2$ and (\ref{A.20}).} we have that $|Dis({\si})|<2\ \ep_0$.
Hence we deduce
\be
\label{A.22}
\lf\{
\begin{array}{l}
\ds|Dis(\si^{-1})|<2\ep_0\quad\quad\mbox{ and }\\[5mm]
\ds\|\log|\nabla\si^{-1}|_{g_0}\|_{L^\infty(\Omega)}\le C\ \|\log|\nabla\vec{\xi}|\|_{L^\infty(D^2)}+C\ \int_{D^2}|\nabla\vec{n}_{\vec{\xi}}|^2\ dvol_g\quad.
\end{array}
\rg.
\ee
Let $h$ be the solution to the Riemann Mapping Theorem for $\Omega$ :
\[
\begin{array}{l}
h\ :\ \Omega\longrightarrow D^2\quad\mbox{ is holomorphic }\\[5mm]
h(0)=0\quad\quad\mbox{ and }\quad h'(0)\in{\R}\quad\mbox{ with }\quad h'(0)>0\quad.
\end{array}
\]
The Riemann Mapping Theorem asserts that $h$ is bi-holomorphic and we will denote $k$ it's holomorphic inverse from $D^2$ into $\Omega$. 
Finally let $\zeta:=\si^{-1}\circ k$. (We shall often see $\zeta$
as a ${\C}-$valued map).
$\zeta$ satisfies
\be
\label{A.22a}
\p_{\overline{z}}\zeta=\nu(z)\ \overline{\p_z\zeta}\quad,
\ee
where $\nu\circ h:=\p_{\overline{z}}\si^{-1}/\overline{\p_{z}\si^{-1}}=H(\nabla\si^{-1})/[|\nabla\si^{-1}|^2+2\ det\nabla\si^{-1}]$. Since $det\nabla\si^{-1}\ge 0$, we obtain
\be
\label{A.23}
\|\nu\|_{L^{\infty}(D^2)}\le 2\ep_0\quad .
\ee
Let $\delta:=dist(0,\p\Omega)$. Integrating $d\si^{-1}$ on a segment $S$ connecting $0$ to one of it's nearest  point $P$  on $\p \Omega$ gives $1=|\int_S\ d\si^{-1}|\le \delta\ \|\nabla \si^{-1}\|_\infty$
and integrating now $d\si$ on a ray $R$ issued from zero and connecting $\si^{-1}(P)$ gives $\delta=|P-0|=|\int_{R}d\si|\le \|\nabla\si\|_\infty$. Hence we have
\be
\label{A.24}
\frac{1}{\|\nabla \si^{-1}\|_\infty}\le\delta\le \|\nabla\si\|_\infty\quad.
\ee
Combining (\ref{A.21}), (\ref{A.22}) and (\ref{A.24})  gives
\be
\label{A.25}
|\log\delta|\le C\ \|\log|\nabla\vec{\xi}|\|_{L^\infty(D^2)}+C\ \int_{D^2}|\nabla\vec{n}_{\vec{\xi}}|^2\ dvol_g
\ee
Since $h$ is holomorphic and $h(\Omega)=D^2$ we have
\be
\label{A.28}
\|h'\|_{L^{\infty}(B_{\delta/2}(0))}\le\ \frac{C_1}{\delta}\ \|h\|_\infty\le\  \frac{C_1}{\delta}\quad.
\ee
This implies that $h(B^2_r(0))\subset B^2_{C_1r/\delta}(0)$ and hence for instance $k(\p B_{1/4}(0))\subset \Omega\setminus B_{\delta/4 C_1}(0)$
Hence
\be
\label{A.28a}
\begin{array}{rl}
\ds\|\log|\xi|\|_{L^\infty(\p B_{1/4}(0))}&\ds\le |\log\lf[4 C_1\|\nabla \si\|_\infty\rg]|+|\log\delta|\\[5mm]
 &\ds\le C\ \lf[1+\|\log|\nabla\vec{\xi}|\|_{L^\infty(D^2)}\rg]+C\ \int_{D^2}|\nabla\vec{n}_{\vec{\xi}}|^2\ dvol_g
 \end{array}
 \ee
 We have
\be
\label{A.29}
\|k'\|_{L^\infty(B_{1/2}(0))}\le 2 C_1\ \|k\|_{\infty}=2\ C_1\|\nabla\si\|_\infty\quad,
\ee
This implies that 
\be
\label{A.30}
\begin{array}{rl}
\ds\|\nabla\zeta\|_{L^\infty(B_{1/4}(0))}&\ds\le 2\ C_1\ \|\nabla\si^{-1}\|_\infty\ \|\nabla\si\|_\infty\\[5mm]
 &\ds\le 2\ C_1\ \exp\lf[C\ \|\log|\nabla\vec{\xi}|\|_{L^\infty(D^2)}+C\ \int_{D^2}|\nabla\vec{n}_{\vec{\xi}}|^2\rg]\quad,
 \end{array}
\ee
In $B_1(0)\setminus B_{1/4}(0)$ we write $\zeta i=e^{\la+i\mu}$ where $\la$ is a real-valued function and where $\mu$ takes value into ${\R}/2\pi{\Z}$. 
Using this notation (\ref{A.22a}) becomes
\[
\frac{[1-\nu\ e^{-2i\mu}]}{[1+\nu\ e^{-2i\mu}]}\ \p_{\ov{z}}\la=-i\ \p_{\ov{z}}\mu\quad,
\]
which implies
\be
\label{A.31}
\Re\lf[\frac{\p}{\p z}\lf[(1+\beta)\frac{\p\la}{\p \ov{z}}\rg]\rg]=0\quad,
\ee
where $\beta:=\beta_1+i\beta_2=-2\nu\ e^{-2i\mu}\ [1+\nu\ e^{2i\mu}]^{-1}$. Hence $\la$ satisfies
\be
\label{A.32}
\lf\{
\begin{array}{l}
\ds\sum_{i,j=1}^2\p_{x_i}[a_{ij}\ \p_{x_j}\la]=0\quad\quad\mbox{ in }B^2_1(0)\setminus B^2_{1/4}(0)\\[5mm]
\ds\la=0\quad\quad\quad\mbox{ on }\p B^2_1(0)\\[5mm]
\ds\la=\log|\zeta|\quad\quad\mbox{ on }\p B^2_{1/4}(0)
\end{array}
\rg.
\ee
where $a_{11}=a_{22}=1+\beta_1$ and $a_{12}=-a_{21}=-\beta_2$. From (\ref{A.23}) we have that $|\beta|<4\ep_0$.
Hence for $\ep_0>0$ small enough De Giorgi-Nash result (see for instance \cite{TA3}) gives the existence of $0<\al<1$ such that
\[
\|\la\|_{C^{0,\al}(B^2_1\setminus B^2_{1/4})}\le C\ \| \log|\zeta|\|_{W^{1,\infty}(\p B_{1/4}(0))}\quad .
\]
Using (\ref{A.30}) we then obtain
\be
\label{A.33}
\||\zeta|\|_{C^{0,\al}(B^2_1\setminus B^2_{1/4})}\le 2\ C_1\ \exp\lf[C\ \|\log|\nabla\vec{\xi}|\|_{L^\infty(D^2)}+C\ \int_{D^2}|\nabla\vec{n}_{\vec{\xi}}|^2\rg]\quad.
\ee
Since
\be
\label{A.34}
\Re\lf[\frac{\p}{\p z}\lf[(1+\beta)^{-1}\frac{\p\mu}{\p \ov{z}}\rg]\rg]=0\quad,
\ee
we get a similar control to (\ref{A.33}) for $arg\Psi$ on $B^2_1\setminus B^2_{1/4}$. Hence combining these last estimates together with
(\ref{A.30}) again, we finally obtain
\be
\label{A.35}
\|\zeta\|_{C^{0,\al}(D^2)}\le 2\ C_1\ \exp\lf[C\ \|\log|\nabla\vec{\xi}|\|_{L^\infty(D^2)}+C\ \int_{D^2}|\nabla\vec{n}_{\vec{\xi}}|^2\rg]\quad,
\ee
which finishes the proof of lemma~\ref{lm-A.3}.\hfill $\Box$

\begin{Lma}
\label{lm-A.4}
Let $\vec{\xi}$ be a conformal immersion of $D^2\setminus\{0\}$ into ${\R}^m$ in $W^{2,2}_{loc}(D^2\setminus\{0\},{\R}^m)$ and such that
$\log|\nabla\vec{\xi}|\in L^\infty_{loc}(D^2\setminus\{0\})$. Assume $\vec{\xi}$ extends to a map in $W^{1,2}(D^2)$ and that the corresponding
Gauss map $\vec{n}_{\vec{\xi}}$ also extends to a map in $W^{1,2}(D^2,Gr_{m-2}({\R}^m))$. Then $\vec{\xi}$ realizes a lipshitz conformal immersion 
of the whole disc $D^2$ and there exits a positive integer $m$ and a constant $C$ such that
\be
\label{A.36}
(C-o(1))\ |z|^{m-1}\le\lf|\frac{\p\vec{\xi}}{\p z}\rg|\le (C+o(1))\ |z|^{m-1}\quad.
\ee
 \hfill $\Box$
\end{Lma}
{\bf Proof of Lemma~\ref{lm-A.4}}. We can always localize in order to ensure that
\[
\int_{D^2}|\nabla\vec{n}_{\vec{\xi}}|^2\ dx\,dy<\frac{8\pi}{3}\quad.
\]
Using lemma 5.1.4 of \cite{Hel} we deduce the existence of a framing $\vec{e}:=(\vec{e}_1,\vec{e}_2)$ which is in $W^{1,2}(D^2,S^{m-1}\times S^{m-1})$ such that
\be
\label{A.37}
\vec{e}_1\cdot \vec{e}_2=0\quad\quad,\quad\quad n_{\vec{\xi}}=\vec{e}_1\wedge \vec{e}_2 \quad,
\ee
\be
\label{A.38}
\int_{D^2}\lf[|\nabla \vec{e}_1|^2+|\nabla \vec{e}_2|^2\rg]\ dx\,dy\le 2\int_{D^2}|\nabla\vec{n}_{\vec{\xi}}|^2\ dx\,dy
\ee
and
\be
\label{A.39}
\lf\{
\begin{array}{l}
\ds div(\vec{e}_1,\nabla \vec{e}_2))=0\quad\quad\mbox{ in }D^2\quad,\\[5mm]
\ds\lf(\vec{e}_1,\frac{\p \vec{e}_2}{\p \nu}\rg)=0\quad\quad\mbox{ on }\p D^2\quad .
\end{array}
\rg.
\ee

Similarly as in the proof of lemma~\ref{lm-A.3}, we introduce $\ep_i:=d\vec{\xi}^{-1}\vec{e}_i$ and $\ep_i^\ast$ to be the dual
framing. Denoting $|\p_x\vec{\xi}|^2=|\p_y\vec{\xi}|^2=e^{2\la}$ we have that the metric 
$g_\infty:=\vec{\xi}^\ast\,g_{{\R}^m}$ is given by $g=e^{2\la}\ \lf[dx^2+dy^2\rg]$. Hence with respect to the flat metric $g_0:=\lf[dx^2+dy^2\rg]$ one has 
\[
|\ep_i|^2_{g_0}= g_0({\ep_i},{\ep_i})=e^{-2\la}\ g_\infty(\ep_i,\ep_i)=e^{-2\la}\quad.
\]
and since $<\ep_i,\ep_j^\ast>=\delta_{ij}$ we have that $|\ep_i^\ast|^2_{g_0}= e^{2\la}$. Thus we deduce that the 1-forms $\ep_i^\ast$ are
in $L^2(D^2)$. Since $\vec{\xi}$ is in $W^{1,\infty}\cap W^{2,2}_{loc}(D^2\setminus\{0\},{\R}^m)$ and 
$\log|\nabla\vec{\xi}|\in L^\infty_{loc}(D^2\setminus\{0\})$ we have that the framing given by $\vec{f}_i:=e^{-\la}\ \p_{x_i}\vec{\xi}$ is in 
$L^\infty_{loc}\cap W^{1,2}_{loc}(D^2\setminus\{0\},{\R}^m)$. Since
$\vec{\xi}$ is conformal the unit framing $(\vec{f}_1,\vec{f}_2)$ is Coulomb\footnote{This follows from a straightforward computation 
presented in \cite{Hel} chapter 5.} :
\[
div(\vec{f}_1,\nabla\vec{f}_2)=0\quad\quad\mbox{ in }D^2\setminus\{0\}\quad.
\]
Denoting $e^{i\theta}$ the rotation which passes\footnote{ $e^{i\theta}\,(\vec{f}_1+i\vec{f}_2)=(\vec{e}_1+i\vec{e}_2)$} from $(\vec{f}_1,\vec{f}_2)$ to $(\vec{e}_1,\vec{e}_2)$. The Coulomb condition satisfied by the two framings implies that $d\theta:=(ie^{i\theta},d(e^{i\theta}))$ is an harmonic 1-form on $D^2\setminus\{0\}$ and hence analytic on this domain. This implies that 
\[
\ep_i^\ast\in L^\infty_{loc}\cap W^{1,2}_{loc}(D^2\setminus\{0\})\quad\quad.
\]
 Like again in lemma~\ref{lm-A.3} we introduce\footnote{By virtue of Wente's theorem (see theorem 3.1.2 of \cite{Hel}).} $f\in C^{0}\cap W^{1,2}(D^2)$ to be the solution of
\[
\lf\{ 
\begin{array}{l}
\ds\Delta f=(\nabla^\perp\vec{e_1},\nabla\vec{e}_2)\quad\quad\mbox{ on }D^2\\[5mm]
f=0\quad\quad\quad\mbox{ on }\p D^2\quad.
\end{array}
\rg.
\]
As in lemma~\ref{lm-A.3}, the computations give in $D^2\setminus\{0\}$
\[
\forall i=1,2\quad\quad\quad d[e^{-f}\ep_i^\ast]=0\quad\quad\mbox{ a.e. in }D^2\setminus\{0\}\quad.
\]
By schwartz lemma the distribution $d[e^{-f}\ep_i^\ast]$ is a finite linear combination of successive derivatives of the Dirac Mass
at the origin but since $e^{-f}\ep_i^\ast\in L^2(D^2)$, this linear combination can only be 0. Hence we have
\[
\forall i=1,2\quad\quad\quad d[e^{-f}\ep_i^\ast]=0\quad\quad\mbox{ in }{\mathcal D}'(D^2)\quad.
\]
Hence, by Poincar\'e Lemma, there exists $(\sigma_1,\sigma_2)\in W^{1,2}(D^2,{\R}^2)$ such that $d\sigma_i=e^{-f}\ep_i^\ast$. The dual basis
$(\p/\p\sigma_1,\p/\p\sigma_2)= e^f(\ep_1,\ep_2)$ is positive, orthogonal on $D^2\setminus\{0\}$ and integrable by nature.
 hence $\sigma=\sigma_1+i\sigma_2$ is an holomorphic function on $D^2\setminus\{0\}$ which extends to a $W^{1,2}-$map
 on $D^2$. The classical point removability theorem for holomorphic map implies that $\sigma$ extends to an holomorphic function on $D^2$.
 We can choose it in such a way that $\sigma(0)=0$.
 The holomorphicity of $\sigma$ implies in particular that $|d\sigma|_{g_0}=e^{\la-f}$ is uniformly bounded and, since $f\in L^\infty(D^2)$, we deduce that $\la$ is bounded
 from above on $D^2$. This later fact implies that $\vec{\xi}$ extends to a Lipshitz map on $D^2$. Though $|d\sigma|_{g_0}=e^{\la-f}$
 has no zero on $D^2\setminus\{0\}$, $\sigma'$ might have a zero at the origin : there exists an holomorphic function $h(z)$ on $D^2$
 satisfying $h(0)=0$, a complex number $c_0$ and an integer $m$ such that
 \be
 \label{A.40}
 \sigma(z)= c_0\  z^m\ (1+h(z))\quad\quad.
 \ee
We have that locally 
$$
\frac{\p \vec{\xi}}{\p \sigma}=\p_{\sigma_1}\vec{\xi}-i\p_{\sigma_2}\vec{\xi}=d\vec{\xi}\ e^{f}\ep_1^\ast-id\vec{\xi}\ e^{f}\ep_2^\ast =e^f[\vec{e}_1-i\vec{e}_2]
$$
 Hence, since $f$ is continuous, we have that
 \be
 \label{A.41}
 \lf|\frac{\p \vec{\xi}}{\p \sigma}\rg|= \sqrt{2}\ e^{f(0)}\ (1+o(1))\quad.
 \ee
 Combining (\ref{A.40}) and (\ref{A.41}) gives 
 \be
 \label{A.42}
 \lf|\frac{\p \vec{\xi}}{\p z}\rg|=  \lf|\frac{\p \vec{\xi}}{\p \sigma}\rg|\ \lf|\frac{\p \sigma}{\p z}\rg|=c_0\ m\ \sqrt{2}\ e^{f(0)}\ |z|^{m-1}\ (1+o(1))\quad.
 \ee
 This last identity implies the lemma~\ref{lm-A.4}.\hfill $\Box$

\end{document}